\documentclass[11pt]{article}

\usepackage{amsmath,amsthm} 
\usepackage{amssymb,mathrsfs} 
\usepackage{a4wide} 
\usepackage{graphicx}
\usepackage{color} 
\usepackage{enumerate}
\usepackage{todonotes}
\usepackage[normalem]{ulem}
\usepackage{mathtools}
\usepackage{booktabs}
\usepackage{pifont}
\usepackage{subcaption}
\usepackage{hyperref}

\usepackage[capitalise,nameinlink]{cleveref}

\usepackage{diagbox}
\usepackage{bbm}

\newcommand{\dps}{\displaystyle } 
\newcommand{\rme}{\mathrm{e}}
 
\newcommand{\cL}{\mathcal{L}}
\newcommand{\ocL}{\overline{\mathcal{L}}}
\newcommand{\cX}{\mathcal{X}}
\newcommand{\cE}{\mathcal{E}}
\newcommand{\cD}{\mathcal{D}}
\newcommand{\cR}{\mathcal{R}}
\newcommand{\wcL}{\widetilde{\mathcal{L}}}
\newcommand{\cLham}{\mathcal{L}_{\rm ham}}

\newcommand{\R}{\mathbb{R}}

\renewcommand{\d}{\mathrm{d}}
\renewcommand{\leq}{\leqslant}
\renewcommand{\geq}{\geqslant}

\newcommand\dd{\mathrm{d}}

\newcommand{\rmL}{\mathrm{L}}
\newcommand{\rmR}{\mathrm{R}}
\newcommand{\cA}{\mathcal{A}}

\newcommand{\cS}{\mathcal{S}}

\newcommand{\calD}{\cD}
\newcommand{\calL}{\cL}
\newcommand{\calE}{\cE}
\newcommand{\calX}{\cX}
\newcommand{\bbE}{\mathbb{E}}
\newcommand{\bbV}{\mathbb{V}}
\newcommand{\rmd}{\mathrm{d}}

\newcommand{\E}{\mathbb{E}}
\newcommand{\1}{\mathbbm{1}}

\newcommand{\expect}{\E}

\newcommand{\Niter}{{N_\mathrm{iter}}}
\newcommand{\dt}{{\Delta t}}

\begin{document}


\title{Mathematical analysis and numerical methods for the computation of transport coefficients in molecular dynamics}
\author{No\'e Blassel$^{1}$, Louis Carillo$^{2,3}$, Shiva Darshan$^{4}$, Rapha\"el Gastaldello$^{2,3}$, \\
  Alessandra Iacobucci$^{5}$, Elisa Marini$^5$, Régis Santet$^{6,7,8}$, Xiaocheng Shang$^9$, \\
  Gabriel Stoltz$^{2,3}$, Urbain Vaes$^{3,2}$ \\
  \small 1: MatMat group, Institute of Mathematics, EPFL, Lausanne, Switzerland \\
  \small 2: CERMICS, CNRS, ENPC, Institut Polytechnique de Paris, Marne-la-Vallée, France \\
  \small 3: MATHERIALS project-team, Inria Paris, France \\
  \small 4: ESSEC Business School, Cergy, France \\
  \small 5: CEREMADE, Université Paris Dauphine-PSL \\
  \small 6: CEA, DES, IRESNE, DEC, Cadarache, F-13108 Saint-Paul-Lez-Durance, France \\
  \small 7: CEA, DAM, DIF, F-91297 Arpajon, France \\
  \small 8: Université Paris-Saclay, CEA, Laboratoire Matière en Conditions Extrêmes, 91680 Bruyères-le-Châtel, France \\
  \small 9: School of Mathematics, University of Birmingham, Edgbaston, Birmingham, B15 2TT, United Kingdom
}
 
\maketitle

\begin{abstract}
  We review various numerical approaches to compute transport coefficients in molecular dynamics. These approaches can be broadly classified into three groups: (i) nonequilibrium methods based on applying an external driving field to the system, measuring the average response in the system, and evaluating the related linear response coefficient; (ii) approaches reformulating the transport coefficient of interest through a time correlation function for the equilibrium dynamics (the most popular instances being Green--Kubo and Einstein formulas); (iii) transient techniques, where the transport coefficient can be computed by monitoring the return to the steady state of a dynamics perturbed off its stationary distribution. For all three classes of methods, we provide elements of numerical analysis, allowing to estimate or at least quantify the level of numerical errors in the estimator of the transport coefficient; and also briefly present recent attempts to more efficiently compute transport coefficients with variance reduction approaches such as control variates, importance sampling and coupling methods. The computation of transport coefficients remains nonetheless challenging and will continue requiring research efforts in the foreseeable future.
\end{abstract}

\section{Why this review?}

Transport coefficients are key quantities to understand how physical properties of a dynamical system behave under perturbations, with applications in various fields ranging from biophysics to materials science. For instance, the mobility, which is proportional to the diffusion coefficient, plays a key role in the design of ionic batteries; while the thermal conductivity of atom chains is essential to understanding the necessary and/or sufficient conditions for the validity of Fourier's law in theoretical physics. However, obtaining exact values of these coefficients, either from direct experimental measurements or by simulations, has always been a challenge. With the rapid increase of computational power, molecular dynamics (MD) has proven to be highly effective, by capitalizing on scalable parallelized algorithms~\cite{Frenkel2001,Allen2017,Tuckerman2023}. A historical perspective on the numerical methods devised to estimate transport coefficients can be found in~\cite{Ciccotti_history_MD}.

At the macroscopic scale, a transport coefficient quantifies the relationship between the average flux (or response) of some transported physical quantity in the steady state of a system driven out of equilibrium by an external forcing (or perturbation), and the magnitude of that forcing. Such a system is known as a nonequilibrium system, and standard examples of transport coefficients, which we present in detail below, include the mobility and the shear viscosity of a fluid, or the thermal conductivity in atom chains and lattices.

The computation of transport coefficients is a crucial step in parameterizing continuum models of fluids and materials at a macroscopic scale, such as the Navier--Stokes or heat equation. However, their estimation from trajectory data can be very costly with naive methods, prompting the need for innovative numerical strategies. This topic is still the object of active research, both in the MD and mathematical communities, see~\cite{stoltz2024} and references therein for another recent overview.

Broadly speaking, computational methods to measure transport coefficients from MD simulations fall in one of two categories. The first relies on the analysis of time-dependent signals in equilibrium systems, following the pioneering theoretical work of Kubo~\cite{kubo1957,kubo1957b}. The second is the family of nonequilibrium molecular dynamics (NEMD), which relies on direct measurement of various response functions in a system subjected to perturbations of the equilibrium dynamics. We refer to~\cite{ciccotti1975,ciccotti1976} for early accounts of NEMD experiments, as well as~\cite{ciccotti1979,hoover1986,hoover1993,ciccotti2005,ciccotti2016} for reviews of the NEMD approach at various stages of historical development, and the books~\cite{evans_morriss_2007,todd_daivis_2017}. Early measurements of transport coefficients from equilibrium MD trajectories can be found in~\cite{alder1970,verlet1970,verlet1973}. 

The equilibrium approach is subject to large statistical fluctuations in measuring the transport coefficient, due to the difficulty of accurately evaluating long time correlations; therefore, the signal-to-noise ratio is particularly unfavorable at long times, where there may be a significant (but unwanted) contribution of random fluctuations to the integral defining the transport coefficient. On the other hand, to numerically detect the induced system response in the NEMD approach, one typically applies perturbations that are orders of magnitude larger than those used in real experiments, thereby enhancing the signal-to-noise ratio of the measured response. However, perturbations must remain sufficiently small to avoid nonlinear effects in the response, and so a careful balance must be maintained. The limitations of standard numerical techniques call for dedicated alternative approaches, particularly those aimed at variance reduction. 

While there may be lots of advice and insights on computing/debugging static properties in MD textbooks, much less material is available for transport coefficients. Particularly, in contrast with methods for equilibrium sampling, the error quantification of transport coefficients is not well developed. This is a major motivation for this review, which aims at providing elements of numerical analysis.


\paragraph{Outline of the review.}
We first provide in Section~\ref{sec:presentation} a brief review on the computation of equilibrium properties in MD, describing paradigmatic systems which are used throughout this work to illustrate the methods and concepts at hand. We then turn to nonequilibrium dynamics and the definition of transport coefficients in Section~\ref{sec:transport_coefficients}. We next successively present the most important classes of methods to compute transport coefficients, discussing in all cases standard approaches as well as recent methods aiming at improving the numerical performance of the estimators at hand. Nonequilibrium methods are considered in Section~\ref{sec:NEMD}, fluctuation formulas for equilibrium dynamics are discussed in Section~\ref{sec:GK_and_co}, and transient methods are presented in Section~\ref{sec:transient}. We conclude the review with some perspectives in Section~\ref{sec:extensions}.

\section{A brief presentation of molecular dynamics}
\label{sec:presentation}

We provide in this section a brief presentation of sampling methods in MD, with the aim of computing static average properties. We start by reviewing the most important elements of statistical physics in Section~\ref{sec:stat_phys}, and then present in Section~\ref{sec:example_systems} the three paradigmatic systems we will consider to illustrate the various concepts we discuss in this review. We next discuss in Section~\ref{sec:sampling_methods} methods to approximate in practice ensemble averages given by statistical physics, relying on ergodic averages of stochastic dynamics, such as underdamped and overdamped Langevin dynamics. Error estimates quantifying the quality of the approximation are studied in Section~\ref{sec:equilibrium_error_estimates}. We conclude in Section~\ref{sec:equilibrium_variance_reduction} by presenting various ways to reduce the most important source of error, namely the statistical error.

\subsection{Some elements of statistical physics}
\label{sec:stat_phys}

MD has to be understood through the lens of statistical physics. Consider a system composed of~$N$ particles living in~$D$ spatial dimensions, so that the total dimension of the system is~$d = ND$. The physical state of the particles can be described by an element~$x$ of the so-called configuration space~$\cX$. A state~$x$ can be either the positions of the particles~$q \in \cD$, or their positions and momenta~$(q, p)\in \cE = \cD \times \R^d$. A typical choice for the position space is the torus~$\cD = (L\mathbb{T})^d$ (with~$L>0$ being the box size), which corresponds to using periodic boundary conditions (see however~\cref{sec:atom_chains} for another choice of position space). The other element to describe the microscopic configuration of the system is the energy, which can be either the potential energy~$V(q)$ for systems described by their positions only, or, for systems described by positions and momenta, the Hamiltonian
\begin{equation}\label{eq:hamiltonian}
    H(q,p) = \frac{1}{2}p^\top M^{-1} p + V(q),
\end{equation}
with $M = \mathrm{Diag}(m_1, \, \dots, \, m_N)$ the mass matrix. The Hamiltonian is here considered in its simplest and most common form, namely separable (sum of kinetic and potential energies) with a quadratic kinetic energy. The potential energy function~$V$ incorporates all the physics of the model. It should ideally be obtained by resorting to {\em ab-initio} computations~\cite{CDK+2003}, or machine learned potentials fitted on \emph{ab-initio} data~\cite{JMA+2025}. In practice, for large systems as those considered in statistical mechanics, in particular for biological applications, one still often relies on parameterized empirical formulas for the potential energy function, whose physical parameters are fitted to best reproduce properties of interest. 

Let us next turn to the macroscopic description of physical systems, which relies on a probability measure~$\pi$ on the configuration space. This probability measure is the least biased distribution consistent with the invariants of the system~\cite{Balian}. Macroscopic properties are obtained via the average of an observable~$\varphi \in L^1(\pi)$ with respect to the probability measure under consideration: 
\begin{equation}\label{intro:averages}
  \mathbb{E}_\pi(\varphi) = \int_{\cX} \varphi(x) \, \pi(\rmd x).  
\end{equation}
For the sake of simplicity, we introduce only the canonical ensemble (fixed number of particles, volume and temperature), as it is the playground for a large number of problems in MD (if not most of them); see for instance~\cite{Tuckerman2023} for a presentation of other ensembles used in MD simulations. The probability measure describing the canonical ensemble is the Boltzmann--Gibbs measure
\begin{equation}\label{eq:canonical measure}
  \pi(\d x) = \left\{ \begin{aligned}
    \mu(\d q \, \d p) & = Z_\mu^{-1} \rme^{-\beta H(q, p)} \, \d q \,  \d p \qquad \qquad 
    && Z_\mu = \int_{\cE}\rme^{-\beta H},\\
    \nu(\d q) & = Z_\nu^{-1} \rme^{-\beta V(q)} \, \d q \qquad \qquad  
    && Z_\nu= \int_{\cD}\rme^{-\beta V},
  \end{aligned} \right.
\end{equation}
where $\beta = 1 / (k_\mathrm{B} T)$ is proportional to the inverse temperature, with $k_\mathrm{B} = 1.38 \times 10^{-23} \mathrm{J} \cdot \mathrm{K}^{-1}$ the Boltzmann constant; and the partition functions~$Z_\mu$ and~$Z_\nu$ are normalization constants (assumed to be finite, which is the case when~$\rme^{-\beta V} \in L^1(\cD)$). Note that~$\nu$ is the marginal distribution in positions of~$\mu$, and that it is easy to sample the marginal distribution in the momenta of~$\mu$ since this amounts to sampling a Gaussian measure.

Let us already emphasize that the probability measures on the configuration space, such as~\eqref{eq:canonical measure}, are living in a very high dimensional space. The computation of averages with respect to them, as in~\eqref{intro:averages}, is therefore very costly, if at all possible. A key idea in MD is to avoid computing these high dimensional averages by replacing them by a one dimensional average, corresponding to averaging over time an ergodic dynamics (see Section~\ref{sec:sampling_methods}). 
 
\subsection{Paradigmatic systems}
\label{sec:example_systems}

We present in this section the three systems that will be used throughout this review to illustrate the methods we discuss.

\begin{figure}[h]
    \centering
    \includegraphics[width=0.7\linewidth]{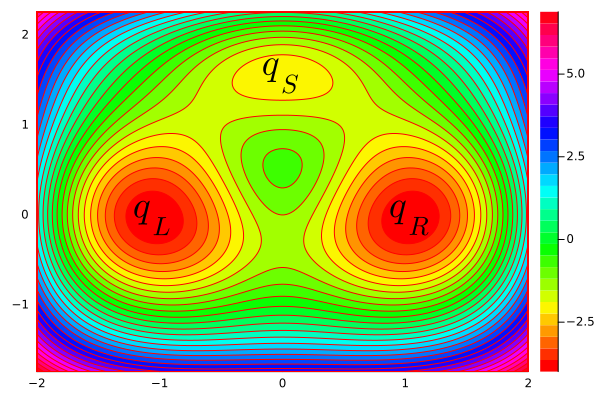}
    \caption{Isolines of the potential energy function~$V$ defined in~\eqref{eq:entropic switch}.}
    \label{fig:entropic switch}
\end{figure}

\paragraph{Two-dimensional entropic switch.}
We start by introducing a toy low-dimensional model for metastable systems \cite{Park_2003,Metzner_2006,C_rou_2011}. A single particle in a two-dimensional space (i.e.,~$\cD = \R^2$) has the potential energy 
\begin{equation}\label{eq:entropic switch}
  \begin{aligned}
    V(x, y) & = 3 \, \rme^{-x^2}\left( \rme^{-\left(y - 1/3 \right)^2} - \rme^{-\left(y - 5/3 \right)^2} \right) - 5 \, \rme^{-y^2} \left( \rme^{-(x-1)^2} + \rme^{-(x+1)^2}\right) \\
    & \ \ + 0.2 \, x^4 +  0.2 \left( y - \frac13\right)^4.
  \end{aligned}
\end{equation}
This energy function, depicted in Figure~\ref{fig:entropic switch}, has two wells associated with global minima, at~$q_\mathrm{L} = (-1.048,-0.0421)$ and $q_\mathrm{R} = (1.048,-0.0421)$, and a saddle point between the two wells at~$q_\mathrm{S} = (0, 1.5371)$. The particle diffuses from~$q_{\rm L}$ to~$q_{\rm R}$ either by passing through~$q_{\rm S}$ or by overcoming the energy barrier between~$q_\mathrm{L}$ and~$q_\mathrm{R}$ directly, the latter possibility becoming more and more probable as the temperature increases.

\paragraph{One-dimensional atom chains.}
\label{sec:atom_chains}
One-dimensional atom chains are commonly used models to understand thermal transport and the validity of Fourier's law (see Section~\ref{subsec:transport_coefficient_examples} for more details). To simplify the presentation, we consider~$N$ classical point particles of identical masses~$m_i=1$.
The configuration of the system is described by the momenta~$p = (p_1, p_2, \dots, p_N) \in \mathbb R^N$, and the vector~$q =(q_1, q_2, \dots, q_N) \in \mathcal D \subseteq \mathbb R^N$ of displacements from equilibrium positions on the lattice $(0, a, 2a,\dots,(N-1)\,a)$ (the value~$a>0$ being irrelevant). 

The specific form of the potential energy function depends on the interatomic potential, the possible interaction with the environment (mimicked by an on-site pinning potential), the boundary conditions (free, fixed or periodic), and any additional perturbation applied to the system. In what follows, we consider the simplest case corresponding to nearest-neighbor interactions only, and no pinning potential. To take into account various boundary conditions, it is convenient to consider two additional fixed fictitious particles at the edges of the chain, whose coordinates are~$(q_i, 0)$ for~$i\in\{0, N+1\}$. Fixed boundary conditions on the left and right correspond to~$q_0=0$ and $q_{N+1}=0$, respectively; free boundary condition on the left and right correspond to~$q_0=q_1$ and $q_N = q_{N+1}$, respectively; while periodic boundary conditions are obtained with the choice~$q_0 = q_N$ and $q_1 = q_{N+1}$. In any case, the Hamiltonian of the system reads
\begin{equation}
  \label{eq:H(q,p)_general}
  H(q,p) = \sum_{i=1}^N \frac{p^2_i}{2 m_i} + \sum_{i=1}^{N+1} v(q_i-q_{i-1}),
\end{equation}
for some elementary interaction potential~$v : \mathbb{R} \to \mathbb{R}$. Without loss of generality, we assume that the interaction potential~$v$ has its global minimum at~0 and that $v(0)=0$. For later use, we define the energy of the $i$-th particle in the bulk:
\begin{equation}\label{eq:ei}
  e_i = \frac{p^2_i}{2 m_i} + \frac{1}{2}( V(q_{i} - q_{i-1}) + V(q_{i+1} - q_i) ),\qquad i \in \{1,\dots,N\}.
\end{equation}
Among the most paradigmatic models in this setting are the Fermi--Pasta--Ulam--Tsingou\footnote{Originally introduced to investigate thermalization in nonlinear lattices, but instead revealed unexpected quasi-recurrences, raising deep questions about the mechanisms of thermalization and the validity of statistical assumptions in nonlinear systems.} (FPUT) chain~\cite{FPU1955,Gallavotti2008}, for which the elementary interaction potential is 
\begin{equation}\label{eq:V_FPUT}
  v(r) = \frac{r^2}{2} + \frac{\alpha r^3}{3}  + \frac{\beta r^4}{4}, \qquad \alpha,\beta >0,
\end{equation}
\begin{figure}[t]
  \centering
  \includegraphics[width = 0.9\textwidth]{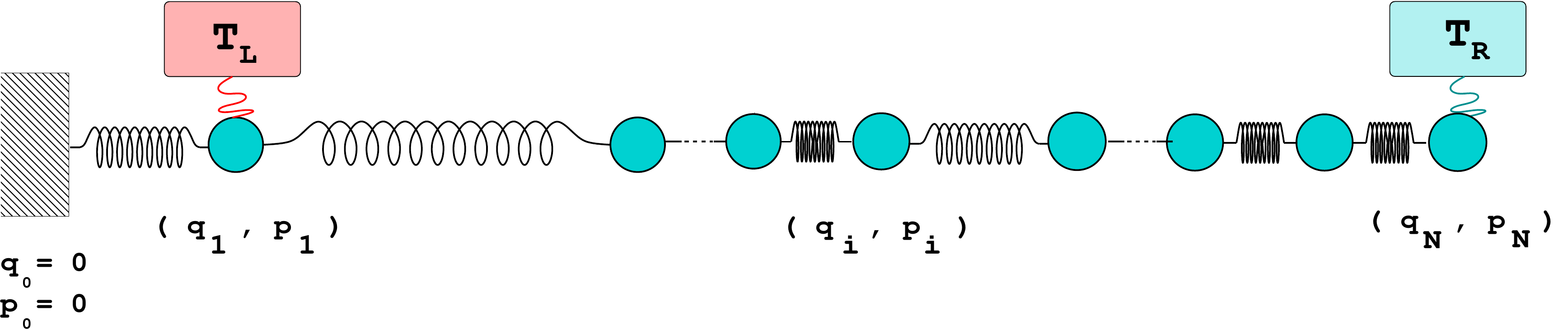}
  \caption{\label{fig:FPUT} An example of oscillator chain, with fixed boundary conditions on the left, free on the right, heat baths at both boundaries, with temperatures $T_\rmL$ and $T_\rmR$ and a mechanical forcing on the rightmost atom.}
\end{figure}
and the rotor chain~\cite{GS00,GLPV00}, for which~$q \in (2 \pi \mathbb{T})^N$, and  
\begin{equation}\label{eq:V_rotor}
  v(r) = 1 - \cos{(r)}.
\end{equation} 
 \begin{figure}[t]
  \centering
  \includegraphics[width = 0.9\textwidth]{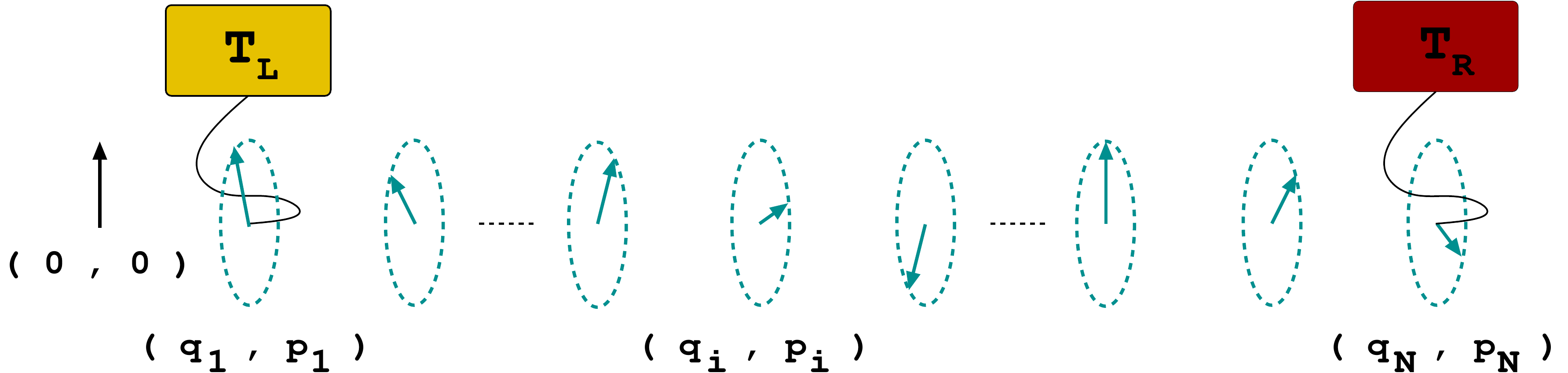}
  \caption{\label{fig:rotor} An example of rotor chain, with free boundary conditions, heat baths with temperatures $T_\rmL$ and $T_\rmR$ and mechanical forcings at both boundaries. This model has been studied for instance in~\cite{IOS2021}.}
 \end{figure}
See Figures~\ref{fig:FPUT} and~\ref{fig:rotor} for illustrations of these systems.

\paragraph{Lennard-Jones fluid.}
The Lennard-Jones potential, introduced in~\cite{lennard-jones_1924_i,lennard-jones_1924_ii}, is a pair potential that describes the basic interaction features between two particles: they repel one another at short range (modeling the London dispersion force as part of the van der Waals interactions) and attract each other at moderate range. It realistically models interactions between atoms in simple liquids and noble gases such as Argon~\cite{rahman_1964,barker_1971,rowley_1975}. 
It is also used as a building block for more complex potentials, e.g.,~to model non-bonded interactions. The form of the potential function is simple, as it is based on only two parameters, and has therefore been extensively studied in the literature, see~\cite{stephan_2020,fischer_2023} and references therein.

\begin{figure}
    \centering
    \includegraphics[width=\columnwidth]{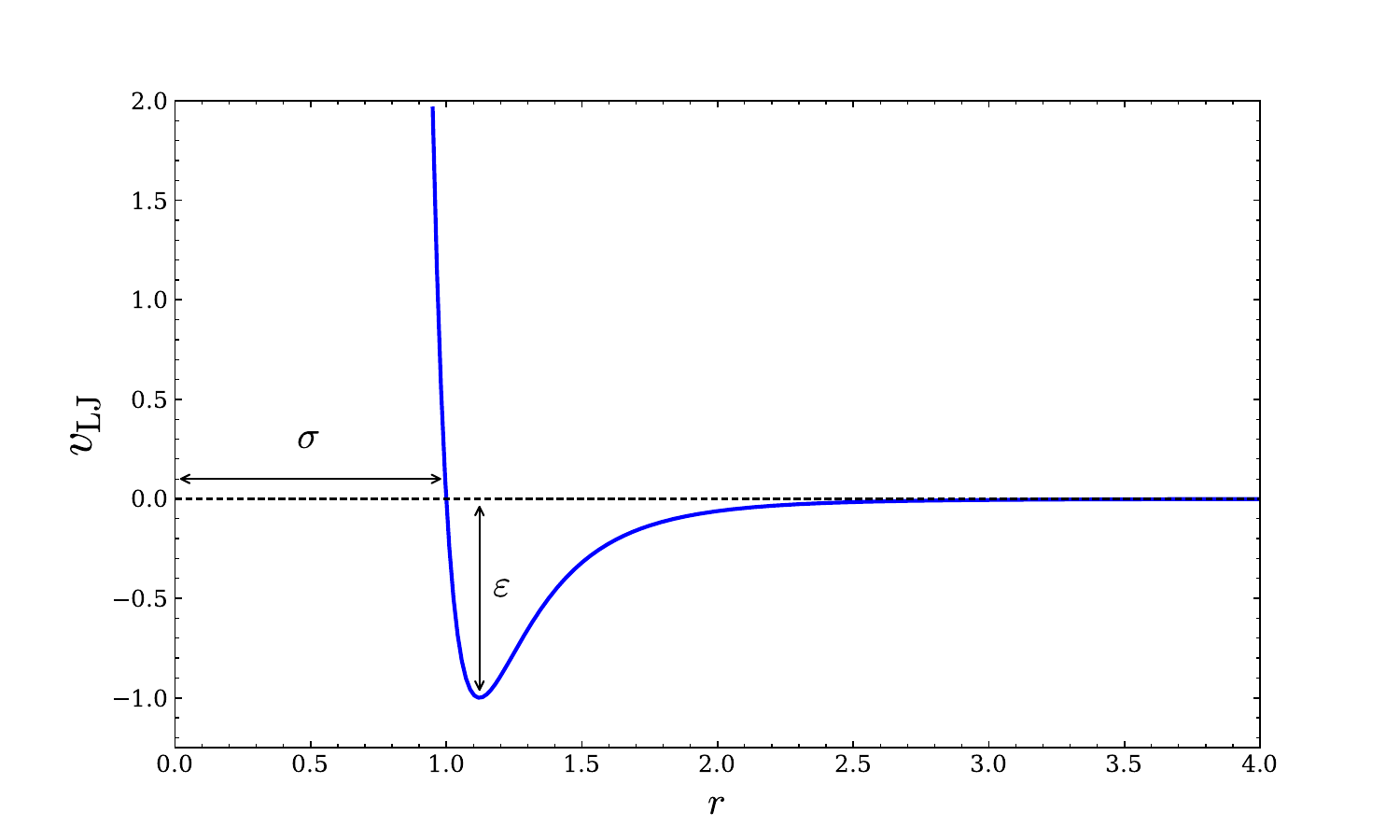}
    \caption{Lennard-Jones potential~$v_{\rm LJ}$ with~$\varepsilon=1$ and~$\sigma=1$.}
    \label{fig:LJ_potential}
\end{figure}

For simplicity, we only consider three-dimensional fluids formed by~$N$ atoms of a single species of atomic mass~$m$. The positions and momenta are denoted by
\begin{equation}
    q = \left(q_{x,j},q_{y,j},q_{z,j}\right)_{1\leq j\leq N},\qquad p = \left(p_{x,j},p_{y,j},p_{z,j}\right)_{1\leq j\leq N}.
\end{equation}
The potential energy function is the sum of pairwise elementary contributions: 
\begin{equation}
  V(q) = \sum_{1\leqslant i<j\leqslant N}v_{\rm LJ}\left(\left\lVert (q_{x,i},q_{y,i},q_{z,i})-(q_{x,j},q_{y,j},q_{z,j})\right\rVert\right),
\end{equation}
for the elementary interaction
\begin{equation}
  v_{\rm LJ}(r)=4\varepsilon\left[\left(\frac{\sigma}{r}\right)^{12}-\left(\frac{\sigma}{r}\right)^6\right].
\end{equation}
The Lennard-Jones potential~$v_{\rm LJ}$ has value~0 at~$r=\sigma$ and is minimal at~$r=2^{1/6}\sigma$, where it has value~$-\varepsilon$. A graphical illustration of this potential is presented in Figure~\ref{fig:LJ_potential}. In practice, in order to comply with periodic boundary conditions (to avoid particles interacting with more than one of their periodic images), the elementary interaction is smoothly truncated to~0 at a finite range of the order of a few multiples of~$\sigma$. 

\subsection{Sampling methods}
\label{sec:sampling_methods}

Several methods have been developed to sample the canonical probability distribution~\eqref{eq:canonical measure} to approximate the high dimensional integrals~\eqref{intro:averages}. Probabilistic methods are commonly used, relying for instance on Markov chain Monte Carlo (MCMC) algorithms~\cite{metropolis_1953,hastings_1970,robert_2004,brooks_2011}, as they do not suffer from the curse of dimensionality as much as deterministic methods (e.g.,~quadratures). We focus here on methods relying on the overdamped and underdamped Langevin dynamics, which are the basis for the nonequilibrium dynamics we consider. Let us emphasize that there are numerous other methods to sample canonical probability distributions, including piecewise-deterministic Markov processes~\cite{davis_1984,michel_2014} or generative models (such as Generative Adversarial Networks~\cite{goodfellow_2014}, Variational Autoencoders~\cite{kingma_2013}, energy-based models~\cite{ackley_1985}, normalizing flows including Boltzmann generators~\cite{tabak_2010,noe_2019,wirnsberger_2020,kobyzev_2021,papamakarios_2021}, score-based models such as diffusion models~\cite{sohl-dickstein_2015} and stochastic interpolants~\cite{albergo_2023}).

Langevin dynamics are a particular class of stochastic differential equations (SDEs). SDEs on~$\calX$ are of the general form
\begin{equation}\label{eq:formal SDE}
  \d x_t = b(x_t)\,\d t + \sigma(x_t) \, \d W_t,
\end{equation}
where~$(x_t)_{t \geqslant 0}$ is a realization of the overdamped or underdamped Langevin dynamics, $(W_t)_{t \geq 0}$ is a standard Brownian motion, and~$b$ and~$\sigma$ are drift and noise functions, respectively, which are chosen so that the canonical measure~\eqref{eq:canonical measure} is invariant for the process~\eqref{eq:formal SDE}. Under additional assumptions on~$b$ and~$\sigma$, the process~\eqref{eq:formal SDE} is ergodic (see for instance~\cite{ikeda_1981,kliemann_1987,khasminskii_2012_i}). In our context this means that, for any observable~$\varphi\in L^1(\pi)$ and initial condition of the dynamics, it holds (see for instance~\cite{kliemann_1987})
\begin{equation}\label{eq:ergodic_theorem}
  \lim_{T \to +\infty} \frac{1}{T}\int_0^{T} \varphi(x_t) \, \d t = \mathbb{E}_\pi[\varphi] \qquad \text{almost surely}.
\end{equation}
The spatial average~$\bbE_\pi\left[\varphi\right]$ can therefore be approximated by using finite-time trajectorial averages. In the following, we introduce Langevin dynamics both in its underdamped and overdamped versions, and present standard numerical schemes to discretize them.

\paragraph{Underdamped Langevin dynamics.}
The underdamped Langevin dynamics is defined on~$\calE$ as 
\begin{equation}\label{eq:langevin}
  \left\{
  \begin{aligned}
    \d q_t & = M^{-1}p_t\,\d t,\\
    \d p_t & = -\nabla V(q_t)\,\d t - \gamma M^{-1} p_t \,\d t + \sqrt{\frac{2\gamma}{\beta}}\,\d W_t, 
  \end{aligned}
  \right.
\end{equation}
where $\gamma>0$ is a friction parameter and~$(W_t)_{t\geqslant 0}$ is a standard~$d$-dimensional Brownian motion. This dynamics is simply the Hamiltonian dynamics for the Hamiltonian~\eqref{eq:hamiltonian} coupled with an Ornstein--Uhlenbeck (OU) process on the momenta whose coefficients are dictated by the fluctuation-dissipation relation in order to ensure that the Boltzmann--Gibbs probability measure~$\mu$ in~\eqref{eq:canonical measure} is preserved. The infinitesimal generator of~\eqref{eq:langevin} reads
\begin{equation}\label{eq:langevin_generator}
  \calL_0 = \calL_\mathrm{ham} + \gamma \calL_\mathrm{OU},
  \qquad
  \calL_\mathrm{ham} = p^\top M^{-1}\nabla_q - \nabla V^\top\nabla_p,
  \qquad
  \calL_\mathrm{OU} = -p^\top M^{-1}\nabla_p + \frac{1}{\beta}\Delta_p.
\end{equation}
This implies that~$\calL_0^\dagger\mu=0$, where~$\calL_0^\dagger$ denotes the~$L^2$-adjoint of~$\calL_0$, which shows the invariance of~$\mu$ with respect to the underdamped Langevin dynamics as a consequence of the Fokker--Planck equation (see for instance~\cite[Section~2.1]{lelievre2016}). The condition~$\calL_0^\dagger\mu=0$ can be equivalently reformulated as follows: for any test function~$\varphi$ (for instance smooth and compactly supported),
\begin{equation}\label{eq:invariance_mu_by_L}
  \int_{\calE} \calL_0\varphi \, \rmd\mu = 0.
\end{equation}
This fact, together with the hypoellipticity of the dynamics ensures by~\cite{kliemann_1987} that the process~\eqref{eq:langevin} is ergodic with respect to~$\mu$ in the sense of~\eqref{eq:ergodic_theorem}. For further use, it is useful to rewrite the generator on the space~$L^2(\mu)$ using
\begin{equation}\label{eq:generator_with_adjoints}
  \calL_\mathrm{ham} = \frac1\beta \sum_{i=1}^d \partial_{p_i}^\star \partial_{q_i}-\partial_{q_i}^\star \partial_{p_i},
  \qquad
  \calL_\mathrm{OU} = \frac1\beta \sum_{i=1}^d \partial_{p_i}^\star \partial_{p_i},
\end{equation}
where~$A^\star$ denotes the adjoint of the (closed) operator~$A$ on~$L^2(\mu)$. In particular, $\partial_{q_i}^\star=-\partial_{q_i}+\beta\partial_{q_i}V$ and $\partial_{p_i}^\star=-\partial_{p_i}+\beta(M^{-1}p)_i$. This highlights the fact that the generator is the sum of an antisymmetric term (the Hamiltonian part) and a symmetric one (the fluctuation/dissipation part). From the viewpoint of the stochastic process, this translates into reversibility up to momentum reversal (as the adjoint~$\cLham^\star = -\calL_\mathrm{ham}$ of the Hamiltonian part can be equivalently obtained by composing on both sides~$\cLham$ with the operator of momentum reversal). 

Classical discretizations of~\eqref{eq:langevin} arise from splitting strategies, namely integrating the OU process and the Hamiltonian dynamics successively at each time step. We refer to~\cite{leimkuhler_2015} and~\cite[Section~3.3.3]{lelievre2016} for further details on Strang and Lie--Trotter splittings for the Langevin dynamics, and simply present one possible discretization here. The OU process can be analytically integrated. As for the Hamiltonian dynamics, the Störmer--Verlet integrator is mainly used, combining symplectic Euler schemes~\cite{hairer_2006} in order to obtain a second-order accuracy with respect to the time step. As such, a currently popular numerical integrator is the so-called BAOAB scheme~\cite{leimkuhler_2015}: for a given starting configuration~$(q^0,p^0) \in \mathcal{E}$, 
\begin{equation}\label{eq:BAOAB}
  \left\{
  \begin{aligned}
    p^{1/3} & = p^{0} - \frac{\Delta t}{2}\nabla V(q^0),\\
    q^{1/2} & = q^0 + \frac{\Delta t}{2}M^{-1}p^{1/3},\\
    p^{2/3} & = \rme^{-\gamma M^{-1}\Delta t}p^{1/3} + \sqrt{\frac{ \left( 1-\rme^{-2\gamma M^{-1}\Delta t}\right)M}{\beta}} \, G^{1},\\
    q^1 & = q^{1/2} + \frac{\Delta t}{2}M^{-1}p^{2/3},\\
    p^{1} & = p^{2/3} - \frac{\Delta t}{2}\nabla V(q^1),
  \end{aligned}
  \right.
\end{equation}
where $G^{1}$ is a standard Gaussian vector. Note that only one force computation per step is required when forces are stored.

\paragraph{Overdamped Langevin dynamics.}
The overdamped Langevin dynamics on~$\calD$ reads
\begin{equation}
  \label{eq:overdamped_langevin}
  \d q_t = -\nabla V(q_t)\,\d t + \sqrt{\frac{2}{\beta}}\,\d W_t,
\end{equation}
which can be seen as a perturbed gradient descent on~$V$. It is obtained from the underdamped Langevin dynamics either through a large friction limit of~$\gamma\to+\infty$ with time rescaling or a small mass limit (see for instance~\cite[Section~2.2.4]{lelievre_2010}). Its infinitesimal generator writes
\begin{equation}\label{eq:overdamped_generator}
  \cL_0 = -\nabla V \cdot \nabla + \beta^{-1}\Delta.
\end{equation}
A simple computation shows that~$\cL_0^\dagger \nu = 0$, which implies that the overdamped Langevin dynamics preserves the probability measure~$\nu$ defined in~\eqref{eq:canonical measure}. Similarly to Langevin dynamics, it is useful to write the generator~\eqref{eq:overdamped_generator} in the weighted space~$L^2(\nu)$ (with adjoints taken on this space):
\begin{equation}\label{eq:overdamped_generator_weighted_L2}
  \cL_0=-\beta^{-1}\nabla^\star\nabla=-\beta^{-1}\sum_{i=1}^{d}\partial_{q_i}^{\star}\partial_{q_i}.
\end{equation}
This highlights the fact that the generator is symmetric and in fact self-adjoint~\cite{bakry-gentil-ledoux-14}, which corresponds to the stochastic dynamics being reversible.

To discretize~\eqref{eq:overdamped_langevin}, a simple Euler--Maruyama discretization starting from~$q^0\in\calD$ with time step~$\Delta t$ can be used:
\begin{equation}\label{eq:overdamped Euler Maruyama}
  q^{1} = q^0 - \Delta t \nabla V(q^0) + \sqrt{\frac{2\Delta t}{\beta}}G^1.
\end{equation}

\paragraph{Invertibility of~$\calL_0$.}
In order to establish error estimates (see Section~\ref{sec:equilibrium_error_estimates}) and derive linear response formulas (see Section~\ref{subsec:linear_response_theory}), a key ingredient is the invertibility of the generator~$\calL_0$ on a suitable functional space. It is useful to work on a subspace of~$L^2(\pi)$, as this is the functional space which naturally arises for central limit theorems for time averages (see Section~\ref{sec:equilibrium_error_estimates}).

Let us first motivate that the correct functional space to consider is a subspace of~$\pi$-centered observables~$L^2(\pi)$, namely
\begin{equation}
  L^2_0(\pi) = \left\{ \varphi \in L^2(\pi) \, \middle | \, \int_\cX \varphi \, \rmd\pi = 0 \right\}. 
\end{equation}
Indeed, a necessary condition for~$\Phi \in L^2(\pi)$ to be a solution to the Poisson equation~$-\calL_0 \Phi = \varphi$ is that~$\varphi$ has average~0 with respect to~$\pi$ in view of the invariance of~$\pi$ as expressed by conditions such as~\eqref{eq:invariance_mu_by_L}. This motivates looking for decay estimates on the semigroup~$\left(\rme^{t\calL_0}\right)_{t\geqslant0}$ as an operator on~$\mathcal{B}(L^2_0(\pi))$, the Banach space of bounded linear operators on~$L^2_0(\pi)$, typically of the form
\begin{equation}\label{eq:exp_decay_semigroup}
  \left\|\rme^{t \calL_0} \right\|_{\mathcal{B}(L^2_0(\pi))} \leq C \rme^{-\kappa t},
\end{equation}
for~$\kappa > 0$ and~$C \in \R_+$. When such decay estimates are available (and more generally when the operator norm of the semigroup is integrable), the generators in~\eqref{eq:langevin_generator} and~\eqref{eq:overdamped_generator} are invertible on~$L^2_0(\pi)$ with
\begin{equation}\label{eq:generator_inverse}
  \calL_0^{-1}=-\int_0^{+\infty}\rme^{t\calL_0} \, \rmd t.
\end{equation}
Decay estimates such as~\eqref{eq:exp_decay_semigroup} hold for overdamped Langevin dynamics when the configurational Boltzmann--Gibbs measure~$\nu$ defined in~\eqref{eq:canonical measure} satisfies a Poincar\'e inequality, see for instance~\cite[Corollaries~2.4 and~2.21]{lelievre2016}; and hold as well for underdamped Langevin dynamics under mild additional assumptions on the potential energy function~$V$ thanks to hypocoercive estimates (as popularized by Villani's monograph~\cite{Villani09}; see the introduction of~\cite{BFLS20} for a review and a historical perspective on these methods). The identity~\eqref{eq:generator_inverse} is useful to derive the Green--Kubo formula for nonequilibrium dynamics, see in particular Section~\ref{subsec:linear_response_theory} below.

\paragraph{Extensions of Langevin dynamics.}
Although we focus on (overdamped) Langevin dynamics with constant diffusion coefficients in this review, let us mention various other dynamics that can be considered, as they may be more relevant from a modeling viewpoint for some applications. A first class of dynamics are overdamped Langevin dynamics with a position-dependent diffusion coefficient (with an extra term in the drift part in order to ensure the invariance of the Boltzmann--Gibbs measure)~\cite{Bennett_1975,jardat_1999,lelievre2024,lelievre_2025_ii,lelievre_2025}. Another class of dynamics are underdamped Langevin dynamics with position-dependent friction and diffusion coefficient, as in dissipative particle dynamics (DPD) for instance~\cite{hoogerbrugge_1992,espanol_1995}. One can also consider generalized Langevin dynamics, which are even more inertial than underdamped Langevin dynamics (see for instance~\cite{OP11}). Finally, Langevin-like dynamics can be restricted to live on submanifolds~\cite{lelievre_2012,lelievre_2019} to describe physical systems with constraints. The constraints correspond to fixed values of molecular constraints (such as bond lengths and bond angles), in which case the dynamics are typically integrated using implicit numerical schemes such as SHAKE or RATTLE~\cite{ryckaert_1977,andersen_1983}.

\subsection{Error estimates}
\label{sec:equilibrium_error_estimates}

We discuss in this section the numerical errors which arise when estimating in practice ergodic averages of the form~\eqref{eq:ergodic_theorem}. In order to better distinguish between the various sources of errors, we start by making precise in Section~\ref{eq:sampling_error_continuous} the bias and statistical error due to the finite integration times. We then discuss in Section~\ref{sec:error_dt_equilibrium} extra errors related to the time discretization of the dynamics.

\subsubsection{Error at the continuous time level}
\label{eq:sampling_error_continuous}

We consider in this section the following estimator of the ergodic average~\eqref{eq:ergodic_theorem}:
\begin{equation}\label{eq:estimator_continuous_time}
  \widehat{\varphi}_T = \frac{1}{T}\int_0^{T} \varphi(x_t)\,\d t.
\end{equation}
The main result is that the latter estimator has a bias of order~$\mathrm{O}(1/T)$ and a statistical error of order~$\mathrm{O}(1/\sqrt{T})$ (i.e.,~the variance of~$\widehat{\varphi}_T$ scales as~$1/T$ for long times). The key computation to prove these results is based on It\^o calculus on~$\Phi(x_t)$, with~$\Phi$ the solution to the Poisson equation (see the discussion around~\eqref{eq:generator_inverse} for sufficient conditions on the well posedness of this equation in the case of underdamped and overdamped Langevin dynamics) 
\begin{equation}\label{eq:Poisson_equation_continuous_time_average}
  - \cL_0 \Phi = \Pi\varphi \in L^2_0(\pi),
\end{equation}
with
\begin{equation}\label{eq:Pi_centering}
  \Pi \varphi = \varphi - \int_\cX \varphi \, \rmd \pi.
\end{equation}
More precisely, for the dynamics~\eqref{eq:formal SDE}, it holds $\rmd \Phi(x_t) = \cL_0 \Phi(x_t) \, \rmd t + \nabla \Phi(x_t)^\top \sigma(x_t) \, \rmd W_t$, so that
\begin{equation}\label{eq:Ito_time_averages}
  \frac{1}{T}\int_0^T \left( \varphi(x_t) - \int_\cX \varphi \, \rmd\pi \right) \, \d t = \frac{\Phi(x_0) - \Phi(x_T)}{T} + \frac{1}{T} \int_0^T \nabla \Phi(x_t)^\top \sigma(x_t) \, \d W_t.
\end{equation}
As the second term is a martingale, taking the expectation directly gives the following estimate on the bias of the estimator (upon deriving moment estimates on~$\Phi$ and~$\nabla \Phi$, typically through Lyapunov estimates~\cite{bellet_2006,kopec1,kopec2}):
\begin{equation}
\mathbb E {\left[ \widehat{\varphi}_T \right]} - \int_\cX \varphi \, \rmd\pi = \mathrm{O} \left( \frac{1}{T} \right).
\end{equation}
To quantify the statistical error, we rely on a central limit theorem for martingales (see for instance~\cite{KomorowskiLandimOlla2012} or~\cite{EthKur86}), after noting that the first term on the right-hand side of~\eqref{eq:Ito_time_averages}, whose variance scales as~$1/T^2$, is negligible with respect to the Brownian integral, whose variance scales as~$1/T$. More precisely,
\begin{equation}
  \mathrm{Var} \left( \widehat{\varphi}_T \right) = \frac{v_\varphi}{T},
  \qquad
  v_\varphi = \int_\cX \nabla \Phi^\top \sigma\sigma^\top \nabla \Phi \, \rmd \pi.
\end{equation}
The asymptotic variance~$v_\varphi$ can in fact be reformulated using a time-integrated autocorrelation (similar to the one appearing in the Green--Kubo formula~\eqref{eq:green_kubo}) as
\begin{equation}\label{eq:asymptotic_variance_continuous}
  v_\varphi = 2\int_\cX \Phi (-\cL_0 \Phi) \, \rmd \pi = 2\int_\cX \varphi (-\cL_0^{-1} \Pi \varphi) \, \rmd \pi = 2\int_0^{+\infty} \mathbb{E}(\Pi\varphi(x_t)\Pi\varphi(x_0)) \, dt.
\end{equation}
In fact, it can be shown that a central limit theorem holds provided~$\Phi \in L^2(\mu)$ (see for instance~\cite{Bhattacharya_1982}).

\subsubsection{Error at the discrete time level}
\label{sec:error_dt_equilibrium}

\begin{figure}
    \centering
    \includegraphics[width=0.9\columnwidth]{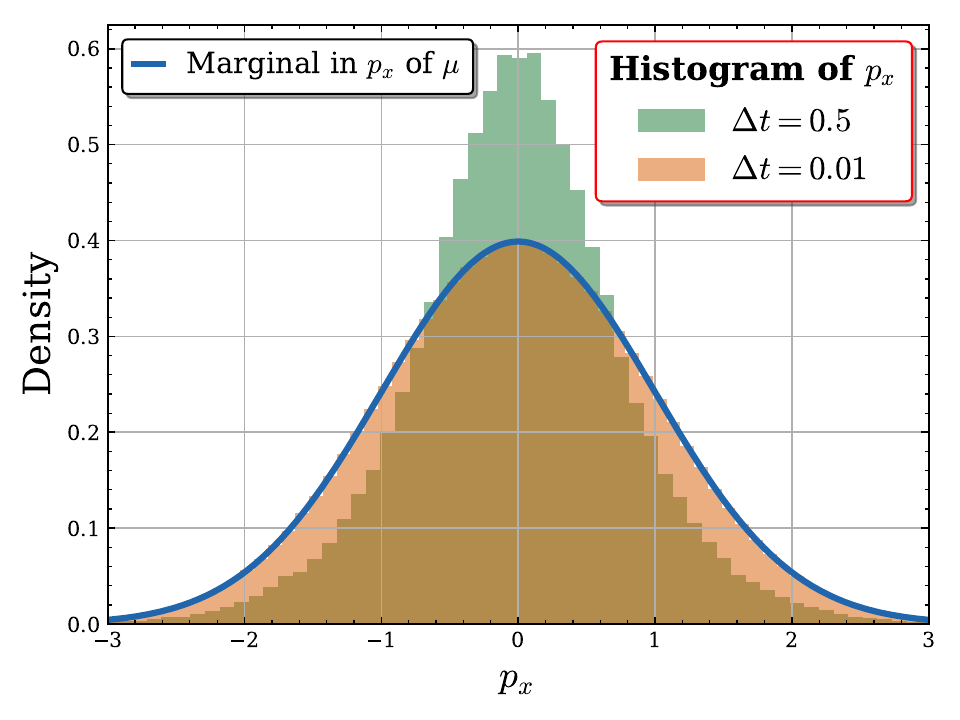}
    \caption{Comparison of the sampled $x$-component~$p_x$ of the momenta distribution for two time steps~$\Delta t$ with the analytical distribution. The dynamics are run with the BAOAB scheme~\eqref{eq:BAOAB} with~$\gamma=1$ and~$\beta=1$, integrating up to~$T=100{,}000$.}
    \label{fig:equilibrium_timestep_bias_momenta_distribution}
\end{figure}

\begin{figure}
    \centering
    \includegraphics[width=0.9\columnwidth]{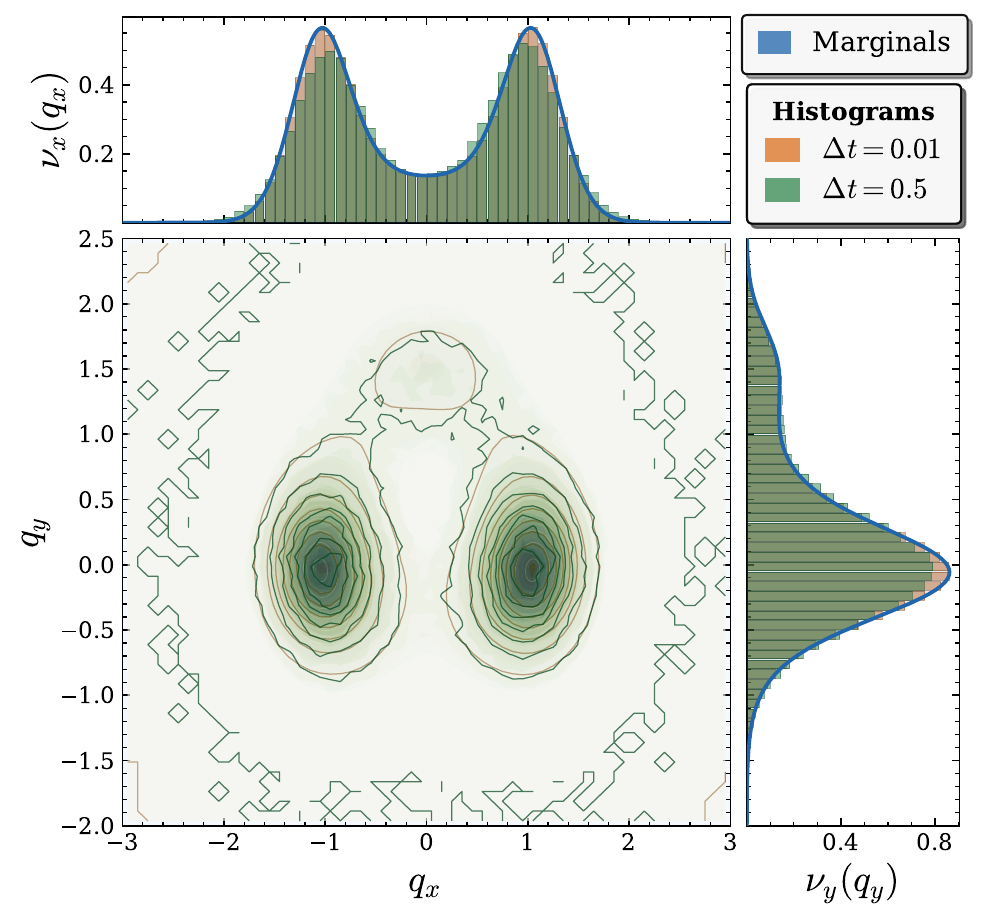}
    \caption{Comparison of the sampled positions in a similar setting as in Figure~\ref{fig:equilibrium_timestep_bias_momenta_distribution}. The 2D histograms have been smoothed using a Gaussian kernel of standard deviation~$1.5$ for better readability.}
    \label{fig:equilibrium_timestep_bias_positions_distribution}
\end{figure}

\begin{figure}
    \centering
    \begin{subfigure}[t]{0.475\columnwidth}
        \includegraphics[width=\columnwidth]{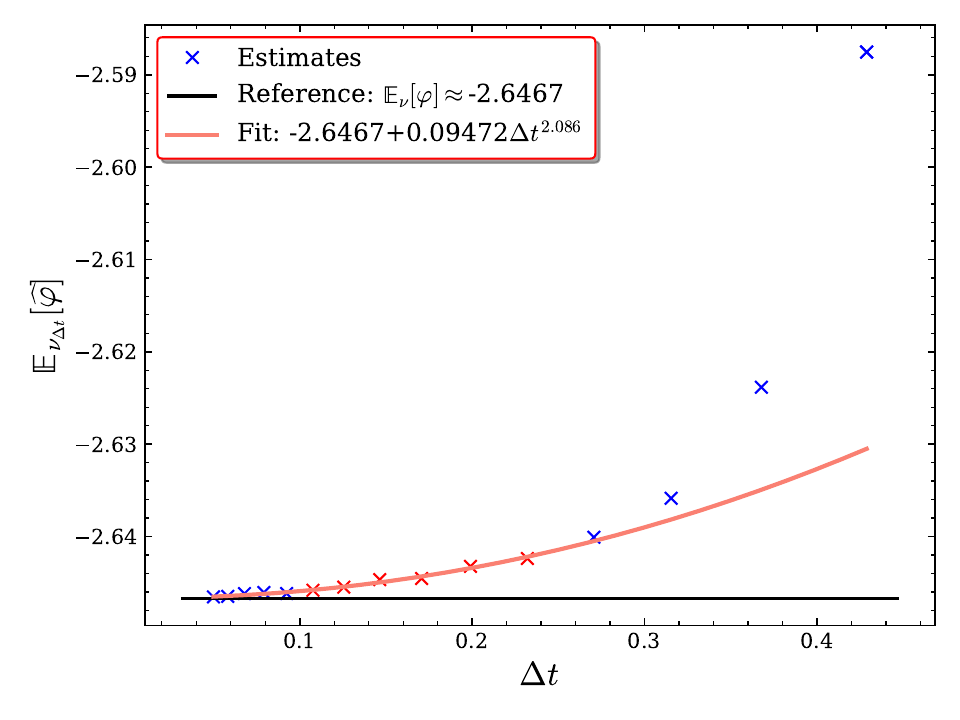}
        \caption{Estimates as a function of the time step~$\Delta t$.}
    \end{subfigure}\hfill
    \begin{subfigure}[t]{0.475\columnwidth}
        \includegraphics[width=\columnwidth]{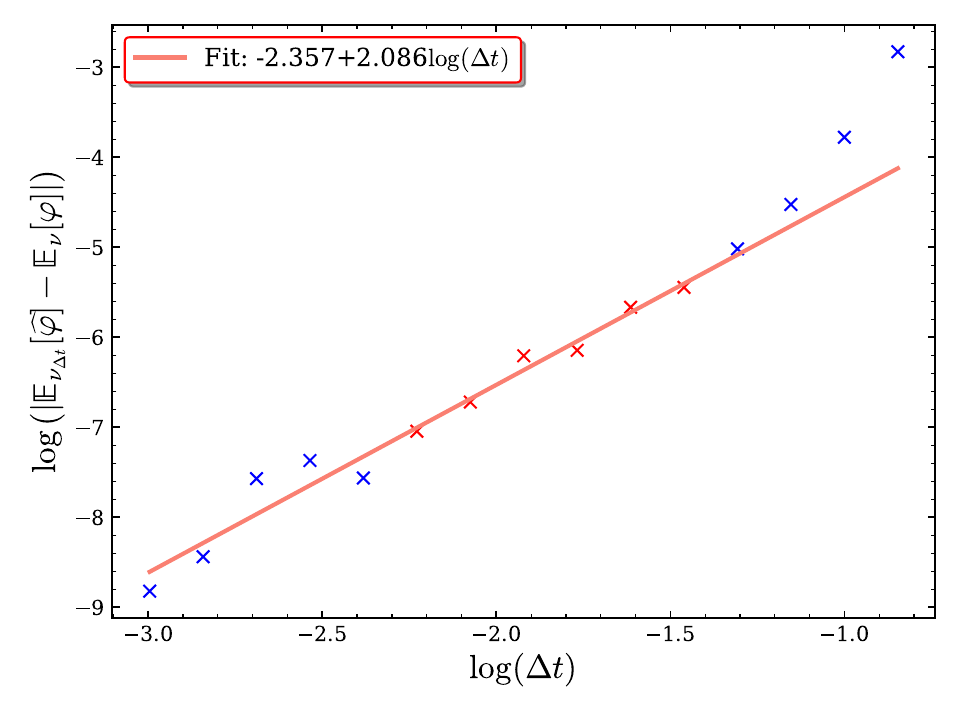}
        \caption{Error between the estimates and the reference value with the time step~$\Delta t$ in log-scale.}
        \label{fig:equilibrium_timestep_bias:b}
    \end{subfigure}
    \caption{Estimates of~$\bbE_\nu\left[\varphi\right]$ with~$\varphi=V$ defined in~\eqref{eq:entropic switch} as a function of~$\Delta t$. Each estimate is obtained by averaging 100 independent estimates obtained by integrating the Langevin dynamics with the BAOAB scheme~\eqref{eq:BAOAB} with~$\gamma=1$ and~$\beta=1$, integrating up to~$T=1{,}000{,}000$. A linear fit is performed in log scale (right panel) on the data points shown in red, which illustrates that the BAOAB scheme leads to a second-order consistent approximation of the target value. Data associated with smaller (resp. larger) time steps were not used in the fitting procedure because of larger statistical errors (resp. bias).}
    \label{fig:equilibrium_timestep_bias}
\end{figure}

In practice, in addition to running simulations in finite time, one also needs to discretize the dynamics. This adds another source of error, related to the time step~$\Delta t$, as the invariant probability measure~$\pi$ of the continuous dynamics~\eqref{eq:formal SDE} is not invariant by the discretized dynamics in general. The discretized dynamics admits, under some conditions on $\Delta t$, 
$b$ and~$\sigma$, a unique invariant probability measure~$\pi_{\Delta t}$, which depends on the time step, and is ergodic for this measure (see for instance~\cite{MSH02,BO10,kopec1,kopec2,DEMS25} and references therein for some results in this direction in the context of Langevin dynamics). When discretizing time, one replaces the estimator~\eqref{eq:estimator_continuous_time} by
\begin{equation}\label{eq:estimator_discrete_time}
  \widehat{\varphi}_{\Niter,\dt} = \frac{1}{\Niter}\sum_{n=1}^\Niter \varphi(x^n),
\end{equation}
for a discrete trajectory~$(x^n)_{n \geq 1}$. Ergodicity means that the following limit holds for any initial condition~$x^0$:
\begin{equation}
  \widehat{\varphi}_{\Niter,\dt} \xrightarrow[\Niter \to +\infty]{} \int_\cX \varphi \, \d \pi_{\Delta t} \qquad \mathrm{a.s.}.
\end{equation}
Typically, $\pi_{\Delta t}$ is a small perturbation (in $\Delta t$) of $\pi$. More precisely, for an observable $\varphi$, a so-called Talay--Tubaro expansion holds~\cite{TT90} (see for instance~\cite{leimkuhler2016computation} for splitting schemes of Langevin dynamics):  
\begin{equation}\label{eq:error_discretised_pi}
  \int_\cX \varphi \, \d \pi_{\Delta t} = \int_\cX \varphi \, \d \pi + \Delta t^\alpha \int_\cX \varphi f_{\alpha + 1} \, \d \pi + \Delta t^{\alpha+1} R_{\varphi, \Delta t},
\end{equation}
for some remainder term~$R_{\varphi, \Delta t}$ controlled uniformly in $\Delta t$, and $f_{\alpha + 1}$ a function fully determined by the discretization scheme used, with $\alpha \geq 1$ an integer. The proof of such an expansion relies, from an algebraic viewpoint, on an expansion for small~$\dt$ of the discrete transition kernel; while moment estimates on~$\pi_\dt$ and convergence results on the transition operator~$P_\dt$, uniform in the time step, are needed to control remainders. The transition operator is defined as 
\begin{equation}
  (P_{\Delta t}\varphi)(x) = \E\left[\varphi(x^{n+1})\middle|x^n=x\right].
\end{equation}
For underdamped Langevin dynamics, discretized with~\eqref{eq:BAOAB}, a second-order accuracy is obtained, i.e.,~$\alpha = 2$ (see Figure~\ref{fig:equilibrium_timestep_bias:b}), while for overdamped Langevin dynamics with~\eqref{eq:overdamped Euler Maruyama} it holds~$\alpha = 1$. We refer to \cite{lelievre2016} for more precise results. The bias with respect to~$\Delta t$ in the sampled distribution is illustrated in Figure~\ref{fig:equilibrium_timestep_bias_momenta_distribution} and~\ref{fig:equilibrium_timestep_bias_positions_distribution}, and for the estimate~\eqref{eq:error_discretised_pi} in Figure~\ref{fig:equilibrium_timestep_bias}, for the two-dimensional entropic switch system.

\begin{figure}
    \centering
    \includegraphics[width=0.75\columnwidth]{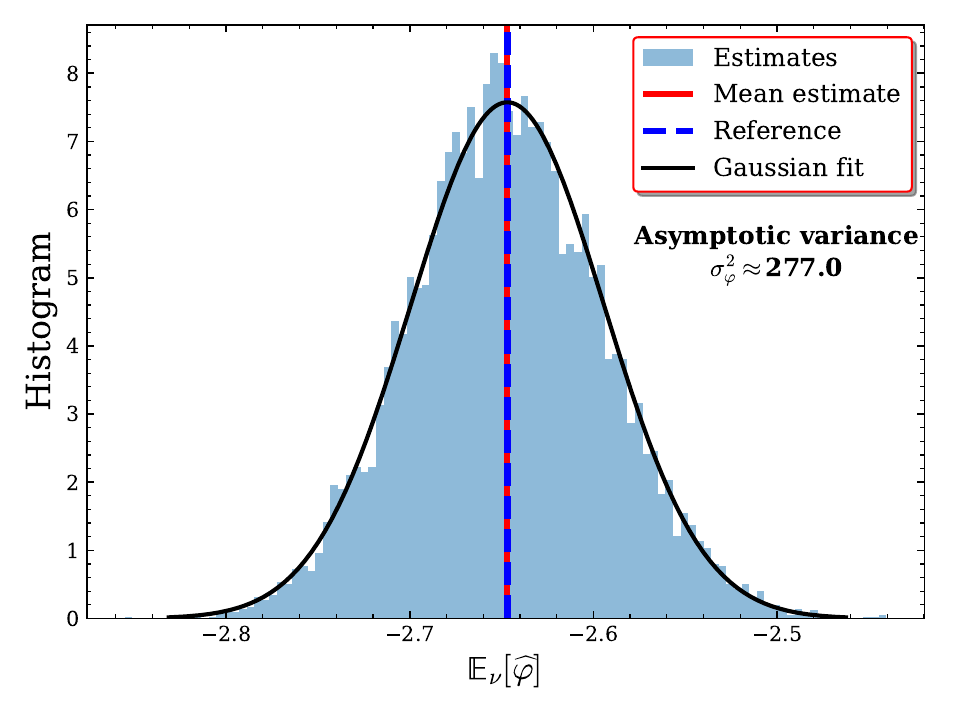}
    \caption{Histogram of 10,000 estimates of~$\bbE_\nu[\varphi]$ with~$\varphi=V$ defined in~\eqref{eq:entropic switch}, obtained by integrating the Langevin dynamics using independent initial configurations under~$\mu$, with~$\Delta t=0.01$, the BAOAB scheme defined in~\eqref{eq:BAOAB} with~$\gamma=1$ and~$\beta=1.0$ and integrating the dynamics up to time~$T=1000$.}
    \label{fig:equilibrium_estimator_variance}
\end{figure}

In practice, one cannot consider the infinite-time limit of~$\Niter\to+\infty$ for~$\widehat{\varphi}_{\Niter,\dt}$. Similarly to the results obtained for the continuous-time dynamics in Section~\ref{eq:sampling_error_continuous}, two additional sources of error should be considered: a bias of order~$1/(\Niter \dt)$, and a variance, of order~$1/(\Niter\dt)$ as well. In order to derive the result, one introduces the discrete counterpart to the Poisson equation~\eqref{eq:Poisson_equation_continuous_time_average}, namely
\begin{equation}
\frac{1-P_\dt}{\dt} \Phi_\dt = \varphi - \int_\cX \varphi\,\d\pi_\dt, 
\end{equation}
and reexpresses the trajectory average by rewriting~$\left[(1-P_\dt)\Phi_\dt\right](x^n)$ as the sum of a telescopic term $\Phi_\dt(x^n)-\Phi_\dt(x^{n+1})$ and a discrete martingale increment
\begin{equation}
  M_n = \Phi_\dt(x^{n+1}) - (P_\dt\Phi_\dt)(x^n),
\end{equation}
which is a standard strategy, used for instance in~\cite{plechac_nemd}. This leads to
\begin{equation}
\begin{aligned}
  \widehat{\varphi}_{\Niter,\dt} - \int_\cX \varphi\,\d\pi_\dt
  & = \frac{1}{\Niter\dt} \sum_{n=1}^{\Niter} \left[(1-P_\dt)\Phi\right](x^n) \\
  & = \frac{\Phi_\dt(x^1)-\Phi_\dt(x^{\Niter+1})}{\Niter\dt} + \frac{1}{\Niter\dt} \sum_{n=1}^\Niter M_n.
\end{aligned}
\end{equation}
The first term on the last right-hand side is a bias, of order~$1/(\Niter\dt)$, while the sum of the discrete martingale increments has a variance of order~$\Niter\dt$ (as can be seen from the fact that~$\E[M_n^2] = \mathrm{O}(\dt)$, an equality which follows from a Taylor expansion), so that the second term on the last right-hand side is of order~$\mathrm{O}(1/\sqrt{\Niter\dt})$. In fact, for any fixed~$\dt>0$, it can be shown that the asymptotic variance of the estimator is (see Section~\ref{sec:dt_bias_GK})
\[
\Niter\dt \mathrm{Var}\left(\widehat{\varphi}_{\Niter,\dt}\right) \xrightarrow[\Niter\to+\infty]{} v_{\varphi,\dt},
\]
with
\[
\dt \int_\cX (\Pi_\dt \varphi)^2 \, \d\pi_\dt + 2 \dt \sum_{n=1}^{+\infty} \E\left[(\Pi_\dt \varphi)(x^n)(\Pi_\dt \varphi)(x^0)\right]
\]
where
\[
\Pi_\dt \varphi = \varphi - \int_\cX \varphi \, \d\pi_\dt.
\]
The results of Section~\ref{sec:dt_bias_GK} show that~$v_{\varphi,\dt}$ agrees at dominant order with~$v_\varphi$; more precisely,
\begin{equation}
  \label{eq:scaling_asymptotic_variance}
  v_{\varphi,\dt} = 2 \int_0^{+\infty} \E\left[(\Pi \varphi)(x_t)(\Pi \varphi)(x_0)\right] dt + \mathrm{O}(\Delta t),
\end{equation}
where~$\Pi$ is the centering operator defined in~\eqref{eq:Pi_centering}. Figure~\ref{fig:equilibrium_estimator_variance} illustrates the statistical fluctuations of a discrete estimator for the two-dimension entropic switch system.

\subsection{Variance reduction methods}
\label{sec:equilibrium_variance_reduction}

A great deal of effort has been devoted to reducing the variance of estimators associated with Markov processes and Markov chains~\cite{fishman_1996,Caflisch_1998,lapeyre_2003,liu_2008,rubinstein_2016}. In view of~\eqref{eq:asymptotic_variance_continuous} and~\eqref{eq:scaling_asymptotic_variance}, the statistical error of averages over discrete trajectories is dictated, to first order in~$\dt$, by the asymptotic variance of averages over continuous trajectories. We present in this section five ways to reduce the asymptotic variance of Monte Carlo estimators in general, and trajectory averages in particular:

\begin{itemize}
\item \textit{stratification} aims to decompose the phase space into several regions, which are sampled independently and reweighted according to the canonical weight of the region itself. Related numerical methods include bridge sampling methods~\cite{bennett_1976,shirts_2008} and thermodynamic integration~\cite{kirkwood_1935,carter_1989,lelievre_2010}; 
\item \textit{control variates} are observables~$\phi$ of known averages under~$\nu$ which are used to reduce the variance of estimators: if~$\bbE_\pi\left[\varphi\right]$ is the target quantity, one chooses~$\phi$ so that the variance of the time average of~$\varphi-\phi$ is minimal. A convenient class of control variates are the images of the generator~$\calL$, i.e., functions of the form~$\phi = \calL\Phi$, which are by construction such that~$\bbE_\pi\left[\phi\right]=0$; see for instance~\cite{lelievre2016,RS19} and references therein. Another possibility is to build control variates based on coupled dynamics, see Section~\ref{sec:noneq_control_variates} for examples related to the estimation of the response of nonequilibrium systems. For further details on control variates for Markov chains, we refer to~\cite{hammersley_1964,andradottir_1993,henderson_2002};
\item \textit{preconditioning} considers a diffusion (equivalent to the inverse of a mass tensor) in Langevin dynamics in order to remove or at least reduce the anisotropy of the potential energy surface~\cite{Bennett_1975,goodman_2010,girolami_2011,patterson_2013,Leimkuhler2018,li_2022,Casas2022,lelievre2024,Tran2024}. The adaptive construction of preconditioners has been suggested in the literature~\cite{Haario1999,haario_2001,Roberts_2009,lelievre_2025}, and the optimization of position-dependent diffusions has been recently studied for overdamped Langevin dynamics~\cite{lelievre_2025_ii,cui_2025};
\item \textit{importance sampling} corresponds to sampling a probability distribution different from~$\nu$ and reweighting the obtained samples accordingly~\cite{chen_1997,liu_2008,tokdar_2010,rubinstein_2016}. A standard procedure in computational statistical physics is to consider an additional potential energy function~$U$ that removes the metastable features of the standard Langevin dynamics associated with the baseline potential energy function~$V$, e.g.,~by ``flattening'' the potential energy surface~$V+U$ in order to remove the original energy barriers. For the method to be effective, the weights~$\rme^{\beta U}$ associated with the additional potential should not blow up, which can be monitored by computing effective sample sizes, namely the number of approximately independent samples obtained from~$\nu$; see for instance~\cite[Section~3.4.3]{lelievre2016} for further details. A typical biasing potential~$U$ is (a fraction of) the free energy associated with some collective variable that describes the difficult (nonlinear) directions to sample; see for instance~\cite{lelievre_2010};
\item{
  \textit{symmetrization}, in the spirit of antithetic variates, exploits symmetries in the system to construct new estimators from existing ones, by relying on the identity
  \begin{equation}
    \E_\pi[\varphi] =  \E_{\pi}\left[\int_G\varphi\circ \psi_g\,\mu_G(\d g)\right],
    \end{equation}
    valid whenever~$g\mapsto \psi_g$ is a group action of $G$ on~$\cX$ leaving~$\pi$ invariant, and $\mu_G$ is any probability measure on~$G$. For physical systems of indistinguishable particles, $G$ is typically generated by a set of rigid motions of the torus on positions, rotations on the momenta, and permutations of particles; see for instance~\cite{MSEDDT22,MESVDDT24} in the context of the computation of shear viscosities. This family of methods is a simple way to ensure some amount of variance reduction, since, by Jensen's inequality,
    \begin{equation}
        \E_{\pi}\left[\left(\int_G\varphi\circ \psi_g\,\mu_G(\d g)\right)^2\right] \leq \E_{\pi}\left[\int_G (\varphi\circ \psi_g)^2\,\mu_G(\d g)\right] = \E_\pi[\varphi^2].
    \end{equation}
    The same method allows to use symmetries in the distribution of trajectories of the dynamics, which provides a simple approach to reducing the variance of dynamical estimators.
    }
\end{itemize}

\section{Transport coefficients}
\label{sec:transport_coefficients}

We now present a mathematical framework for the computation of transport coefficients. We first discuss transport coefficients and their numerical computation in the framework of the stochastic dynamics introduced in Section~\ref{sec:sampling_methods}. For more in-depth discussions of the mathematical aspects of this topic, we refer to~\cite[Section 5]{lelievre2016} and~\cite[Chapter 8]{StoltzM2}, and also the recent review~\cite{stoltz2024}. In Section~\ref{subsec:nemd_dynamics}, we define the nonequilibrium systems we consider in this work. We next give, in Section~\ref{subsec:linear_response_theory}, an overview of linear response theory. Finally, we present in Section~\ref{subsec:transport_coefficient_examples} some typical examples of nonequilibrium systems and their associated transport coefficients.

\subsection{Nonequilibrium dynamics}
\label{subsec:nemd_dynamics}

To formalize the physical notion of the nonequilibrium systems used in the NEMD method, we introduce perturbations of the equilibrium dynamics~\eqref{eq:formal SDE}, given by a parametric family of SDEs evolving on~$\cX$,
\begin{equation}\label{eq:nemd_dynamics}
  \d x^\eta_t = b_\eta(x_t^\eta)\,\d t + \sigma_\eta(x_t^\eta)\,\d W_t,
\end{equation}
indexed by a parameter~$\eta>0$ modulating the strength of the perturbation. The perturbations should be understood as being of order~$\eta$, i.e.,
\begin{equation}\label{eq:linear_forcing}
  \forall\,x\in\cX,\qquad\left|b(x)-b_\eta(x)\right|,\, \left|\sigma(x)-\sigma_\eta(x)\right| = \mathrm{O}\left(\eta\right)\,\text{as }\eta\to 0.
\end{equation}

\paragraph{Two typical situations.}
To be more specific, we focus on the two cases which are the most relevant in practice, although one could in principle consider more general classes of perturbations.
\begin{itemize}
\item \textit{Non-conservative force}. We often consider systems subjected to an external, non-gradient driving force $\eta F$, where~$F:\cD\to\R^d$ is a fixed forcing field. More precisely, this amounts to replacing~$\nabla V(q)$ with
\begin{equation}\label{eq:non_conservative_force}
  \nabla V(q) + \eta F(q)
\end{equation}
in the equilibrium dynamics.
\item \textit{Temperature profile}. The second typical situation consists of fixing a position-dependent temperature field, which amounts to replacing the temperature parameter~$\beta$ with a positive scalar-valued function
\begin{equation}\label{eq:nemd_temperature_profile}
  \beta_\eta(q) = \frac{1}{T_0 +\eta\delta T(q)},\qquad \delta T:\cD\to\R,
\end{equation}
in the equilibrium dynamics, where~$T_0=\beta^{-1}$. The magnitude~$\eta$ should be small enough in order for~$T_0 +\eta\delta T(q)$ to remain positive.
\end{itemize}

\paragraph{Nonequilibrium generator.}
In both examples given above, the generator~$\cL_\eta$ of the nonequilibrium dynamics is a linear perturbation of the equilibrium generator:
\begin{equation}\label{eq:nemd_generator}
  \cL_\eta = \cL_0 + \eta\wcL,
\end{equation}
where the expression of the perturbation~$\wcL$ depends both on the dynamics and on the type of perturbation. Expressions for~$\wcL$ for the two examples above and Langevin-type dynamics are listed in Table~\ref{tab:nemd_generator_perturbations}.
\begin{table}[h]
    \centering
    \begin{tabular}{|c|c|c|}
        \hline
        \backslashbox{Perturbation}{Dynamics} & Overdamped Langevin & Underdamped Langevin\\
        \hline
        Non-conservative force & $F\cdot\nabla$ & $F\cdot\nabla_p$\\
        \hline
        Temperature profile & $\delta T\Delta$ & $\gamma\delta T\Delta_p$ \\
        \hline
    \end{tabular}
    \caption{Nonequilibrium perturbation~$\wcL$ of the generator for usual dynamics and perturbation types.}
    \label{tab:nemd_generator_perturbations}
\end{table}

\paragraph{Nonequilibrium steady states.}
A steady state for the nonequilibrium dynamics~\eqref{eq:nemd_dynamics} is, by definition, an invariant probability distribution~$\pi_\eta(\d x)$, which we also write~$\mu_\eta(\rmd q \, \rmd p)$ or~$\nu_\eta(\rmd q)$ depending on the underlying equilibrium dynamics. The existence and uniqueness of the steady state can typically be proven with Lyapunov techniques, in the spirit of~\cite{hairer_2011}, see~\cite[Section 5]{lelievre2016} and~\cite{bellet_2006,spacek2023} for precise statements.

In turn, the results of~\cite{kliemann_1987} can be used, under suitable conditions on the coefficients~$\sigma_\eta,\,b_\eta$ ensuring the hypoellipticity of the generator~\eqref{eq:nemd_generator} and the positivity of the steady state density, to obtain the following pathwise ergodic property: for any initial condition~$x\in\mathcal X$ and observable~$\varphi\in L^1(\pi_\eta)$,
\begin{equation}
    \label{eq:nemd_ergodicity}
    \frac1T\int_0^T \varphi(x_t^\eta)\,\d t \xrightarrow[T\to+\infty]{} \E_{\pi_\eta}[\varphi]\qquad \mathbb{P}_x\text{-almost surely}.
\end{equation}

From the analytical point of view, the steady state~$\pi_\eta$ is a weak solution of the stationary Fokker--Planck equation
\begin{equation}
    \label{eq:fp_nemd}
    \cL_\eta^\dagger \pi_\eta = 0,
\end{equation}
where~$\cL_\eta^\dagger$ denotes the~$L^2(\cX)$-adjoint of the nonequilibrium generator~\eqref{eq:nemd_generator}. Regularity properties of~$\pi_\eta$ can usually be obtained by applying standard elliptic regularity theory, or H\"ormander's hypoelliptic theory~\cite{hormander_67} if the diffusion matrix~$\sigma_\eta$ is degenerate.

We stress that~$\pi_\eta$ typically does not have an explicit expression: since the non-conservative force~$F$ cannot be written as the gradient of a potential on~$\cD$, and since the fluctuation-dissipation relation is generally not satisfied for a position-dependent temperature profile, one cannot write~$\pi_\eta$ as a Boltzmann--Gibbs measure for an explicit energy function. We can nevertheless, in view of~\eqref{eq:nemd_ergodicity}, sample from the nonequilibrium steady state, by considering sufficiently long trajectories of~\eqref{eq:nemd_dynamics}.

Similarly to the equilibrium case, discretizations of the dynamics lead to a bias on the invariant probability measure sampled by the numerical scheme and the resulting trajectory averages, but, while in the equilibrium case one can correct for this bias via Metropolization, this option is unavailable in the nonequilibrium case in all but the simplest settings, since computing the nonequilibrium probability density up to normalization is equivalent to finding a solution to the Fokker--Planck equation~\eqref{eq:fp_nemd}.

\begin{figure}
    \centering
    \includegraphics[width=0.8\textwidth]{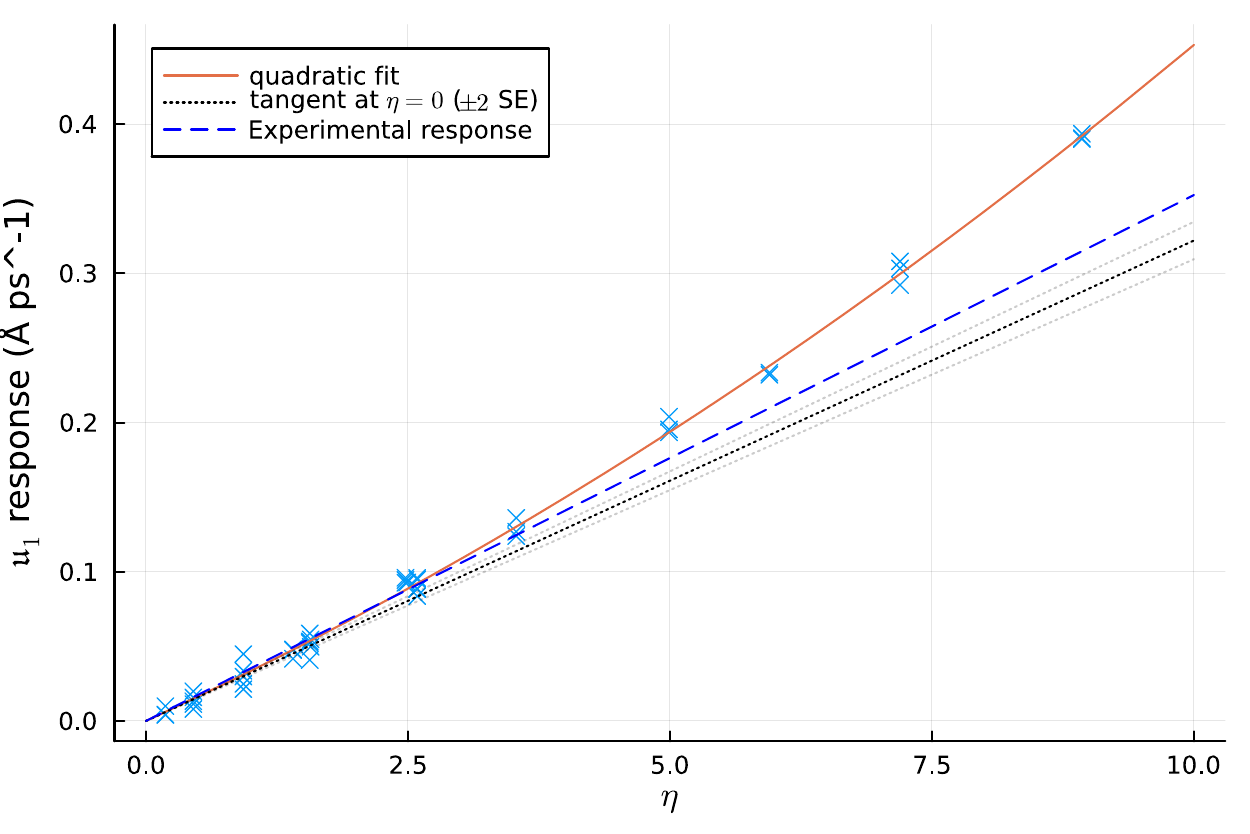}
    \caption{The NEMD fitting procedure used to compute the shear-viscosity (see Section~\ref{subsec:transport_coefficient_examples}) of liquid argon.
    }
    \label{fig:alt_linear_response}
\end{figure}

\paragraph{Fluxes and linear response.}
We now formally define transport coefficients, which measure the relative magnitude of a nonequilibrium flux with respect to the magnitude of the perturbation. The response, or flux, of the nonequilibrium system is measured by a scalar-valued observable~$R:\cX\to\R$, which we assume to vanish on average at equilibrium:
\begin{equation}\label{eq:response_function}
  \E_{\pi}\left[R\right]=0.
\end{equation}
For small forcing magnitudes, it is expected that the average response~$\E_{\pi_\eta}\left[R\right]$ is proportional to~$\eta$ (at least to leading order in~$\eta$), the proportionality constant being the transport coefficient~$\alpha$. The latter quantity is therefore defined as the derivative of the average flux with respect to the perturbation magnitude (provided the limit exists):
\begin{equation}\label{eq:linear_response}
  \alpha=\underset{\eta\to 0}{\lim}\frac{\E_{\pi_\eta}\left[R\right]}{\eta}.
\end{equation}
This definition motivates the following natural NEMD approach to estimate~$\alpha$:
\begin{itemize}
    \item{Pick a set of perturbation sample points~$\boldsymbol{\eta}=(\eta_k)_{1\leq k\leq K}$, and estimate the corresponding nonequilibrium steady state flux using trajectory averages of the nonequilibrium dynamics~\eqref{eq:nemd_dynamics}:
    \begin{equation}
      \forall 1\leq k\leq K,\qquad\widehat{\pi}_{\eta_k,T}(R) := \frac1{T}\int_{0}^T R(x_t^\eta)\,\d t,
    \end{equation}
    for~$T>0$ sufficiently large.
    }
    \item{Fit an assumed functional form~$\eta\mapsto \widehat{R}_{\boldsymbol{\eta},T}(\eta)$ to the data points~$\left\{\left(\eta_k,\widehat{\pi}_{\eta_k,T}(R)\right)\right\}_{1 \leq k \leq K}$ (e.g., a polynomial function), chosen so that~$\widehat{R}_{\boldsymbol{\eta},T}(0)=0$.}
    \item{The derivative
    \begin{equation}
       \widehat{\alpha}_{\boldsymbol{\eta},T} := \widehat{R}'_{\boldsymbol{\eta},T}(0) 
    \end{equation}
    is then used as an estimator of the transport coefficient~$\alpha$.
    }
\end{itemize}
A typical nonlinear response profile is depicted in Figure~\ref{fig:alt_linear_response}, together with the estimators associated to this fitting procedure. The target transport coefficient is the shear viscosity of liquid argon, at a temperature $T=85\,\mathrm{K}$ and a density~$1.41\,\mathrm{g/cm}^3$ using~$N=216$ atoms. Individual NEMD simulations were performed for a time~$T_{\mathrm{sim}}=100\,\mathrm{ps}$, using the BAOAB scheme~\eqref{eq:BAOAB} with a time step~$\Delta t = 2\,\mathrm{fs}$, and a friction coefficient~$\gamma = \sqrt{\epsilon m}/\sigma$ (where~$m$ is the mass of an Argon atom and~$(\sigma,\epsilon)$ are the parameters of the Lennard--Jones potential~\eqref{fig:LJ_potential}).

Noisy estimates (blue crosses) of the nonlinear response, obtained via NEMD simulations, are used to fit a quadratic model (red line). The linear response is estimated from the tangent to this model at $\eta=0$ (dotted black line) with the gray dotted lines indicating~$\pm 2$ standard errors in ordinary least squares regression. The dashed blue line corresponds to a reference linear response, obtained using~\eqref{eq:viscosity_equation} and experimental measurements of the viscosity provided by NIST.

The discrepancy between the computed and experimental responses can be attributed to a variety of error sources. Modelling errors are present, owing to the choices of friction coefficient~$\gamma$, regression model for the non-linear response, and empirical interaction potential. For this illustrative example, numerical parameters were also chosen to favor computational speed over accuracy: the number of particles~$N$, time step~$\Delta t$, potential cutoff parameter and simulation time would all be set to more conservative values in production, to reduce the associated biases and statistical error. The code used to perform these simulations is available in a public repository~\cite{github_cecam}.

\subsection{Connection with equilibrium fluctuations}
\label{subsec:linear_response_theory}

We briefly review linear response results giving alternative expressions for the coefficient~$\alpha$ in terms of equilibrium dynamical averages, and which form the basis of numerical methods such as the celebrated Green--Kubo formula~\cite{green1954,kubo1957,kubo1957b}. Here, we only give a somewhat informal presentation in the~$L^2(\pi)$ framework, but we stress that similar results can be proven for a broad class of systems and various functional settings, see for example~\cite{hairer_2010} or~\cite[Section~5.2]{lelievre2016}.

We present two derivations of the Green--Kubo formula. While the first one is conceptually simple, it is less general than the second, which in particular allows to cover the case of the nonequilibrium dynamics arising from perturbations of the diffusion coefficient as in~\eqref{eq:nemd_temperature_profile}.

\paragraph{Expansion of the nonequilibrium steady state.}
This derivation assumes that the nonequilibrium steady state~$\pi_\eta$ admits a probability density with respect to the equilibrium steady state~$\pi$, which can be perturbatively expanded (for sufficiently small~$|\eta|$) into a power series as
\begin{equation}
    \label{eq:nemd_steady_state_expansion}
    \pi_\eta=\pi\sum_{k=0}^{+\infty} \eta^k\psi_k,
\end{equation}
where~$\psi_k\in L^2(\pi)$ for all~$k\geq 0$.

Setting~$\eta=0$ in~\eqref{eq:nemd_steady_state_expansion},  it necessarily holds that~$\psi_0=\1_{\cX}$. The formal expansion~\eqref{eq:nemd_steady_state_expansion} can be shown to be valid when the nonequilibrium perturbation is sufficiently small, e.g., when~$\wcL$ is~$\cL_0$-bounded on~$L^2(\pi)$ and~$|\eta|$ is small enough; see for instance the proof of~\cite[Theorem~5.2]{lelievre2016}. This is the case for perturbations of the dynamics by a non-conservative force~$F$, under mild assumptions on~$F$ (see Table~\ref{tab:nemd_generator_perturbations} for the nomenclature).

Assuming the validity of such an ansatz, the stationary Fokker--Planck equation~\eqref{eq:fp_nemd} writes
\begin{equation}
    \label{eq:fp_equation_formal}
    \left(\cL_0+\eta\wcL\right)^\star\sum_{k=0}^{+\infty} \eta^k\psi_k=0,
\end{equation}
where adjoints are taken on~$L^2(\pi)$. Matching terms with the same powers in~$\eta$, it therefore holds
\begin{equation}
  \forall k\geq 1, \qquad -\cL_0^\star\psi_k = \wcL^\star\psi_{k-1},
\end{equation}
so that in particular~$-\cL_0^\star\psi_1 = \wcL^\star\1_{\cX}$. Provided the so-called conjugate flux
\begin{equation}
    \label{eq:conjugate_flux}
  S=\wcL^\star\1_{\cX}
\end{equation} 
belongs to the space~$L_0^2(\pi)$ of~$\pi$-centered observables, we can write
\begin{equation}\label{eq:psi_1}
    \psi_1=\left(-\cL_0^{-1}\right)^\star S,
\end{equation}
which in turn implies, by the definition~\eqref{eq:linear_response}, the following expression for the transport coefficient:
\begin{equation}\label{eq:linear_response_eq}
  \alpha = \underset{\eta\to 0}{\lim}\,\frac{1}{\eta} \int_{\cX}\left(\1_{\cX}-\eta(\cL_0^{-1})^\star S\right)R\,\d \pi = \int_{\cX}S(-\cL_0^{-1}R)\,\d\pi.
\end{equation}
Note that the expression for the conjugate response~$S$ can be computed explicitly by integration by parts. The expressions of the conjugate fluxes corresponding to the perturbations listed in Table~\ref{tab:nemd_generator_perturbations} are given in Table~\ref{tab:conjugate_response} below, and can easily be checked to belong to~$L_0^2(\pi)$ under mild assumptions on~$V,\delta T$ and~$F$.

The formulation~\eqref{eq:linear_response_eq} shows that the linear response~$\alpha$ can be expressed as the equilibrium average of the observable~$S(-\cL_0^{-1}R)$. Unfortunately, the various equilibrium sampling methods described in Section~\ref{sec:sampling_methods} cannot be applied outright, since they require the evaluation of the solution~$-\cL_0^{-1}R$ to a high-dimensional Poisson equation. Instead, one can reformulate~\eqref{eq:linear_response_eq} as a dynamical average using the expression of the inverse of~$\calL_0$ in~\eqref{eq:generator_inverse}. This leads to the celebrated Green--Kubo formula
\begin{equation}
  \label{eq:green_kubo}
  \alpha = \int_{0}^{+\infty}\E_{\pi}{\left[R(x_t)S(x_0)\right]}\,\d t,
\end{equation}
where~$\E_\pi$ denotes the expectation with respect to initial conditions~$x_0 \sim \pi$ and over all realizations of the stochastic dynamics at hand. Note that the Green--Kubo formula~\eqref{eq:green_kubo} allows for the computation of multiple transport coefficients from equilibrium trajectories of the dynamics~\eqref{eq:formal SDE} (by changing~$R$), which is convenient from a practical point of view.

\begin{table}[h]
    \centering
    \begin{tabular}{|c|c|c|}
        \hline
        \backslashbox{Perturbation}{Dynamics}& Overdamped Langevin&Underdamped Langevin\\
        \hline
        Non-conservative force & $\beta F\cdot \nabla V - \nabla\cdot F$ & $\beta F\cdot M^{-1}p$\\
        \hline
        Temperature profile & $\substack{\dps \Delta\delta T-2\beta\nabla V\cdot \nabla\delta T \\[3pt] \dps + \delta T\left(\beta^2|\nabla V|^2-\beta\Delta V\right)}$ & $\beta\gamma\delta T\left(\beta\left|M^{-1}p\right|^2-\mathrm{Tr}\,M^{-1}\right)$ \\
        \hline
    \end{tabular}
    \caption{Expressions for the conjugate response~$S$ for usual dynamics and perturbation types.}
    \label{tab:conjugate_response}
\end{table}

\paragraph{First-order expansion of the average flux.}
In some cases, in particular when the perturbation~$\widetilde\cL$ is not~$\cL_0$-bounded, the expansion~\eqref{eq:nemd_steady_state_expansion} is not valid. This is for instance the case with perturbations of the diffusion coefficients (see the last line of Table~\ref{tab:nemd_generator_perturbations}). In this case, one can still derive a finite-order expansion with respect to~$\eta$ for averages with respect to~$\pi_\eta$ of observables~$R \in L^2_0(\pi)$ whose derivatives satisfy certain integrability conditions. More precisely, one can show, under suitable assumptions, expansions of the form
\begin{equation}
  \label{eq:finite_eta_expansion}
  \frac{1}{\eta} \int_\cX R \, \rmd \pi_\eta = \alpha + \eta \int_\cX R \psi_2 \, \rmd\pi + \dots + \eta^{k-1} \int_\cX R \psi_k \, \rmd\pi + \eta^k r_{R,\eta},
\end{equation}
where~$r_{R,\eta}$ is a remainder term uniformly bounded in~$\eta$ for sufficiently small~$\eta$, and the functions~$(\psi_j)_{2\leq j\leq k}$ are the ones previously defined. The first-order expansion is in fact sufficient to obtain the Green--Kubo formula.

The argument to establish this formula is formally the following. We start from the equality
\[
\int_{\mathcal X} \left( \mathcal L_0 + \eta \widetilde {\mathcal L} \right) \varphi \, \d \pi_{\eta} = 0.
\]
Upon replacing~$\varphi$ by~$\mathcal L_0^{-1} \Pi R$, this leads to
\[
\int_{\mathcal X} R \, \d \pi_{\eta} = \int_{\mathcal X} R \, \d \pi + \eta \int_{\mathcal X} \widetilde {\mathcal L} \mathcal L_0^{-1} \Pi R \, \d \pi_{\eta}.
\]
We next iterate the reasoning for the second term term on the right hand side:
\begin{align*}
    \int_{\mathcal X} R \, \d \pi_{\eta} 
    &= \int_{\mathcal X} R \, \d \pi + \eta \int_{\mathcal X} \widetilde {\mathcal L} \mathcal L_0^{-1} \Pi R \, \d \pi + \eta^2 \int_{\mathcal X} (\widetilde {\mathcal L} \mathcal L_0^{-1} \Pi)^2 R \, \d \pi_{\eta} \,  \\
    &= \dots \\
    &= \int_{\mathcal X} R \, \d \pi + \sum_{j=1}^{k-1} \eta^j \int_{\mathcal X} (\widetilde {\mathcal L} \mathcal L_0^{-1} \Pi)^j R \, \d \pi + \eta^{k} \int_{\mathcal X} (\widetilde {\mathcal L} \mathcal L_0^{-1} \Pi)^{k+1} R \, \d \pi_{\eta}.
\end{align*}
Making this argument rigorous requires sufficient integrability properties in order to control the error terms, which involve derivatives of both~$R$ and of various solutions to Poisson equations associated to~$\cL_0$. More precisely, the response~$R$ should belong to a class of observables~$\mathcal S_0 \subset L_0^2(\pi)$ which is stable by~$\cL_0^{-1}$ and~$\Pi \widetilde{\cL}$; and moreover~$\pi_\eta$ should integrate elements of~$\cS_0$. Typically, one takes $\mathcal S_0$ to be the class of smooth observables whose derivatives grow at most polynomially at infinity, see~\cite[Section~3.3 and Remark~5.5]{lelievre2016} for further details.

\subsection{Examples}
\label{subsec:transport_coefficient_examples}

We provide examples of nonequilibrium systems that are used in practice to compute physically relevant transport coefficients. We focus the presentation on the prototypical potentials discussed in Section~\ref{sec:example_systems}, and assume that the system evolves according to the underdamped Langevin dynamics~\eqref{eq:langevin} in a periodic configurational domain~$\cD = (L\mathbb T)^d$. The first two examples are concerned with Lennard-Jones fluids, and the last one with atom chains.

\paragraph{Mobility of a Lennard-Jones particle.} Arguably the simplest example of transport coefficient is provided by the mobility of a particle moving in a fluid. Physically, this quantity measures how easily mass is transported through the fluid in response to an external driving field. It is closely related to the self-diffusion coefficient, obtained by the Einstein relation (see~\eqref{eq:einstein_relation} below and~\cite{rodenhausen1989}).

For simplicity, we consider the case in which all particles are identical, so that the mass matrix in~\eqref{eq:hamiltonian} is given by~$M=\mathrm{Diag}(m,\dots,m)$, where~$m>0$ is the mass of a single particle.
In the framework described in Section~\ref{subsec:nemd_dynamics}, and keeping the same notation, this corresponds to taking a constant forcing~$F\in\R^d$ with~$|F|=1$, and measuring the average particle flux through the hyperplane perpendicular to~$F$ (which is the velocity in the direction~$F$). In other words, we choose respectively
\begin{equation}\label{eq:pert_resp_LJ_mobility}
  R(q,p) = m^{-1}F^\top p, \qquad \wcL = F^\top\nabla_p,
\end{equation}
for the response observable and nonequilibrium perturbation of the generator, so that the conjugate response is
\begin{equation}
  S(q,p) = m^{-1}\beta F^\top p.
\end{equation}
The mobility~$\alpha_F$ is defined as the linear response~\eqref{eq:linear_response} associated with the choice~\eqref{eq:pert_resp_LJ_mobility}. In view of the Green--Kubo formula~\eqref{eq:green_kubo}, it can also be written as
\begin{equation}
  \alpha_F = \beta F^\top \mathfrak{C}F,\qquad\mathfrak{C} := m^{-2}\int_0^{+\infty} \E_{\mu}\left[p_0p_t^\top\right] \rmd t \in \R^{d\times d}
\end{equation}
in terms of the velocity autocovariance matrix~$\mathfrak{C}$. Note that by isotropy, $\alpha_F$ does not depend on the particular choice of direction~$F$. The common value is called the mobility.

The mobility can be related to the self-diffusion coefficient~$\mathfrak{D}$ entering in Fick's law, defined from the mean squared displacement based on the displacement
\begin{equation}\label{eq:unperiodized_displacement}
  Q_t - Q_0 = \int_0^t m^{-1}p_s \, \rmd s
\end{equation}
corresponding to unperiodized evolutions of particles. More precisely,
\begin{equation}
  \mathfrak{D} = \lim_{T \to +\infty} \frac{\dps \E_\mu\left[\left|Q_T-Q_0\right|^2\right]}{2d T},
\end{equation}
which corresponds to the Einstein relation~\cite{rodenhausen1989} (see also equation~\eqref{eq:ein_diff} below). In fact, 
\begin{equation}\label{eq:einstein_relation}
  \mathfrak{D} = \frac1d \mathrm{Tr} \, \mathfrak{C} = \mathfrak{C}_{11},
\end{equation}
so that, in particular, $\alpha_F = \beta \mathfrak{D}$ for~$|F|=1$. 

\paragraph{Shear viscosity of a Lennard-Jones fluid.}
The second prototypical example of transport coefficient is the shear viscosity of a Newtonian fluid. Here, we present the computation of the shear viscosity in monoatomic Lennard-Jones fluids, following the method described in~\cite{js12}, inspired by the sinusoidal transverse force (STF) method~\cite{gms73}. Another class of NEMD algorithms to measure the shear viscosity consists of boundary-driven methods, such as the direct simulation of Couette flows via shearing boundary conditions~\cite{le72}; see~\cite[Section 6.3]{evans_morriss_2007} and~\cite[Section 9.3]{todd_daivis_2017} for a presentation of these methods, and~\cite{DG23} for their mathematical analysis.

The STF method proceeds by analogy with Newton's macroscopic law of viscosity, for a fluid subjected to a shear force~$f$ directed along the longitudinal $x$-coordinate, which varies in intensity in the transverse $y$-coordinate. At the continuum level, the shear viscosity~$\xi$ is defined via the constitutive relation
\begin{equation}\label{eq:newton_equation}
  \sigma_{xy} = -\xi\frac{\d u_x}{\d y},
\end{equation}
where~$\sigma_{xy}$ is the~$(x,y)$ component of the local stress tensor, and~$u_x$ is the local $x$-velocity field of the fluid. Both~$\sigma_{xy}$ and~$u_x$ are functions of the $y$ position in the fluid.

\begin{figure}
    \centering
    \includegraphics[width=0.7\textwidth]{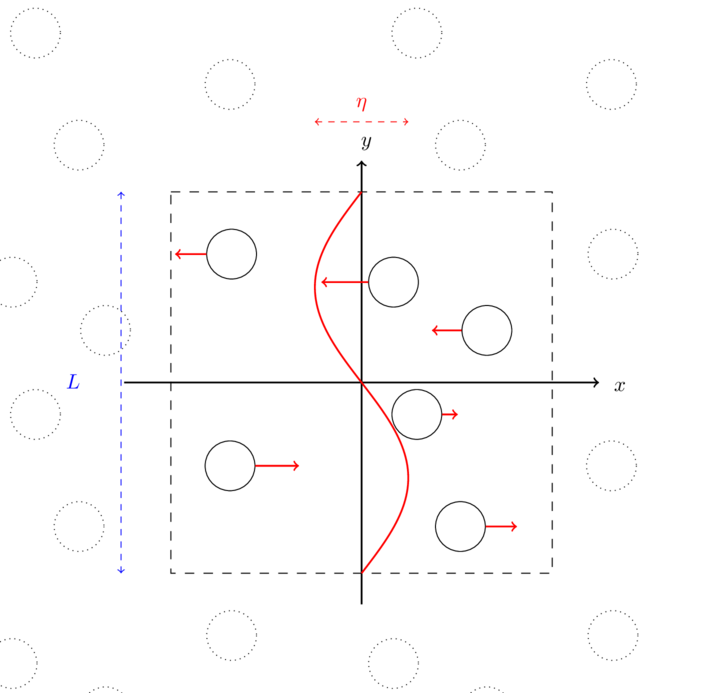}
    \caption{Schematic representation of the sinusoidal transverse field method in a two-dimensional fluid. The black dashed square is the unit cell of the one-particle configurational domain~$L\mathbb{T}^{2}$. Particles, represented by black circles, are subjected to an external forcing field in the longitudinal direction, represented by red arrows. The field is a function of a periodic transverse profile, plotted in the solid red line, whose amplitude is proportional to the forcing parameter~$\eta$. Periodic images of the system are represented by dotted circles. An example trajectory of the empirical Fourier response to this forcing is represented in Figure~\ref{fig:nemd_response}.}
    \label{fig:stf}
\end{figure}

Microscopically, the action of the shear forcing on the fluid particles is defined by the following non-conservative force field:
\begin{equation}\label{eq:stf_forcing}
  \forall 1\leq j\leq N,\qquad F(q)_{j,x} = f(q_{j,y}),\qquad F(q)_{j,y}=F(q)_{j,z}=0.
\end{equation}
The forcing field~\eqref{eq:stf_forcing} acts on each component of the $x$-momentum variable, in a way that is dictated by the corresponding component of the $y$-position variable, according to a fixed forcing profile~$f$. The STF method derives its name from the standard choice $f(y) = \sin(2\pi y/L)$, although other profiles can also be considered. A schematic representation of the transverse forcing procedure is provided in Figure~\ref{fig:stf}.

The microscopic formulation of the relation~\eqref{eq:newton_equation} relies on appropriate definitions of the velocity profile~$u_x$ and shear-stress profile~$\sigma_{xy}$. These are defined as the linear responses~\eqref{eq:linear_response} for (a limit of)~$y$-dependent observables defined in~\cite{js12} following a mathematically rigorous version of the Irving--Kirkwood procedure~\cite{irving_kirkwood_1950}. Using the linear response formula~\eqref{eq:linear_response_eq}, one can show that the velocity and shear-stress linear response profiles~$u_x,\sigma_{xy}$ and forcing profile~$f$ are related as (see~\cite[Proposition 1]{js12})
\begin{equation}
  \frac{1}{\rho}\frac{\d\sigma_{xy}(y)}{\d y} + \gamma u_x(y) = f(y),
\end{equation}
where we recall that~$\gamma>0$ is the friction parameter in the underdamped Langevin dynamics~\eqref{eq:langevin}, and~$\rho=N/L^3$ is the particle density.

Formally substituting~$\sigma_{xy}$ from the Newton relation~\eqref{eq:newton_equation}, one arrives at the following differential equation for~$u_x$:
\begin{equation}
  -\frac{\xi}{\rho}\frac{\d^2 u_x(y)}{\d y^2} + \gamma u_x(y) = f(y),
\end{equation}
from which, passing to Fourier series, one can express the viscosity~$\xi$ as
\begin{equation}\label{eq:viscosity_equation}
  \xi = \rho\left(\frac{\mathfrak{f}_1}{\mathfrak{u}_1}-\gamma\right)\left(\frac{L}{2\pi}\right)^2,
\end{equation}
where~$\mathfrak{f}_1$ and $\mathfrak{u}_1$ are the first Fourier coefficients in~$y$ of the forcing profile~$f$ and the velocity linear response profile~$u_x$, respectively, on the periodic one-dimensional torus $L\mathbb T$. Given that~$\mathfrak{f}_1$ is analytically known, the only unknown quantity on the right-hand side of~\eqref{eq:viscosity_equation} is the Fourier coefficient~$\mathfrak{u}_1$, but one can formally show that it is a transport coefficient in its own right, associated with the response function corresponding to the ``empirical Fourier flux''
\begin{equation}\label{eq:fourier_flux}
  R(q,p) = \frac1{Nm}\sum_{j=1}^N p_{j,x}\exp\left(\frac{2\mathrm{i}\pi q_{j,y}}{L}\right).
\end{equation}
This linear response can be estimated from nonequilibrium trajectory averages, or using the Green--Kubo formula with the conjugate flux
\begin{equation}
  S(q,p) = \frac{\beta}{m} F(q)^\top p. 
\end{equation}

Various other approaches, relying on discretized estimates of~$\sigma_{xy},u_x$ or transformations thereof, via NEMD or equilibrium-fluctuation formulas, and based on the constitutive equation~\eqref{eq:newton_equation}, are of course possible.

\paragraph{Thermal transport in atom chains.}
\label{ss:thermtransp}

Let us consider the atom chains described in Section~\ref{sec:example_systems}, with free boundary conditions. The system evolves according to a purely Hamiltonian dynamics in the bulk, with temperatures maintained to~$T_{\rm L} = T + \Delta T /2$ on the left end and~$T_{\rm R} = T - \Delta T /2$ on the right end through Ornstein--Uhlenbeck processes on the momenta with friction intensities~$\gamma_{\rm L}$ and~$\gamma_{\rm R}$. More precisely, the equations of motion of the dynamics read
\begin{equation}\label{eq:dyn_eqns}
  \left\{ \begin{aligned}
    \dd q_i & = \frac{p_i}{m_i} \, \d t,& \qquad i &\in \{1, \dots, n\}, \\
    \dd p_i & = \left[ v'(q_{i+1} - q_i) - v'(q_i-q_{i-1}) \right] \dd t, &\qquad i &\in \{2, \dots, n-1\}, \\
    \dd p_1 & =  v'(q_2 - q_1) \, \dd t - \gamma_\rmL p_1 \, \dd t + \sqrt{2\gamma_\rmL T_\rmL} \, \dd W_1(t),&&\\
    \dd p_n & =  -v'(q_n - q_{n-1}) \, \dd t- \gamma_\rmR p_n \, \dd t + \sqrt{2\gamma_\rmR T_\rmR} \, \dd W_n(t),&&
  \end{aligned} \right.
\end{equation}
where~$\{W_i(t)\}_{i\in \{1,n\}}$ are independent standard Brownian motions. In the absence of thermal forcing (i.e.,~$T_\rmL=T_\rmR=T$),~\eqref{eq:dyn_eqns} describes the dynamics of a system in equilibrium at temperature~$T$. The infinitesimal generator of the process~\eqref{eq:dyn_eqns} can be written as
\begin{equation}\label{eq:L_gen}
  \cL = \cA + \gamma_\rmL \mathcal{S}_{T_\rmL} + \gamma_\rmR \mathcal{S}_{T_\rmR},
\end{equation}
where $\cA$ is the generator of the Hamiltonian dynamics, and $\cS_{T_\rmL}$ and $\cS_{T_\rmR}$ are respectively the generators of the Ornstein--Uhlenbeck processes at the boundaries: 
\begin{align}
  \cA & = \sum_{i=1}^{n+1} \left[\frac{p_i}{m_i} \partial_{q_i} + \left( v'\left(q_{i+1} - q_i\right) - v'\left(q_i - q_{i-1}\right) \right) \partial_{p_i}\right], \label{eq:A_gen}\\ 
  \cS_{T_\rmL} & = -p_1 \partial_{p_1} + T_\rmL \partial_{p_1}^2, \qquad \cS_{T_\rmR} = -p_n \partial_{p_n} + T_\rmR \partial_{p_n}^2.\label{eq:S_gen}
\end{align}
In fact, in order to identify linear response properties, we write, following~\eqref{eq:nemd_generator},
\[
  \cL = \cL_0 + \frac{\Delta T}{2} \widetilde{\cL},
  \qquad
  \cL_0 = \cA + T(\gamma_\rmL \partial_{p_1}^\star\partial_{p_1} + \gamma_\rmR \partial_{p_n}^\star\partial_{p_n}), \qquad \widetilde{\cL} = \gamma_\rmL \partial_{p_1}^2 - \gamma_\rmR \partial_{p_n}^2,
\]
where adjoints are taken on~$L^2(\mu)$ with~$\mu$ the canonical probability measure~\eqref{eq:canonical measure} at temperature~$T$ associated with the Hamiltonian~\eqref{eq:H(q,p)_general}. Note that the dynamics~\eqref{eq:dyn_eqns} locally conserves the energy in the bulk. This conservation is expressed via the following homogeneous continuity equation:
\[
\frac{\dd}{\dd t} e_i(t) = j_{i-1,i}(t) - j_{i,i+1}(t), \qquad i \in \{2,\dots n-1\},
\]
where the expression of $e_i$ is given in~\eqref{eq:ei} and
\begin{equation}\label{eq:instcurr_bulk}
  j_{i,i+1}(t) = -\frac{p_i(t) + p_{i+1}(t)}2\,v'(q_{i+1}(t) - q_i(t)), \qquad i \in\{1,\dots,n-1\},
\end{equation}
is the instantaneous energy current flowing from the $i$-th to the $(i+1)$-th atom.

Under suitable conditions on~$v$ (see~\cite{Carmona2007}) the dynamics reaches a steady state, denoted by~$\mu_{T_\rmL,T_\rmR}$, characterized by a non-zero stationary current of energy flowing through the system. In view of Fourier's law
\begin{equation}\label{eq:Fourier}
  j(x) = -\kappa(T(x)) \nabla T(x),
\end{equation}
energy flows from the hotter to the colder thermostat and its density is related to the temperature gradient by a positive constant $\kappa$, called the thermal conductivity. The computation of the thermal conductivity for a chain of $n$~atoms in the linear response regime can be performed by means of the two possible strategies presented in Section~\ref{sec:transport_coefficients}, namely the NEMD approach (Section~\ref{subsec:nemd_dynamics}) and the Green--Kubo approach (Section~\ref{subsec:linear_response_theory}). Introducing the total energy current in the bulk
\begin{equation}\label{eq:atom_chain_observable}
  J_n(q,p) = \sum_{i=1}^{n-1} j_{i,i+1}(q,p),
\end{equation}
the thermal conductivity of an atom chains of length~$n$ is expressed following~\eqref{eq:linear_response} as
\begin{equation}\label{eq:thermcondNEMD}
  \kappa_n = \lim_{\Delta T \to 0} \,\frac{\mathbb E_{\mu_{T_\rmL,T_\rmR}}[J_n]}{\Delta T}, \qquad \Delta T = T_\rmL-T_\rmR.
\end{equation}
This definition is motivated by an analogy with Fourier's law. Indeed, in the stationary state, the value of~$\E_{\mu_{T_\rmL,T_\rmR}}[j_{i,i+1}]$ does not depend on~$i \in \{1,\dots,n-1\}$, and the spatially averaged energy current~$\E_{\mu_{T_\rmL,T_\rmR}}[J_n]/(n-1)$ is also equal to this value. Moreover, assuming a linear temperature profile across the chain (which is the case for sufficiently small temperature differences), the quantity $(T_\rmL-T_\rmR)/(n-1)$ gives the temperature gradient. The thermal conductivity is finally defined as the opposite of the ratio of the energy current and the temperature gradient, in accordance with~\eqref{eq:Fourier}.

For the equilibrium computation of $\kappa$ via the Green--Kubo formula~\eqref{eq:green_kubo}, one chooses for the response function the energy current, namely~$R(q,p) = J_n(q,p)$. Moreover, a straightforward computation of the conjugate response~\eqref{eq:conjugate_flux} gives $S(q,p) = (\gamma_\rmL p_1^2 - \gamma_\rmR p_n^2)/T^2 - (\gamma_\rmL-\gamma_\rmR)/T$. Hence,
\begin{equation}\label{eq:thercond_GK}
  \kappa_n = \int_0^{+\infty} \mathbb{E}_\mu
  \left[ J_n(q_t,p_t) S(q_0,p_0) \right] \dd t = \frac{1}{(n-1)T^2}\int_0^{+\infty} \mathbb{E}_\mu  \left[J_n(q_t,p_t)J_n(q_0,p_0)\right] \dd t.
\end{equation}
The final equality is nontrivial; we refer the interested reader to~\cite{KDN09} for instance for details of the computation. Fourier's law is satisfied if $\kappa_n$ has a finite limit as the system size goes to infinity (the so-called thermodynamic limit), that is
\begin{equation}
  \label{eq:thermo_limit_1D}
  \lim_{n \to +\infty} \kappa_n  = \kappa < +\infty.
\end{equation}
While for common three-dimensional solids, numerical simulations have shown that the thermal conductivity has a well defined thermodynamic limit, the validity of Fourier's law in one-dimensional systems remains unsettled, see~e.g., \cite{BLR2000,LLP2003,Dhar2008,LLP2016}. One-dimensional atom chains are fundamental models for studying whether the thermal conductivity~$\kappa_n$ remains finite as the system size increases. In purely harmonic chains (i.e., when $v(r)$ is quadratic in~$r$), energy travels ballistically so that~$\kappa \sim n$; while in chains with anharmonic interaction potentials~$v$, the transport is typically anomalous ($\kappa \sim n^{\alpha}$ for some~$0<\alpha<1$), rotor chains being the unique nonlinear model currently known for which the thermal conductivity has a well defined thermodynamic limit. Normal diffusion ($\alpha=0$) is recovered only when the total momentum preservation is broken, either by the addition of an on-site pinning potential or a bulk noise to the dynamics (see for instance~\cite{LLP2016}). In any case, simulations of long chains are needed to decide whether the thermal conductivity has a well defined large system size limit.

\section{Nonequilibrium molecular dynamics}
\label{sec:NEMD}

We focus in this section on nonequilibrium methods to compute transport coefficients. We start by quantifying numerical errors in the standard NEMD approach in Section~\ref{sec:NEMD_error_analysis}. We next consider various ways to reduce the dominant error, namely the statistical error:
\begin{itemize}
\item \textit{synthetic forcings}, considered in Section~\ref{subsec:synth_for}, can potentially increase the range of the linear response, thereby allowing for larger forcing magnitudes and leading to a better signal to noise ratio in the computation of the response of the system under an external forcing;
\item \textit{control variates}, discussed in Section~\ref{sec:noneq_control_variates}, for which two cases are considered. 

\textit{Coupling methods} use averages of the response function under the equilibrium dynamics as a dynamical control variate to reduce the statistical error as long as the equilibrium and nonequilibrium trajectories remain sufficiently close (Section~\ref{subsec:static_control_variates});

\textit{Static control variates} are constructed by replacing the response function~$R$ by~$R - \mathcal{L}_\eta \Phi$, with~$\Phi$ being chosen in order to reduce the variance of the related estimator for the transport coefficient (Section~\ref{subsec:NEMD_coupling});

\item \textit{Norton dynamics}, considered in Section~\ref{sec:NEMD_Norton}, relies on a dual viewpoint: instead of fixing the magnitude of the forcing and measuring the average response, one fixes the response of the system and measures the average forcing needed to maintain it.
\end{itemize}

\subsection{Error analysis}
\label{sec:NEMD_error_analysis}

The principle of NEMD is to approximate the limit~$\eta \to 0$ in~\eqref{eq:linear_response} by considering a small but finite value of~$\eta$ and replacing the expectation with respect to the steady-state probability measure by a time average. More precisely, for the generic nonequilibrium dynamics~\eqref{eq:nemd_dynamics}, an estimator of the linear response is, for an observable~$R$ with average zero under the equilibrium measure:
\begin{equation}\label{Stoltz:eq:as_limit_NEMD}
  \widehat{\alpha}_{\eta,t} = \frac{1}{\eta t}\int_0^t R(x_s^\eta) \, \rmd s \xrightarrow[t\to+\infty]{\mathrm{a.s.}} \alpha_\eta := \frac1\eta \int_\mathcal{X} R \, \rmd \pi_\eta = \alpha + \mathrm{O}(\eta).
\end{equation}
The various sources of error for the estimator~$\widehat{\alpha}_{\eta,t}$, made precise below, are the following: 
\begin{enumerate}[(i)]
\item A statistical error with asymptotic variance $\mathrm{O}(t^{-1}\eta^{-2})$, much larger than the usual asymptotic variance of order~$1/t$ associated with standard time averages which are not divided by a factor~$\eta$ (see Section~\ref{T_finite_NEMD}).
\item A bias of order~$\mathrm{O}(\eta)$ due to the fact that~$\eta \neq 0$, as made apparent on the right-hand side of~\eqref{Stoltz:eq:as_limit_NEMD}. This bias is directly obtained from the results of Section~\ref{sec:transport_coefficients}, see for instance~\eqref{eq:finite_eta_expansion}.
\item A bias arising from the finiteness of the integration time~$t$ in the estimator~$\widehat{\alpha}_{\eta,t}$ (see Section~\ref{T_finite_NEMD}).
\item A bias arising from the discretization in time when implementing nonequilibrium dynamics in computer simulations (see Section~\ref{sec:dt_bias_NEMD}).
\end{enumerate}
Let us emphasize that there is a balance between taking~$\eta$ sufficiently small in order to limit the bias~$\alpha_\eta - \alpha = \mathrm{O}(\eta)$, and~$\eta$ not being too small so that the asymptotic variance controlling the magnitude of the statistical error is not too large. In practice, the largest source of error usually is the statistical error.

\subsubsection{Analysis of the variance and the finite integration time bias}
\label{T_finite_NEMD}

Let us start by stating the two important results concerning errors related to the finiteness of the integration time. First, the statistical error is dictated by a central limit theorem: it is shown below that
\begin{equation}\label{Stoltz:eq:CLT_NEMD}
  \sqrt{t} \left(\widehat{\alpha}_{\eta,t} - \alpha_\eta \right) \xrightarrow[t \to +\infty]{\mathrm{law}} \mathcal{N}\left(0,\frac{v_{R,\eta}}{\eta^2}\right),
\end{equation}
where
\begin{equation}
v_{R,\eta} = v_{R,0} + \mathrm{O}(\eta),
\qquad
v_{R,0} = 2  \int_0^{+\infty} \mathbb{E}_0\left[R(x_t) R(x_0)\right] \rmd t.
\end{equation}
This quantifies the fact that the statistical error~$\widehat{\alpha}_{\eta,t} - \alpha_\eta$ is of order~$1/(\eta \sqrt{t})$. Note also that the asymptotic variance agrees, at dominant order in~$\eta$, with the one of the equilibrium dynamics. This result shows that long simulation times of order~$t \sim \eta^{-2}$ are required in order to have an asymptotic variance of order~1. Balancing the bias of order~$\eta$ arising from~$\alpha_\eta - \alpha = \mathrm{O}(\eta)$ and the statistical error requires~$t \sim \eta^{-3}$. Concerning the finite time integration bias, the following estimate holds:
\begin{equation}\label{Stoltz:eq:FT_bias}
  \left| \mathbb{E}\left(\widehat{\alpha}_{\eta,t}\right) - \alpha_\eta \right| \leq \frac{K}{\eta t}.
\end{equation}
The bias due to $t < +\infty$ is therefore of order~$\mathrm{O}\left(t^{-1}\eta^{-1}\right)$, typically smaller than the statistical error.

The key tool for proving these results is the Poisson equation
\begin{equation}\label{eq:Poisson_eq_eta}
  -\left(\mathcal{L}_0+\eta\widetilde{\mathcal{L}}\right) \mathscr{R}_\eta = R - \int_\mathcal{E} R \, \rmd \pi_\eta.
\end{equation}
A simple computation based on It\^o calculus gives
\begin{equation}\label{Stoltz:eq:A_eta_t_Ito}
    \widehat{\alpha}_{\eta,t} - \frac1\eta \int_\cX R \, \rmd \pi_\eta = \frac{\mathscr{R}_\eta(x_0^\eta) - \mathscr{R}_\eta(x_t^\eta)}{\eta t} + \frac{\sqrt{2\gamma}}{\eta t\sqrt{\beta}} \int_0^t \nabla \mathscr{R}_\eta(x_s^\eta)^\top \sigma_\eta(x_s^\eta) \, \rmd W_s.
\end{equation}
The limit~\eqref{Stoltz:eq:CLT_NEMD} then follows from a central limit theorem for martingales, while~\eqref{Stoltz:eq:FT_bias} is obtained by taking expectations in the above equality~\eqref{Stoltz:eq:A_eta_t_Ito}. The prefactor of~$\eta^{-2}$ in the expression of the asymptotic variance in~\eqref{Stoltz:eq:CLT_NEMD} turns out to be
\begin{equation}
  v_{R,\eta} = 2 = \int_\cX \nabla\mathscr{R}_\eta^\top \sigma_\eta\sigma_\eta^\top \nabla\mathscr{R}_\eta\, \rmd \pi_\eta = 2\int_\cX R \mathscr{R}_\eta\, \rmd \pi_\eta,
\end{equation}
which can be expanded in powers of~$\eta$ by the same techniques as the ones used to prove~\eqref{eq:finite_eta_expansion} (expanding both the steady state probability measure and the solution to the Poisson equation~\eqref{eq:Poisson_eq_eta}). At leading order one retrieves the asymptotic variance~$v_{R,0}$ for time averages of~$R$ under the equilibrium dynamics. Of course, some care has to be taken in order to make these computations rigorous, as one needs to ensure that the solution~$\mathscr{R}_\eta$ to the Poisson is sufficiently regular, and has good integrability properties; see~\cite{spacek2023} for details.

\subsubsection{Analysis of the time step discretization bias}
\label{sec:dt_bias_NEMD}

In practice, one needs to discretize the nonequilibrium dynamics, similarly to what is done in Section~\ref{sec:sampling_methods} for equilibrium dynamics. Error estimates similar to the ones considered in Section~\ref{sec:error_dt_equilibrium} can be derived to make precise the bias on the estimated transport coefficients.

To simplify the presentation, we consider overdamped Langevin like dynamics perturbed by a non-gradient drift~$\eta F(q)$, and discretized with a Euler--Maruyama scheme:
\begin{equation}\label{eq:noneq_ovd}
  q^{n+1} = q^n - \Delta t (\nabla V(q^n) + \eta F(q^n)) + \sqrt{\frac{2\Delta t}{\beta}} G^n,
\end{equation}
where~$G^n$ is a standard Gaussian vector. The results established in~\cite[Appendix~C]{darshan2024} provide error estimates à la Talay--Tubaro~\cite{TT90}. More precisely, denoting by~${\nu}_{\eta,{\Delta t}}$ the invariant probability measure of~\eqref{eq:noneq_ovd}, which is well defined for compact position spaces,
\begin{equation}
  \int_\mathcal{D} R \, \rmd{\nu}_{\eta,{\Delta t}} = \int_\mathcal{D} R \Big(1+ \eta f_{0,1} + {\Delta t} f_{1,0} + \eta \Delta t f_{1,1} \Big) \, \rmd{\nu} + r_{R,\eta,{\Delta t}},
\end{equation}
with~$f_{0,1}$ being the perturbation~\eqref{eq:psi_1} of the invariant measure arising from the nonequilibrium forcing, $f_{1,0}$ being the one arising from the time step discretization as in~\eqref{eq:error_discretised_pi}, and where the remainder is compatible with linear response:
\begin{equation}
  \left|r_{R,\eta,{\Delta t}}\right| \leq K(\eta^2 + {\Delta t}^2), 
  \qquad 
  \left|r_{R,\eta,{\Delta t}} - r_{R,0,{\Delta t}}\right| \leq K \eta (\eta + {\Delta t}^2).
\end{equation}
A corollary of this equality is the following error estimate on the numerically computed linear response where, importantly, the linear response is computed from averages with respect to the invariant probability measure of the discretized equilibrium dynamics:
\begin{equation}
\begin{aligned}
\alpha_{{\Delta t}} & = \lim_{\eta \to 0} \frac{1}{\eta} \left(\int_\mathcal{X} R(q) \, \nu_{\eta,{\Delta t}}(\rmd{q}) - \int_\mathcal{X} R(q) \, \nu_{0,{\Delta t}}(\rmd{q}) \right) \\
& = \alpha + \Delta t \int_\mathcal{X} R(q) f_{1,1}(q) \, \nu(\rmd q) + \mathrm{O}(\Delta t^2).
\end{aligned}
\end{equation}
Similar results can be obtained for Langevin dynamics, with results that hold uniformly in the overdamped limit of~$\gamma \to +\infty$; see~\cite{leimkuhler2016computation} and~\cite[Remark~4]{darshan2024}.

\subsection{Synthetic forcings}
\label{subsec:synth_for}

In order to decrease the statistical error, which scales as~$1/\eta$ (as made precise in~\eqref{Stoltz:eq:CLT_NEMD}), one idea is to extend the regime of linear response (that is, the range of values of~$\eta$ for which the nonlinear part of the response~$\E_{\pi_\eta}[R]$ is small), so that larger values of~$\eta$ can be considered for a fixed level of bias. An approach for doing so is to rely on synthetic forcings, an idea initially suggested in~\cite{evans_morriss_2007} and explored more systematically in~\cite{spacek2023}. The key idea behind synthetic forcings is that there are infinitely many forcings which lead to the same transport coefficient. This flexibility can be exploited to develop better numerical methods for the computation of transport coefficients. 

A practical way of making use of synthetic forcings is to add an extra forcing to the physical perturbation of the system, as long as this extra forcing preserves the linear part~$\psi_1$ of the response of the system, given by~\eqref{eq:psi_1}, which determines the value of the transport coefficient we want to estimate; see~\eqref{eq:linear_response_eq}. This amounts to considering a new perturbed dynamics with generator 
\begin{equation}
  \cL_{\eta,a} = \cL_0 + \eta\left( \wcL + a \cL_{\rm{extra}} \right), 
\end{equation}
where $\wcL$ is the generator of the physical perturbation, as given in~\eqref{eq:nemd_generator}, while the extra forcing is chosen in order for its generator~$\cL_{\rm{extra}}$ to satisfy
\begin{equation}\label{eq:synthetic_condition}
  \cL_{\rm{extra}}^\star \mathbbm{1} = 0.
\end{equation}
The condition~\eqref{eq:synthetic_condition} ensures, in view of~\eqref{eq:conjugate_flux}, that the conjugate response~$S$ is unchanged, therefore the linear response of the steady state~$\psi_1$ is not altered in view of~\eqref{eq:psi_1}. We call the resulting perturbation $\wcL + a\cL_{\rm{extra}}$ a synthetic forcing, as it has no physical meaning; it is a mathematical device used to reduce variance in the estimator of the transport coefficient. 

In practice, some extra perturbations that might be added to the perturbed dynamics~\eqref{eq:nemd_dynamics} and that can be easily implemented in Monte Carlo simulations are: 
\begin{enumerate}[(1)]
\item \emph{First-order differential operators} $\cL_{\rm{extra}} = G^\top \nabla_x$, with $G\colon\mathcal{X} \to \R^d$ such that $\nabla\cdot(G\psi_0) = 0$. The latter condition ensures that~\eqref{eq:synthetic_condition} is satisfied since~$\cL_{\rm{extra}}^\star = -G^\top \nabla_x$. This amounts to modifying the drift in~\eqref{eq:nemd_dynamics} by adding a term~$a\eta G$ to it. 
\item \emph{Second-order differential operators} of the form $-\partial_{x_i}^\star\partial_{x_i}$ or more generally
  \begin{equation}
    \cL_{\rm{extra}} = \sum_{i,j=1}^d -\partial_{x_j}^\star D_{ij}\partial_{x_i}
  \end{equation}
  for some symmetric positive semidefinite function $D_{ij}\colon\mathcal{X}\to\R$. Such operators are self-adjoint, i.e.,~$\cL_{\rm{extra}} = \cL_{\rm{extra}}^\star$, so that~\eqref{eq:synthetic_condition} is easily seen to hold. In practice, resorting to such extra forcings amounts to modifying the strength of the fluctuation-dissipation in~\eqref{eq:nemd_dynamics}. 
\end{enumerate}

\begin{figure}
    \centering
    \includegraphics[scale=0.5]{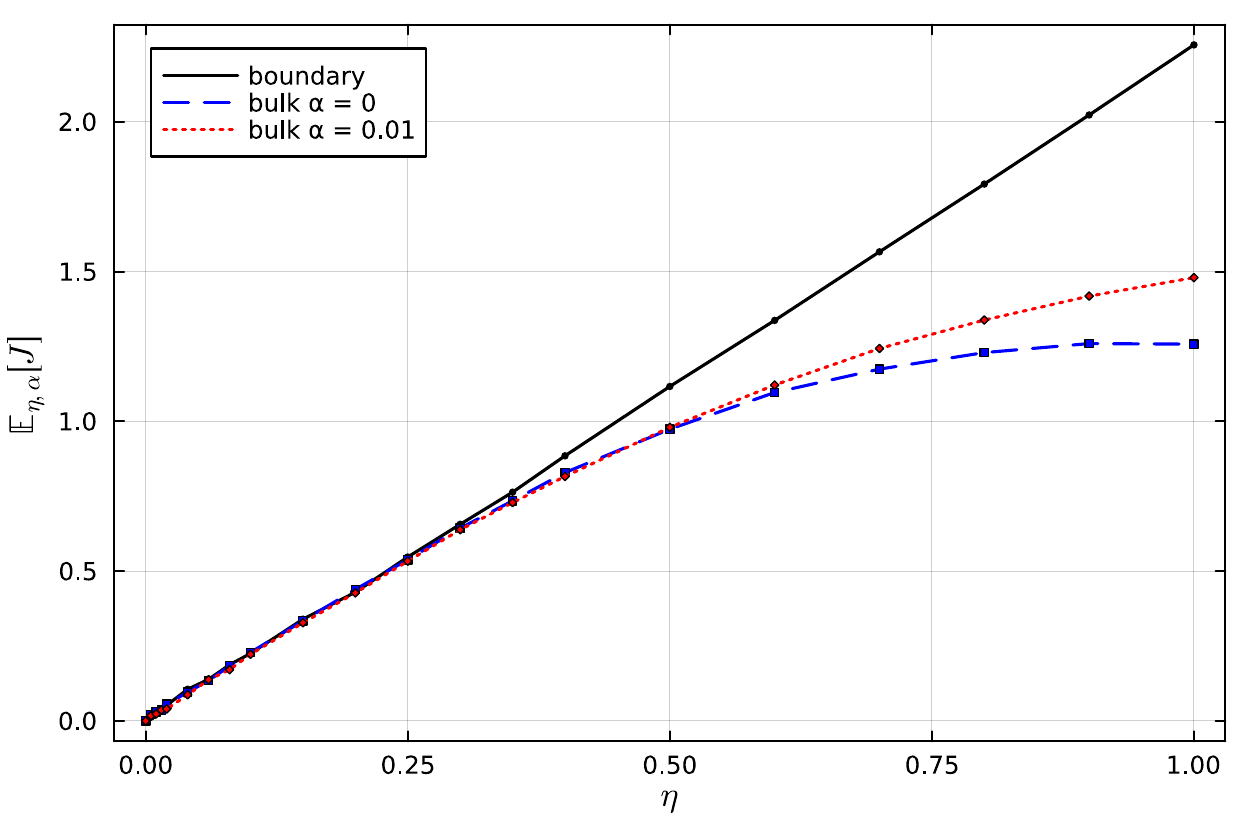}
    \caption{Response of the energy current for a chain of rotors obtained with different forcings: boundary forcing (black solid curve); bulk forcing (blue dashed curve); combination of bulk perturbation and thermostats on each atom (red dotted curve). } 
    \label{fig:synthetic_forcing}
\end{figure}

\paragraph{A concrete example: the rotor chain. }
Figure~\ref{fig:synthetic_forcing} illustrates the concept of synthetic forcing in the context of a chain of rotors (see Figure~\ref{fig:rotor}). Specifically, we consider a chain of~$n$ rotors obeying the reference dynamics~\eqref{eq:dyn_eqns}  
with potential $v(r)$ given in \eqref{eq:V_rotor} and reference temperature $T_\rmL=T_\rmR = T$. We choose the following perturbations of this reference dynamics:
\begin{itemize}
    \item \textbf{Boundary forcing}. We perturb the reference dynamics by imposing a thermal gradient induced by boundary thermostats at different temperatures~$T_\rmL \neq T_\rmR$, producing the black solid curve in Figure~\ref{fig:synthetic_forcing}. The associated linear response is given by~\eqref{eq:thercond_GK}.
    
    \item \textbf{Bulk forcing}. Drawing on~\cite{zhang_ibister_evans}, we simulate the following system, obtained by perturbing the vector field of~\eqref{eq:dyn_eqns}: 
    \begin{equation}\label{eq:dyn_eqns:rotors:bulk_drive}
    \left\{ \begin{aligned}
    \dd r_i & = (p_i - p_{i-1}) \, \d t,& \qquad i &\in \{1, \dots, n\}, \\
    \dd p_i & = \left[ \left(1-\frac{\eta}{n-1}\right)\, v'(r_{i+1}) - \left(1+\frac{\eta}{n-1}\right)\,v'(r_i) \right] \dd t, &\qquad i &\in \{2, \dots, n-1\}, \\
    \dd p_1 & = \left(1-\frac{\eta}{n-1}\right)\, v'(r_2) \, \dd t - \gamma_\rmL p_1 \, \dd t + \sqrt{2\gamma_\rmL T} \, \dd W_1(t),&&\\
    \dd p_n & =  - \left(1+\frac{\eta}{n-1}\right)\,v'(r_n) \, \dd t- \gamma_\rmR p_n \, \dd t + \sqrt{2\gamma_\rmR T} \, \dd W_n(t).&&
    \end{aligned} \right.
    \end{equation}
    The generator associated with this bulk perturbation is
    \begin{equation}
    \label{eq:bulk_forcing_chain}
    \wcL = -\frac{1}{N-1}\sum_{i=2}^N v'(r_i) (\partial_{p_i} + \partial_{p_{i-1}}). 
    \end{equation}
    One can show, from \eqref{eq:green_kubo}, that the associated linear response of the energy current is 
    \begin{equation}\label{eq:alpha_bulk_driven}
	\alpha_{n,\rm{bulk}} =  \frac{ 2}{T^2 (N-1)} \int_{0}^{+\infty}\mathbbm{E}_{\pi}\left[R(t) R(0)\right]\dd t.
    \end{equation}
    Note that $\alpha_{n,\rm{bulk}} = 2 \kappa_n$, the latter quantity being defined in~\eqref{eq:thercond_GK}. We therefore plot with the black solid curve in Figure~\ref{fig:synthetic_forcing} the linear response of the energy current~$J$ defined in~\eqref{eq:atom_chain_observable}, multiplied by~2.
    
    \item \textbf{Combined bulk and thermostat forcing}. We consider the combination of the previous bulk forcing~\eqref{eq:bulk_forcing_chain} with an extra forcing with generator 
    \[
    \cL_{\rm{extra}} = - T\sum_{i=2}^{N-1} \partial^\ast_{p_i}\partial_{p_i} = \sum_{i=2}^N (- p_i \partial_{p_i} + T\partial^\ast_{p_i}).
    \]
    This corresponds to coupling each bulk atom in~\eqref{eq:dyn_eqns:rotors:bulk_drive} with thermostats at temperature~$T$ and coupling strength~$a\eta$, and yields the red dotted curve in Figure~\ref{fig:synthetic_forcing}. 
\end{itemize} 
All three perturbations produce the same linear response in the system, although their linear response regimes differ. The simulations were performed for systems of size~$n=20$, with~$\gamma_\rmL=\gamma_\rmR = 1$ and~$T_\rmL=T_\rmR =T= 0.3$. The numerical scheme is a splitting scheme where the Hamiltonian part of the dynamics is integrated using the Störmer-Verlet method, while the Ornstein-Uhlenbeck part is integrated analytically, with a timestep~$\Delta t = 0.01$. 

\subsection{Control variates}
\label{sec:noneq_control_variates}

The principle of the control variate method is to subtract off a process \(\left(\xi_t\right)_{t\geq 0}\) with mean zero at stationarity from the estimator \(\left(\widehat{\alpha}_{\eta, t}\right)_{t\geq 0}\) defined in~\eqref{Stoltz:eq:as_limit_NEMD} to create a new estimator
\begin{equation}
	\widehat{A}_{\eta, t} = \widehat{\alpha}_{\eta, t} - \xi_t,
\end{equation}
which has lower (asymptotic) variance compared to the original estimator and the same expectation. These methods are not specific to the synthetic forcings described in the previous section and can be easily combined with said methods. In this section, we consider two control variate methods where the process \(\left(\xi_t\right)_{t\geq 0}\) is given either by
\[\xi_t = - \int_0^t \mathcal{L}_\eta\Phi\left(x_s^\eta\right)\rmd s,\]
where the function \(\Phi\) is an approximation of the solution to the Poisson equation \eqref{eq:Poisson_eq_eta}, or by
\[\xi_t = \int_0^t R\left(y_s^0\right)\rmd s,\]
where \(\left(y_t^0\right)_{t\geq 0}\) is the equilibrium dynamics, i.e. \eqref{eq:nemd_dynamics} with \(\eta = 0\), coupled to the perturbed dynamics. For simplicity, we present these two control variate methods in continuous time. In practice, one has to consider appropriate discretizations for implementations. 

\subsubsection{Static control variates}\label{subsec:static_control_variates}
To build a suitable control variate one may use the same simulation trajectory and subtract off \(\phi\left(x^\eta_s\right)\) from~$R(x^\eta_s)$ in the time averages such as the ones appearing in \eqref{Stoltz:eq:as_limit_NEMD}, where \(\phi\) is some well-chosen function with
expectation zero with respect to \(\pi_\eta\). This corresponds to the estimator
\begin{equation}
\label{eq:NEMD_static_CV}
    \widehat{A}_{\eta, t}^\phi = \frac{1}{\eta t}\int_0^t \left[R\left(x_s^\eta\right) - \phi\left(x_s^\eta\right)\right]\rmd s.
\end{equation}
To effectively reduce the asymptotic variance, the so--called zero--variance principle~\cite{AC1999} says the optimal choice is
\begin{equation}\label{eq:optimal_static_control_variate}
	\phi_{\mathrm{ZVP}, \eta} := -\cL_\eta \mathscr{R} = R - \pi_\eta(R),
\end{equation}
with \(\mathscr{R}_\eta\) the solution of \eqref{eq:Poisson_eq_eta}. This observable indeed allows to achieve zero variance with respect to \(\pi_\eta\) since
\[
\mathrm{Var}_{\pi_\eta}\left(R - \phi_{\mathrm{ZVP}, \eta}\right) = \int_{\cX} \left(R - \phi_{\mathrm{ZVP}, \eta}\right)^2 \rmd \pi_\eta - \left(\int_{\cX} \left[R - \phi_{\mathrm{ZVP}, \eta}\right] \rmd \pi_\eta\right)^2 = 0,
\]
hence the name of the method. Approaches based on this principle were initially proposed in computational statistics literature~\cite{andradottir_1993, Henderson1997, henderson_2002, MSI2013} and the name was coined in the molecular simulation literature~\cite{AC1999}. 

While the estimator~\eqref{eq:NEMD_static_CV} with $\phi = \phi_{\rm ZVP}$ has zero variance, actually constructing such the function~$\phi_{\mathrm{ZVP}, \eta}$ is infeasible in practice since~$\pi_\eta(R)$ is unknown (this is the very quantity one wants to estimate). Building an approximation of the control variate in~\eqref{eq:optimal_static_control_variate} requires to solve the very high dimensional Poisson equation~\eqref{eq:Poisson_eq_eta}. In \cite{RS19}, the authors suggest a perturbative approach in which one builds an approximation to~\eqref{eq:optimal_static_control_variate} by solving a surrogate Poisson equation
\begin{equation}\label{eq:surrogate_poisson_equation}
	\widehat{\cL} \Phi = R - \widehat{\pi}(R),
\end{equation}
where~$\widehat{\cL}$ is the generator of a surrogate dynamics, and~\(\widehat{\pi}\) is the invariant probability measure associated with~\(\widehat{\cL}\); and then applies the generator of the original dynamics to this solution in order to obtain a function which has by construction average~0 with respect to~$\pi_\eta$:
\[\phi = \cL_\eta \Phi.\]
The operator \(\widehat{\cL}\) in \eqref{eq:surrogate_poisson_equation} is typically the generator of a simplified dynamics, for instance the equilibrium dynamics, or some approximate dynamics where the interacting potential is replaced by its quadractic approximation; see \cite[Section~3.2 \&~4.2]{RS19} for precise examples. The reduction of the variance depends on how close the solution to the surrogate Poisson equation~\eqref{eq:surrogate_poisson_equation} is to the true solution~\eqref{eq:optimal_static_control_variate}, and also on the accuracy of the numerical solution; see \cite[Theorems~3 \&~4]{RS19} for quantitative statements.

\subsubsection{Coupling based control variates}\label{subsec:NEMD_coupling}
Since a transport coefficient quantifies the difference between the non--equilibrium and equilibrium dynamics, one is tempted to include the equilibrium dynamics \(\left(y_s^0\right)_{s\geq 0}\) in the estimator:
\begin{equation}\label{eq:coupled_estimator}
	\widehat{A}_{\eta, t} = \frac{1}{\eta t} \int_0^t \left[R\left(x_s^\eta\right) - R\left(y_s^0\right)\right] \rmd s.
\end{equation}
Since \(\pi_0\left(R\right) = 0\), the addition of the extra term does not change the asymptotic limit of the estimator. However, if the two dynamics are independent, the resulting estimator \eqref{eq:coupled_estimator} has roughly twice the variance of the standard estimator \eqref{Stoltz:eq:as_limit_NEMD}. To reduce the variance, one should \emph{couple} the two dynamics. A coupling of the two dynamics is a process \(\left(z^\eta_s\right)_{s\geq 0} = \left(x_s^\eta, y_s^0\right)_{s \geq 0}\) on \(\cX \times \cX\) whose first component has the law of the non--equilibrium dynamics and whose second component has the law of the equilibrium dynamics. The process can be represented as the solution to
\begin{equation}\label{eq:coupled_dynamics}
	\begin{aligned}
		\rmd x_s^\eta &= b_\eta\left(x_s^\eta\right)\rmd s + \sigma_\eta\left(x_s^\eta\right) \rmd W_s,\\
		\rmd y_s^0 &= b_0\left(y_s^0\right)\rmd s + \sigma_0\left(y_s^0\right) \rmd \widetilde{W}_s,
	\end{aligned}
\end{equation}
where \(\left(W_s, \widetilde{W}_s\right)_{s \geq 0}\) is a coupling of two standard Brownian motions. At least one coupling always exists: the independent coupling where \(W\) and \(\widetilde{W}\) in \eqref{eq:coupled_dynamics} are independent Brownian motions. However, as mentioned earlier, the resulting estimator \eqref{eq:coupled_estimator} has higher variance than the standard one.

For the coupling based estimator to have significantly lower variance compared to the standard estimator, the distance between the coupled dynamics should be of order \(\eta\) for long times. To simplify the presentation, we suppose that the diffusion matrix \(\sigma_\eta\) is constant and equal to the identity matrix times a constant that does not depend on \(\eta\). A natural first coupling to try is the synchronous coupling where the dynamics are coupled by choosing the same Brownian motion, i.e. \(\widetilde{W} = W\) in \eqref{eq:coupled_dynamics}. This coupling is straightforward to implement in practice since it amounts to using the same sequence of random numbers to simulate the two dynamics, which can easily be done in many MD packages by using the same random seed for two simulations -- one with~$\eta=0$ and one with~$\eta > 0$. If the drift is strongly contractive everywhere, i.e. there exists \(m > 0\) such that
\begin{equation}\label{eq:global_strong_contractivity}
	\bigl \langle x - y, b_0(x) - b_0(y)\bigr\rangle \leq -m \left|x - y\right|^2, \qquad \forall x, y \in \cX,
\end{equation}
and the convergence in \eqref{eq:linear_forcing} holds, for example, uniformly in \(x\in \cX\), then one can show that the variance of the resulting coupling based estimator is bounded as \(\eta \downarrow 0\), see for example~\cite[Section~3]{darshan2024} for details. For drifts of the form \eqref{eq:non_conservative_force}, sufficient conditions are that the potential~\(V\) is~\(m\)--strongly convex and the non--conservative forcing being bounded. The key assumption is the convexity of \(V\) since synchronous coupling uses only the drift to bring the trajectories together. If \(V\) is multimodal, the synchronous coupling will not push the trajectories together when the non--equilibrium and equilibrium dynamics are in different modes. For non--convex potentials more intricate couplings, such as sticky coupling, can be considered; see~\cite[Section~4]{darshan2024}.

\subsection{Norton dynamics}
\label{sec:NEMD_Norton}

We present in this section the stochastic Norton dynamics~\cite{blassel2024}, a recent extension to the setting of NEMD stochastic dynamics of the ideas of Evans and Morriss which were formulated in the deterministic setting, see~\cite{evans1983,hoover1983,evans1985,evans1986,evans1993} and~\cite[Section 6.7]{evans_morriss_2007}. The method is based on the observation that, although physical intuition suggests that the steady state response is~\textit{caused} by the nonequilibrium driving force (this is indeed reflected in the choice of terminology), there is no reason to suppose this is the case from the macroscopic point of view. Instead, the flux and the forcing co-occur, and one could equivalently think of the nonequilibrium flux as creating a resisting force opposing the drive out of equilibrium.

\begin{figure}[!h]
    \centering
    \begin{subfigure}{0.6\textwidth}
        \includegraphics[width=\textwidth]{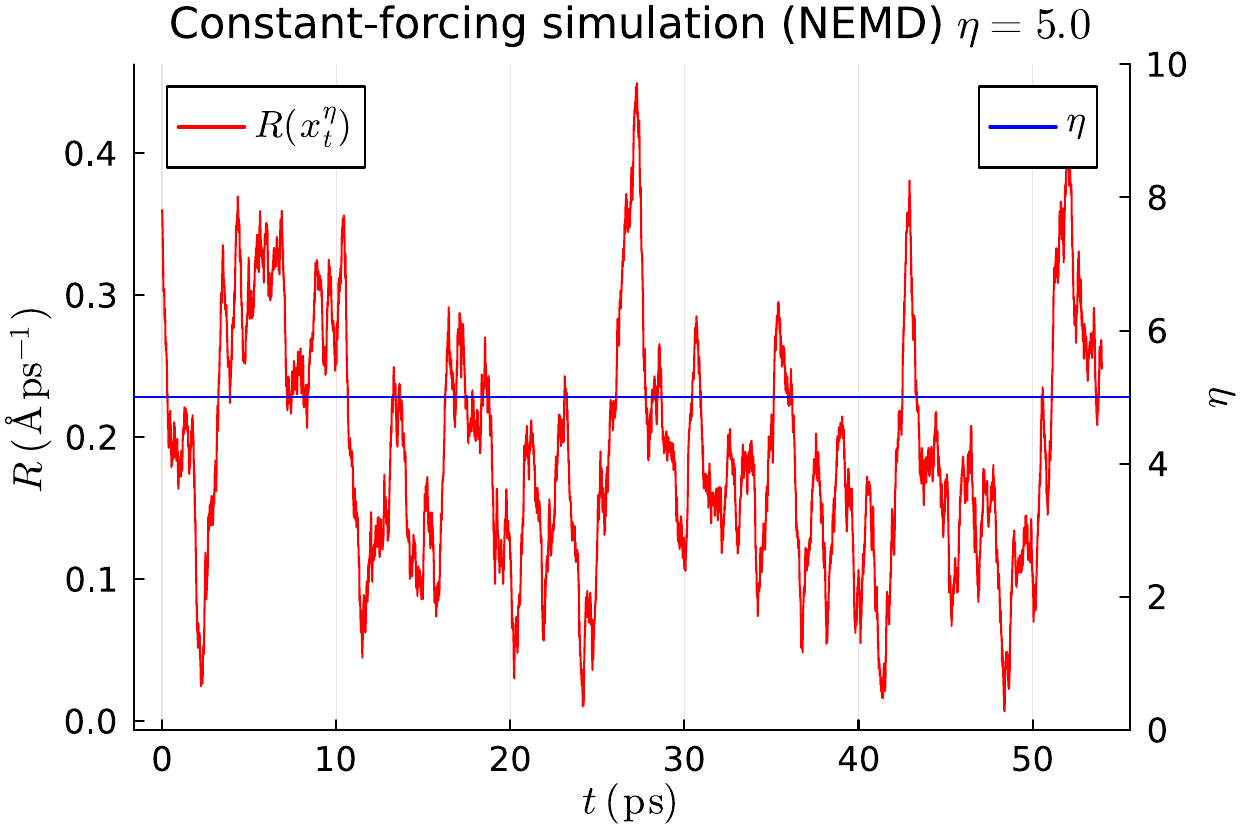}
        \caption{Sample trajectory of the empirical Fourier flux~\eqref{eq:fourier_flux} under the shear forcing~\eqref{eq:stf_forcing}.}
        \label{fig:nemd_response}
    \end{subfigure}
    \hfill
    \begin{subfigure}{0.6\textwidth}
        \includegraphics[width=\textwidth]{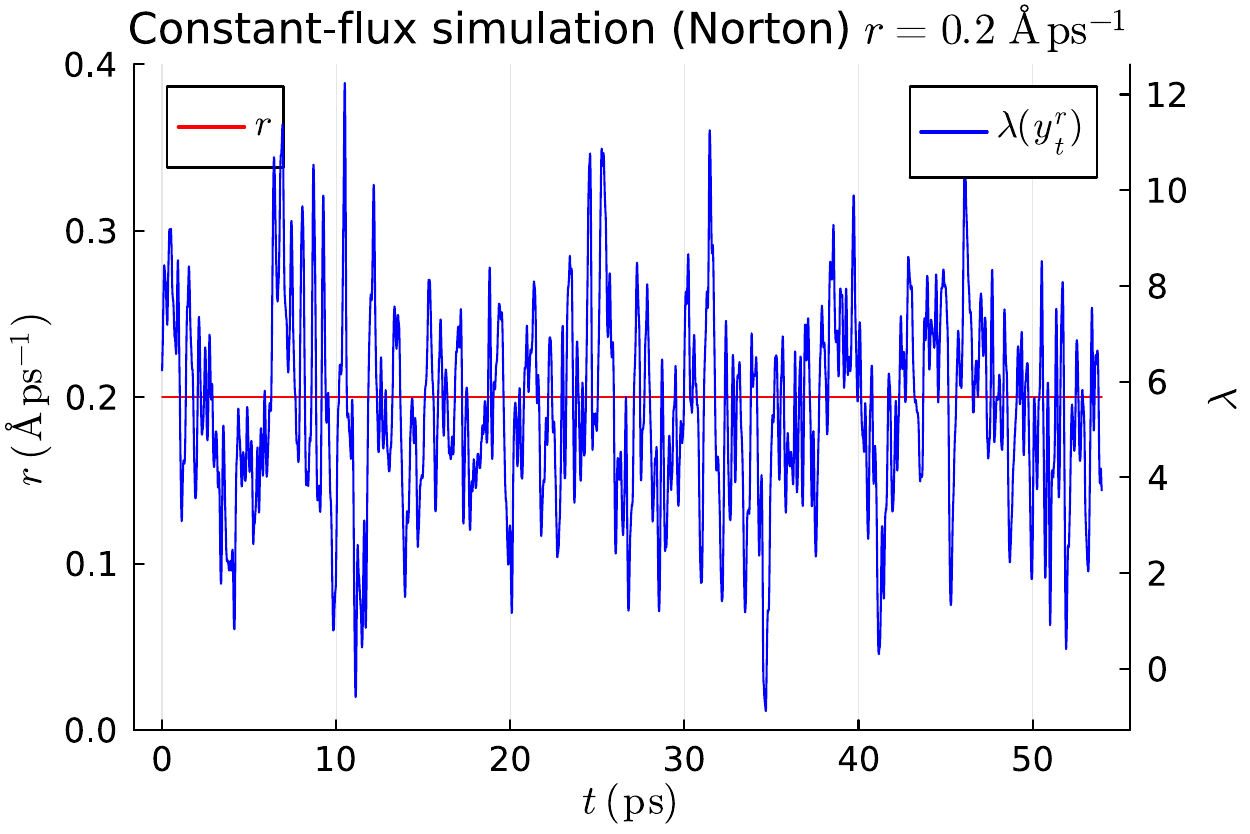}
        \caption{Sample trajectory of the forcing~\eqref{eq:norton_forcing_expr} in the corresponding Norton ensemble.}
        \label{fig:norton_lambda}
    \end{subfigure}
    \caption{Illustration of the Norton-NEMD duality for the liquid Argon system considered in Figure~\ref{fig:alt_linear_response}. Both nonequilibrium simulations are carried out at roughly equivalent state points~$(\eta,r)=(5.0,0.2)$, obtained from the fitted non-linear response. Both simulations were performed in the same physical and numerical conditions as those used in Figure~\ref{fig:alt_linear_response}. The code used to obtain these results is available in a public repository~\cite{github_norton}.}
    \label{fig:norton}
\end{figure}

This shift in perspective is at the core of the Norton method: instead of fixing the magnitude of the nonequilibrium perturbation, and measuring the average flux in the nonequilibrium steady state, as in standard NEMD methods, the Norton approach~\textit{fixes the flux}, and measures the~\textit{average forcing magnitude} needed to sustain it in a corresponding constant-flux ensemble. The Norton and NEMD nonequilibrium ensembles are in a dual relationship with respect to the pair of flux/forcing variables: one of these values is fixed exactly, while the other one is fixed on average by the nonequilibrium steady state. This duality is analogous to the relationship between the microcanonical (where the energy, i.e., the Hamiltonian of the system, is conserved) and canonical (where the temperature is conserved) equilibrium ensembles with respect to the energy/temperature pair of thermodynamic variables. A representative example of this duality is given in Figure~\ref{fig:norton}, where the NEMD and Norton methods are applied to a shear viscosity computation (see Section~\ref{subsec:transport_coefficient_examples}).

We formulate the Norton dynamics for NEMD perturbations of the drift, such as those presented in Section~\ref{subsec:transport_coefficient_examples} for the computation of the mobility and shear viscosity, i.e.,~$b_\eta = b + \eta F$ for some~$F:\cX\to\R^d$ and~$\sigma_\eta = \sigma$ in~\eqref{eq:nemd_dynamics}. The dynamics assumes the same basic form as the NEMD dynamics, except that the perturbation parameter~$\eta$ is replaced by the increment of a stochastic process:
\begin{equation}\label{eq:norton_dynamics}
  \d y_t^r = b(y_t^r) \, \rmd t + \sigma(y_t^r) \, \rmd W_t + F(y_t^r) \, \rmd \Lambda_t^r,
\end{equation}
where the forcing magnitude process $\Lambda_t^r$ is determined by the constant-flux constraint
\begin{equation}\label{eq:norton_constraint}
  R(y_t^r) = R(y_0^r) = r,\qquad\forall\,t\geq 0,
\end{equation}
and the constant parameter~$r>0$ fixes the value of the flux. A simple computation shows that these conditions define~$\Lambda^r$ as an adapted Itô process:
\begin{equation}
  \Lambda_t^r = \int_0^t\lambda(y_s^r)\,\d s + \int_0^t\widetilde{\lambda}(y_s^r)\,\d W_s,
\end{equation}
for explicit observables~$\lambda:\mathcal X\to\R$ and~$\widetilde\lambda:\mathcal X\to\R^{1\times d}$. In particular,~$\lambda$ is given by
\begin{equation}\label{eq:norton_forcing_expr}
  \lambda = -\frac1{F^\top \nabla R}\left(b^\top \nabla R +\frac{1}{2} \mathrm{Tr}\left[\sigma^\top \overline{P}_{F,\nabla R}^{\top}\nabla^2 R \overline{P}_{ F,\nabla R}\sigma\right]\right),\qquad \overline{P}_{F,\nabla R} = \mathrm{Id} - \frac{F\nabla R^\top}{F^\top\nabla R}.
\end{equation}
Note that this observable is only defined when~$F^\top\nabla R\neq 0$. This condition is a controllability constraint, which should be satisfied along the trajectories of~$y^r$.

As a result, one may measure the average forcing magnitude via a trajectory average of the dynamics~\eqref{eq:norton_dynamics}:
\begin{equation}\label{eq:norton_trajectory_average}
  \frac{1}{T}\int_0^T \lambda(y_t^r)\,\d t \xrightarrow [T\to+\infty]{} \E_{\pi^r}{\left[\lambda\right]},
\end{equation}
where~$\pi^r$ denotes the Norton steady state, defined as the invariant probability distribution of the dynamics~\eqref{eq:norton_dynamics}. The Norton dynamics can also be defined for multiple and time-dependent flux constraints, see~\cite{blassel2024}. Note that the trajectory average on the left-hand side of~\eqref{eq:norton_trajectory_average} neglects the contribution of the observable~$\widetilde{\lambda}$, corresponding to the martingale component of~$\Lambda^r$. This amounts to using the predictable process~$-\frac1T\int_0^T\widetilde{\lambda}(Y_t^r)\,\d W_t$ as a control variate in the estimation of the average forcing magnitude. In this sense, the stochastic Norton method incorporates by construction some variance reduction in its estimation of the nonequilibrium response.

By analogy with the macroscopic transport law, the Norton linear response is defined as the reciprocal derivative of the forcing with respect to the flux (provided it is well-defined):
\begin{equation}\label{eq:norton_tp}
  \alpha^*=\underset{r\to 0}{\lim} \frac{r}{\E_{\pi^r}\left[\lambda\right]}.
\end{equation}

At this stage, many foundational questions~(see~\cite[Section 7]{blassel2024}) regarding the Norton method remain open, including the existence and uniqueness of~$\pi^r$, the ergodicity of the Norton dynamics, and sufficient conditions to ensure the equivalence of linear responses~(${\alpha=\alpha^*}$), or even nonlinear responses. More generally, the method prompts the question of the equivalence of nonequilibrium ensembles, which is only expected to hold in specific, but nevertheless practically relevant, limiting regimes. The upcoming work~\cite{darshan2025} resolves some of these issues in a mean-field setting.

It has however been demonstrated numerically in~\cite{blassel2024} that the Norton dynamics does provide an equivalent measure of the mobility and shear viscosity of Lennard-Jones fluids, and that this equivalence also extends to the nonlinear regime. For the case of the shear viscosity, the fixed current was the empirical Fourier flux~\eqref{eq:fourier_flux} described in Section~\ref{subsec:transport_coefficient_examples}. Similar results have recently been obtained for Norton dynamics applied to nonequilibrium DPD~\cite{wu2025}.

Another point is that the fluctuation behavior of the forcing process~$\lambda(y_t^r)$ is often quite different from that of response process~$R(x_t^\eta)$ from NEMD. Roughly speaking, the former tends to display large instantaneous fluctuations over short correlation times, while in NEMD, smaller fluctuations occur, which are correlated over much longer times. This phenomenon has also been observed in recent simulations using deterministic Norton dynamics with Nos\'e--Hoover thermostats~\cite{tee2023,sasaki2025}.

More generally, the idea of choosing from equivalent statistical ensembles which possess the most desirable fluctuation behavior for a target property is classical in equilibrium statistical mechanics~\cite{lebowitz1967,cancrini2017}. The Norton dynamics aims at extending this possibility to the nonequilibrium setting by proposing an equivalent nonequilibrium ensemble to sample from, with different fluctuation properties. Overall, variance reduction for constant-flux ensemble estimators of various transport coefficients has indeed been observed in~\cite{evans1983,tee2023,blassel2024,sasaki2025,wu2025}. Furthermore, anomalous scaling properties of the asymptotic variance of Norton estimators in the thermodynamic limit of~$N\to+\infty$ have been reported in~\cite{blassel2024} (see also~\cite{wu2025}), hinting that the Norton method may be particularly useful in larger systems as the associated asymptotic variance would be much smaller than for standard NEMD approaches.

\section{Equilibrium fluctuation formulas}
\label{sec:GK_and_co}

We focus in this section on equilibrium methods to compute transport coefficients, relying on time correlations. We start by quantifying numerical errors in the numerical implementation of the standard Green--Kubo formula in Section~\ref{subsec:error_green_kubo}. We next consider various ways to reduce the dominant error, namely the statistical error:
\begin{itemize}
\item by resorting to \textit{control variates}, see Section~\ref{sec:control_gk};
\item by relying on \textit{importance sampling in path space} through Girsanov reweighting, see Section~\ref{sec:Girsanov};
\item by making use of \textit{alternative fluctuation formulas}, see Section~\ref{sec:alternative_fluctuation}.
\end{itemize}

\subsection{Error analysis for Green--Kubo formulas}
\label{subsec:error_green_kubo}

We present in this section how to compute transport coefficients based on the Green--Kubo formula~\eqref{eq:green_kubo}:
\begin{equation}\label{eq:transport_gk}
  \alpha = \int_{0}^{+\infty}\E_\pi{\left[R(x_t)S(x_0)\right]} \, \d t,
\end{equation}
where the expectation is over all initial conditions~$x_0 \sim \pi$ and over all realizations of the equilibrium dynamics~\eqref{eq:formal SDE}. Recall that, even if this expression only refers to the equilibrium dynamics, the Green--Kubo formula still implicitly depends on the nonequilibrium forcing~\eqref{eq:non_conservative_force} through the conjugate response~$S$ defined in~\eqref{eq:conjugate_flux}.

\begin{figure}
    \centering

    \begin{subfigure}[t]{0.49\columnwidth}
        \includegraphics[width=\columnwidth]{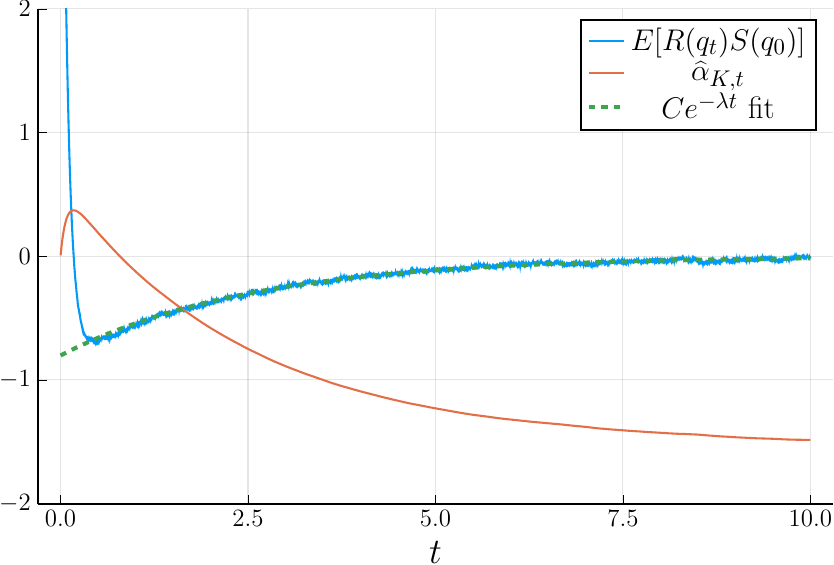}
        \caption{Autocorrelation function associated with the Green--Kubo formula, together with its exponential asymptotic decay of the form~\eqref{eq:exp_decay_semigroup} obtained by least-squares fitting.}
        \label{fig:autocorrelation}
    \end{subfigure}\hfill
    \begin{subfigure}[t]{0.49\columnwidth}
        \includegraphics[width=\columnwidth]{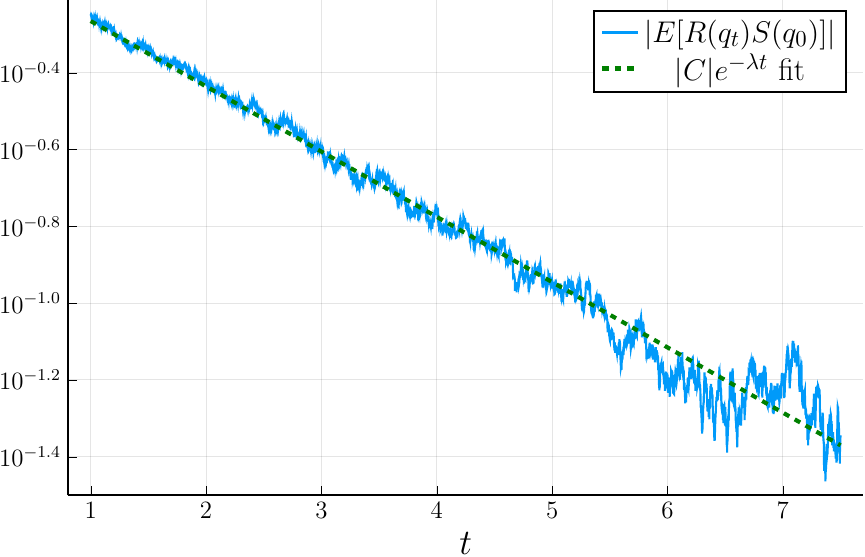}
        \caption{Log-scale representation of the autocorrelation function and its exponential decay between $t = 1$ and $t=7.5$.}
        \label{fig:log_fit_autocorrelation}
    \end{subfigure}

    \begin{subfigure}[t]{0.49\columnwidth}
        \includegraphics[width=\columnwidth]{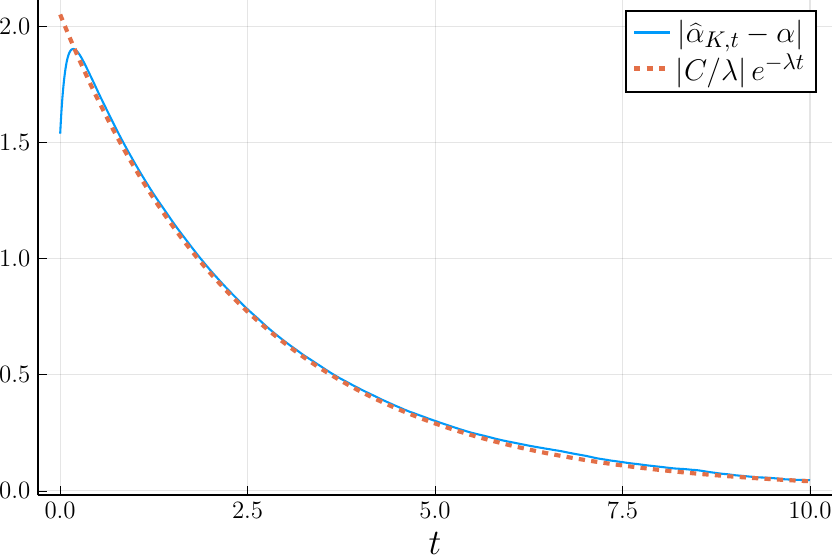}
        \caption{Truncation bias~\eqref{eq:truncation_bias} of the Green--Kubo formula. The decay rate is inferred from the fitted exponential behavior shown in~\cref{fig:autocorrelation}, without additional fitting.}
        \label{fig:truncation_bias}
    \end{subfigure}\hfill
    \begin{subfigure}[t]{0.49\columnwidth}
        \includegraphics[width=\columnwidth]{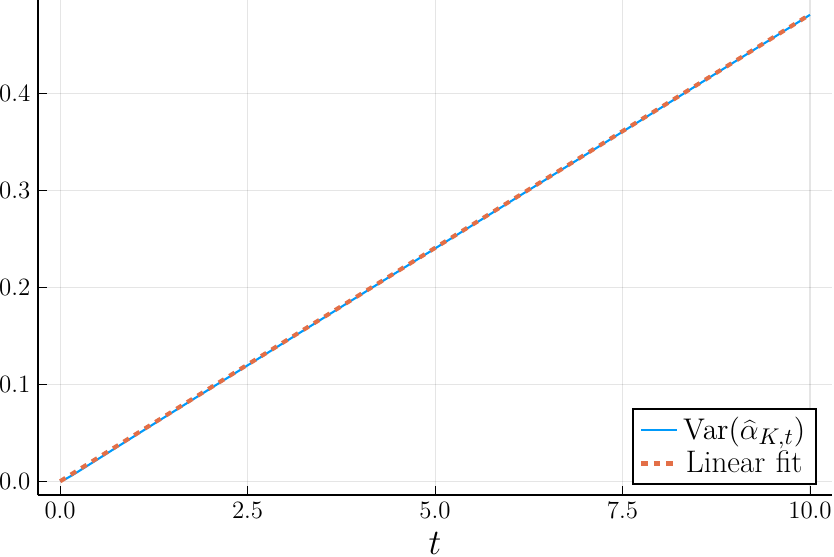}
        \caption{Asymptotic scaling~\eqref{eq:var_green_kubo} of the empirical variance of the Green--Kubo estimator.}
        \label{fig:green_kubo_variance_bias}
    \end{subfigure}

    \caption{Numerical illustration of the Green--Kubo formula and its associated statistical and truncation errors for the two-dimensional overdamped Langevin dynamics with potential~$V$ defined in~\eqref{eq:entropic switch}. The inverse temperature is set to~$\beta = 1$, and the observables are given by~$R = S = \partial_x V$. Simulations are performed up to time~$T = 10$ with time step~$\Delta t = 10^{-3}$, using~$K = 5\times 10^{5}$ independent realizations of the dynamics.}
    \label{fig:green_kubo_bias_full}
\end{figure}

Our aim is to distinguish between the various sources of errors inherent to the numerical estimation. This section is mainly based on~\cite[Section 3.2]{stoltz2024} with, in addition, some precisions on the statistical bias of the numerical estimator. A natural estimator for the right-hand side of~\eqref{eq:green_kubo} is obtained by truncating the time integral to some time~$T < +\infty$ and then replacing the expectation over realizations of the equilibrium dynamics~\eqref{eq:formal SDE} by an average over a finite number~$K$ of such realizations, which leads to
\begin{equation}
  \label{eq:green_kubo_estimator}
  \widehat{\alpha}_{K,T} = \frac{1}{K}\sum_{k=1}^{K} \int_{0}^{T}R{\left(x_t^k\right)}S{\left(x_0^k\right)}\d t.
\end{equation}
In addition, the dynamics needs to be discretized while the time integral needs to be approximated by quadrature formulas. To summarize, the three sources of error of this method are
\begin{itemize}
\item the truncation of the time integral, which leads to an exponentially small bias when the semigroup~$\rme^{t\mathcal{L}_0}$ decays exponentially (see Section~\ref{sec:time_GK});
\item the statistical error in the estimation of the correlation terms, of order~$\sqrt{T/K}$ (see Section~\ref{sec:time_GK});
\item the finiteness of time-step for the discretization of the dynamics and the quadrature formula used for the time integral, which lead to biases of order~$\Delta t^a$ for some integer power~$a \geq 1$ (see Section~\ref{sec:dt_bias_GK}).
\end{itemize}
Usually, the dominant source of error is the statistical error, particularly at long correlation times~\cite{dSOG17}. The computation of~\eqref{eq:green_kubo_estimator} is illustrated in~\cref{fig:green_kubo_bias_full}, where we also investigate the behavior of the truncation and statistical errors. As predicted, the truncation of the time integral leads to an exponentially decaying bias (see Figure~\ref{fig:truncation_bias}) while the variance of the estimator increases linearly in time (see Figure~\ref{fig:green_kubo_variance_bias}). 

\subsubsection{Truncation of time and statistical error}
\label{sec:time_GK}
The truncation bias depends on the asymptotic behavior of correlations. A typical scenario is that the semigroup decays exponentially, as made precise in~\eqref{eq:exp_decay_semigroup}. In this case, 
\begin{equation}
    \label{eq:truncation_bias}
    \left|\E_\pi\left[\widehat{\alpha}_{K,T}\right] - \alpha\right| = \left|\int_{T}^{+\infty} \int_{\cE}\left(\rme^{t\cL_0}R\right)S \,\d\pi \, \d t\right| \leq \frac{C}{\kappa}\|R\|_{L^2(\pi)}\|S\|_{L^2(\pi)} \rme^{-\kappa T}.
\end{equation}
The inequality~\eqref{eq:truncation_bias} suggests to take~$T$ large enough to obtain a small upper bound. However, it has been shown in the recent works~\cite{gastaldello2025dynamical,PSSV24} that the statistical error of~$\widehat{\alpha}_{K,T}$ linearly increases with~$T$. More precisely
\begin{equation}\label{eq:var_green_kubo}
  \lim\limits_{\substack{T \to +\infty}}\frac1T \mathbb{V}\left[\widehat{\alpha}_{K,T}\right] = \frac{2}{K} \|S\|_{L^2(\pi)}^2\left\langle R, -\mathcal{L}_0^{-1}R \right\rangle.
\end{equation}
The proof of this limit is based on the following equality, obtained by using Itô calculus on the solution~$\cR$ of the Poisson equation~$-\cL_0 \cR = R \in L_0^2(\pi)$:
\begin{equation}
  \frac{1}{\sqrt{T}} \int_0^T R(x_t) S(x_0)\, \d t = \frac{S(x_0)}{\sqrt{T}}  \left( \mathcal R(x_0) - \mathcal R(x_T) \right) +
  \frac{S(x_0)}{\sqrt{T}} \int_{0}^{T} \nabla \mathcal R(x_t)^\top \sigma(x_t) \, \d W_t.
\end{equation}
The first term on the left-hand side of the above equality tends to zero in squared expectation in the limit~$T \to +\infty$. Therefore, by using the It\^o isometry for the second term, we obtain
\begin{equation}
\begin{aligned}
  &\lim_{T \to +\infty} \frac{1}{T} \E_{\pi}\left[ \left( \int_0^T R(x_t) S(x_0) \, \d t \right)^2 \right]
  = \lim_{T \to +\infty} \frac{1}{T} \int_{0}^{T} \E_{\pi} \left[ \left| S(x_0) \right|^2 |\sigma^\top \nabla \mathcal R(x_t)|^2 \right]  \, \d t \\
  & \qquad\qquad\qquad\qquad\qquad  = \|S\|^2 \E_{\pi} \left[ \Gamma_{\mathcal R} \right]
  + \lim_{T \to +\infty} \frac{1}{T} \int_{0}^{T} \E_{\pi} \left[ \left| S(x_0) \right|^2 \Pi \Gamma_{\mathcal R}(x_t) \right]  \, \d t,
\end{aligned}
\end{equation}
with~$\Gamma_{\mathcal R} = |\sigma^\top \nabla \mathcal R|^2$. The second term of the above equation is equal to zero. This result can be shown by conditioning over~$x_0$ and using the exponential decay~\eqref{eq:exp_decay_semigroup} of the semigroup. Then, by the fact that~$\E_{\pi} {\left[ \Gamma_{\mathcal R} \right]} = 2\langle R, -\cL_0^{-1} R\rangle$ (see for instance~\cite[Theorem 6.12]{pavliotis2008multiscale}), we have the limit~\eqref{eq:var_green_kubo}.
    
\subsubsection{Time step bias for Green--Kubo formulas}
\label{sec:dt_bias_GK}

We briefly discuss in this section the time-step bias arising from the discretization of the Green--Kubo formula. The result is presented for general diffusions satisfying the following conditions:
\begin{itemize}
\item ergodicity and consistency conditions ensuring a bias of order~$\Delta t^a$ on the invariant probability measure~$\pi_{\Delta t}$ of the numerical scheme (see Section~\ref{sec:error_dt_equilibrium});
\item the transition operator for the numerical scheme can be expanded as~$P_{\Delta t} = \textrm{Id}+\Delta t \cL_0 + \Delta t^2 L_2 +\cdots+\Delta t^{a}L_a + \dots$, with controlled remainder terms;
\item the operator $P_{\Delta t}^{\lceil T/\Delta t \rceil}$ decays exponentially with~$T$ for functions with average zero with respect to~$\pi_{\dt}$, uniformly in~$\Delta t$. This convergence property can be obtained as a consequence of the existence of a Lyapunov function and some minorization conditions.
\end{itemize}
Under such conditions, one can prove the following Riemann-like formula~\cite{leimkuhler2016computation,lelievre2016}:
\begin{equation}\label{eq:time_step_green_kubo}
  \int_{0}^{+\infty}\E_0\left[R(x_t)S(x_0)\right] \, \d t = \Delta t \sum_{n=0}^{+\infty}\E_{\Delta t}\left[\widetilde{R}_{\Delta t}(x^n)S(x^0)\right] + \mathrm{O} (\Delta t^a),
\end{equation}
with
\begin{equation}
  \widetilde{R}_{\Delta t}=\widetilde{\mathcal{R}}_{\Delta t}-\pi_{\Delta t}(\widetilde{\mathcal{R}}_{\Delta t}), \qquad \widetilde{\mathcal{R}}_{\Delta t} = \left(\textrm{Id}+\Delta t L_2 \cL_0^{-1}+\cdots+\Delta t^{a-1}L_a \cL_0^{-1}\right)R.
\end{equation}
Straightforward computations show that the right-hand side of~\eqref{eq:time_step_green_kubo} corresponds to a trapezoidal discretization of the time integral when using schemes of weak order~2. An important consequence of the estimate~\eqref{eq:time_step_green_kubo} is that the asymptotic variance for numerical schemes coincides at dominant order in~$\Delta t$ with the asymptotic variance of the continuous process (as discussed in Section~\ref{sec:error_dt_equilibrium}). These error estimates can also be extended, both from a theoretical and numerical point of view, to Metropolis corrections of discretizations of Langevin-like dynamics, see~\cite{fathi2017improving}. 

\subsection{Control variates}
\label{sec:control_gk}

A possible approach to reduce the variance of Green--Kubo estimators is to use control variates (see Section~\ref{sec:noneq_control_variates}). For simplicity, we describe this approach only in continuous time. To this end, recall that the transport coefficient~\eqref{eq:transport_gk} can be expressed in terms of the solution to an appropriate Poisson equation. More precisely, it holds that
\begin{equation}
\alpha = \int_{0}^{+\infty}\E_\pi\left[R(x_t)S(x_0)\right] \, \d t = \left\langle \mathcal R, S \right\rangle_{L^2(\pi)} = \left\langle R, \mathcal S^\star \right\rangle_{L^2(\pi)},
\end{equation}
where $\mathcal L_0^\star$ denotes the $L^2(\pi)$ adjoint of the generator~$\mathcal L_0$, while $\mathcal R$ and $\mathcal S^\star$ are the unique $L^2_0(\pi)$ solutions to the Poisson equations
\begin{equation}\label{eq:poisson_eqs_cv}
  - \mathcal L_0 \mathcal R = R, 
  \qquad 
  - \mathcal L_0^\star \mathcal S^\star = S. 
\end{equation}
In most settings, the generator $\cL_0$ can be explicitly decomposed into its symmetric and antisymmetric parts in~$L^2(\pi)$. The adjoint operator~$\mathcal L_0^\star$ is then obtained simply by flipping the sign of the antisymmetric part in the generator, while keeping the symmetric part unchanged. Solving Poisson equations associated with~$\cL_0^\star$ therefore presents the same level of difficulty as solving Poisson equations associated with~$\cL_0$.

The control-variate approaches deployed in~\cite{PSSV24} are based on two prerequisites:
\begin{itemize}
\item 
  The first prerequisite is that it is possible to calculate approximate solutions to at least one of the Poisson equations in~\eqref{eq:poisson_eqs_cv}.
  We denote these approximate solutions by~$\mathfrak R$ and~$\mathfrak S^\star$,
  respectively.
  This assumption is somewhat restrictive,
  as solving partial differential equations accurately is usually possible only in low dimensions.
  However, the approximate solutions $\mathfrak R$ and~$\mathfrak S^\star$ do not need to be extremely accurate to provide significant variance reduction,
  and so using control variates may prove useful even in moderately high dimensions. One can for instance resort to spectral methods~\cite{RS19,abdulle2017,pavliotis2023} or rely on machine learning approaches~\cite{MR3874585,MR3881695,MR4807807}.

\item 
        The second prerequisite is that,
        given two functions~$f,g \in L^2(\pi)$,
        it is possible to calculate precisely the quantity
        \(
            \left\langle f, g \right\rangle_{L^2(\pi)}
        \)
        at a computational cost that is negligible, compared to that of calculating Green--Kubo estimators of the type~\eqref{eq:green_kubo_estimator}.
        In low to moderate dimensions, 
        it is usually inexpensive to generate independent and identically distributed (i.i.d.) samples~$\{X^1, X^2, \dotsc\}$ from $\pi$,
        in which case a precise approximation is provided by Monte Carlo sampling:
        \begin{equation}
          \left\langle f, g \right\rangle_{L^2(\pi)}
          \approx
          \frac{1}{J} \sum_{j=1}^{J} f(X^j) g(X^j).
        \end{equation}
\end{itemize}

Having detailed these restrictions,
we now motivate the construction of an improved estimator by using the approximate solutions~$\mathfrak R$ and $\mathfrak S^\star$ in control variates.
The key idea of control-variate approaches for Green--Kubo formulas is to rewrite the transport coefficient as
\begin{equation}
\begin{aligned}
    \left\langle \mathcal R, S \right\rangle_{L^2(\pi)}
    &= 
    \left\langle \mathfrak R, S \right\rangle_{L^2(\pi)}
    + \left\langle \mathcal R - \mathfrak R, S \right\rangle_{L^2(\pi)} \\
    &= 
    \left\langle \mathfrak R, S \right\rangle_{L^2(\pi)}
    + \left\langle R + \cL_0 \mathfrak R, \mathcal S^\star \right\rangle_{L^2(\pi)} \\
    &= 
    \left\langle \mathfrak R, S \right\rangle_{L^2(\pi)}
    + \left\langle R + \cL_0 \mathfrak R, \mathfrak S^\star \right\rangle_{L^2(\pi)}
    + \left\langle - \cL_0^{-1} (R + \cL_0\mathfrak R), S + \cL_0^\star \mathfrak S^\star \right\rangle_{L^2(\pi)}.
\end{aligned}
\end{equation}
In view of the second assumption above, 
the first and second terms on the right-hand side can be calculated precisely.
The remaining term may then be recast as an integrated correlation,
which can be approximated by using the estimator~\eqref{eq:green_kubo_estimator},
only with $R$ and $S$ substituted by the residuals $R + \cL_0\mathfrak R$ and $S + \cL_0\mathfrak S^\star$,
respectively. Specifically, the improved estimator reads
\begin{equation}\label{eq:green_kubo_estimator_cv}
  \widehat{\alpha}_{K,T}^{\rm cv} = 
  \left\langle \mathfrak R, S \right\rangle_{L^2(\pi)}
  + \left\langle R + \cL_0 \mathfrak R, \mathfrak S^\star \right\rangle_{L^2(\pi)}
  + \frac{1}{K}\sum_{k=1}^{N} \int_{0}^{T}(R + \cL_0\mathfrak R){\left(x_t^k\right)} (S + \cL_0^\star \mathfrak S^\star){\left(x_0^k\right)} \, \d t.
\end{equation}
It is now clear from~\eqref{eq:var_green_kubo} that, 
if the residuals are small,
then variance reduction should follow.
A few remarks are given below in order to conclude this section:
\begin{itemize}
    \item
        If either $\mathfrak R = \mathcal R$ or~$\mathfrak S^\star = \mathcal S^\star$,
        and if one furthermore omits any potential statistical error involved in the approximation of the first two terms in~\eqref{eq:green_kubo_estimator_cv},
        then the estimator~\eqref{eq:green_kubo_estimator_cv} has zero variance.

    \item
        The approach we described is useful even when only one of the two Poisson equations in~\eqref{eq:poisson_eqs_cv} is solved approximately.
        In this case we can simply set~$\mathfrak R = 0$ or~$\mathfrak S^\star = 0$ for the other approximate solution,
        and use the same estimator~\eqref{eq:green_kubo_estimator_cv}.

    \item 
        For this approach to yield computational benefits,
        it is crucial that the approximate solutions~$\mathfrak R$ and~$\mathfrak S^\star$ be inexpensive to evaluate.
        This is more important for the solution~$\mathfrak R$,
        which is evaluated at all integration times,
        while the approximate solution~$\mathfrak S^\star$ is evaluated only at the initial time.
        In the numerical experiments conducted in~\cite{PSSV24}, 
        it was observed that letting $\mathfrak R = 0$ and  computing only an approximation~$\mathfrak S^\star$ already yields good results in terms of variance reduction,
        and this gives an estimator that is much less expensive to evaluate than using both~$\mathfrak R$ and~$\mathfrak S^\star$.
\end{itemize}

\subsection{Importance sampling in path space}
\label{sec:Girsanov}

The main result of Section~\ref{sec:time_GK} is that the statistical error of the estimator~\eqref{eq:green_kubo_estimator} for the Green--Kubo formula~\eqref{eq:green_kubo} is proportional to~$\langle R, -\mathcal L_0^{-1} R\rangle$. One idea to reduce this term (which encapsulates the dynamical behavior) is to consider a modified dynamics with generator~$\ocL_\alpha$ such that~$\langle R, -\ocL_\alpha^{-1} R\rangle$ is smaller than~$\langle R, -\cL_0^{-1} R\rangle$. Such a modification was studied in~\cite{gastaldello2025dynamical}, where important sampling using Girsanov theorem is applied to~\eqref{eq:green_kubo} -- an idea inspired by various works in theoretical chemistry~\cite{chodera2011dynamical,donati2017girsanov,donati2018girsanov}. More precisely, consider for~$\alpha \in \R$ and~$u \colon \cX \to \R^d$ the following perturbation of the dynamics~\eqref{eq:formal SDE}:
\begin{equation}
  \d x_t^\alpha = b(x_t^\alpha)\, \d t + \alpha\sigma(x_t^\alpha)u(x_t^\alpha)\, \d t + \sigma(x_t^\alpha)\,\d W_t. 
\end{equation}
In this case, $\ocL_\alpha = \cL_0 + \alpha u^\top \sigma^\top \nabla$. The aim of the additional drift is to reduce the metastability of the system, typically by favoring transitions from one metastable state to another. This can be achieved for instance by free energy biasing, namely identifying a so-called collective variable which captures the slow degrees of freedom of the system and considering for~$u$ the gradient of the free energy associated with the collective variable~\cite{lelievre_2010}. The parameter~$\alpha$ allows to tune the magnitude of the forcing. A traditional compromise is that larger values of~$\alpha$ lead to a better exploration of the configuration space and smaller correlation times, hence smaller values of~$\langle R, -\ocL_\alpha^{-1} R\rangle$; but the associated Girsanov weights are larger. 

To make this discussion more quantitative, we introduce the counterpart of the estimator~$\widehat{\alpha}_{K,T}$ defined in~\eqref{eq:green_kubo_estimator} for the unperturbed dynamics:
\begin{equation}
  \widehat{\alpha}_{K,T}^\alpha = \frac{1}{K}\sum_{k=1}^{N} \exp \left(- \alpha\int_0^T u(x_t^{\alpha , k}) \cdot \d W_t - \frac{\alpha^2}{2} \int_0^T |u(x_t^{\alpha , k})|^2 \, \d t \right) \int_{0}^{T}R\left(x_t^{\alpha , k}\right)S\left(x_0^{\alpha , k}\right) \, \d t.
\end{equation}
The Girsanov theorem ensures that this estimator is consistent: $\E_\pi\left[\widehat{\alpha}_{K,T}\right] = \E_\pi\left[\widehat{\alpha}_{K,T}^\alpha\right]$. Its variance however differs from the one of the estimator~\eqref{eq:green_kubo_estimator}. More precisely, 
\begin{equation}
  \mathbb{V}\left[\widehat{\alpha}_{K,T}^\alpha\right] = \E_\pi\left[\left(\widehat{\alpha}_{K,T}^\alpha\right)^2\right] - \E_\pi\left[\widehat{\alpha}_{K,T}\right]^2 =: F_T(\alpha) - \E_\pi\left[\widehat{\alpha}_{K,T}\right]^2.
\end{equation}
Under suitable assumptions, it is proved in~\cite{gastaldello2025dynamical} that, for any~$T>0$, there exists a unique~$\alpha \in \R$ that minimizes~$F_T(\alpha)$ and therefore the variance. Furthermore, it is shown that the relative reduction of the raw second moment is proportional to~$1/T$ when~$T$ goes to infinity. This result, coupled with the numerical simulations in~\cite{gastaldello2025dynamical}, demonstrates that the benefits in terms of variance reduction provided by this method are in fact rather modest for large~$T$.

\subsection{Alternative fluctuation formulas}
\label{sec:alternative_fluctuation}

We discuss in this section two alternative correlation formulas which can be used to estimate transport coefficients.

\paragraph{Einstein formulas.}
The Einstein relation~\eqref{eq:einstein_relation} establishes a link between the mobility of a particle and its effective diffusion coefficient entering in Fick's law. Therefore, the mobility can be calculated by instead estimating the effective diffusion coefficient,
which is defined as (recall the definition~\eqref{eq:unperiodized_displacement})
\begin{equation}\label{eq:ein_diff}
  \mathfrak D := \lim_{T \to +\infty} \expect \left[ \frac{|Q_T - Q_0|^2}{2dT} \right].
\end{equation}
For an ergodic dynamics, this equality holds for any initial configuration. For simplicity,
we consider from now on the Langevin dynamics~\eqref{eq:langevin} in dimension $d = 1$,
and proceed to rewrite formula~\eqref{eq:ein_diff} in a manner that can be generalized to other transport coefficients.
For Langevin dynamics, the diffusion coefficient~$\mathfrak D$ can be rewritten as a double time integral:
\begin{equation}
\begin{aligned}
  \mathfrak D
  &= \lim_{T \to +\infty} \frac{1}{2T} \expect \left[ \left(\int_{0}^{T} p_t \, \d t\right)^2 \right]  \\
  &= \lim_{T \to +\infty} \frac{1}{2T} \int_{0}^{T} \int_{0}^{T} \expect [p_s p_t] \, \d s \, \d t
  = \lim_{T \to +\infty} \frac{1}{T} \int_{0}^{T} \int_{0}^{t} \expect [p_s p_t] \, \d s \, \d t.
\end{aligned}
\end{equation}
While this derivation is concerned with the particular case of the mobility, it is possible to prove that, for any function $w \colon [0, 1] \to \mathbb R$ that is continuous at zero and satisfies~$w(0) = 1$  and for any observables $R, S \in L^2_0(\pi)$, the transport coefficient~\eqref{eq:transport_gk} can be obtained as
\begin{equation}\label{eq:einstein}
  \alpha = \lim_{T \to +\infty} \frac{1}{T} \int_{0}^{T} \int_{0}^{t} w\left(\frac{t-s}{T}\right) \expect \left[ R(x_t) S(x_s) \right] \, \d s \, \d t.
\end{equation}
We refer to~\cite[Section~2.3]{PSSV24} for a proof of this claim.
By a slight abuse of terminology,
we refer to formulas such as~\eqref{eq:einstein} as \emph{Einstein formulas}.
The reformulation~\eqref{eq:einstein} of the transport coefficient suggests the following estimator, 
as an alternative to the Green--Kubo estimator~\eqref{eq:green_kubo_estimator}:
\begin{equation}\label{eq:einstein_est}
  \widehat \alpha_{K,T}^{\rm Ein}
   = \frac{1}{TK}\sum_{k=1}^K \int_0^T\int_0^t w\left( \frac{t-s}{T} \right) R(x_t^k) S(x_s^k) \, \d s \, \d t.
\end{equation}
Since this estimator includes a time average
in addition to the ensemble average,
we expect it to have lower variance for large~$T$.
In fact, 
it is possible to prove that,
under mild technical conditions detailed in~\cite[Proposition 2.6]{PSSV24},
the variance of the estimator~\eqref{eq:einstein_est} converges to a finite limit as~$T \to +\infty$:
\begin{equation}
  \lim_{T \to +\infty} \mathbb V \left[\widehat \alpha_{K,T}^{\rm Ein}\right]  
  = \frac{4}{K}
  \left(\int_{0}^{1} (1- u) \, w(u)^2 \, \d u\right)
  \left\langle R, - \mathcal L_0^{-1} R \right\rangle \left\langle S, - \mathcal L_0^{-1} S \right\rangle.
\end{equation}
As described in Section~\ref{sec:control_gk}, the variance of the Green--Kubo estimator can be reduced using control variates based on approximate solutions to the Poisson equations~\eqref{eq:poisson_eqs_cv}. The variance of the estimator~\eqref{eq:einstein_est} can be reduced in a similar way. The idea is exactly the same: 
to calculate the first two terms in~\eqref{eq:green_kubo_estimator_cv} using a traditional sampling method,
as these terms do not involve expectations with respect to path measures,
and then to estimate the third term using the estimator~\eqref{eq:einstein_est} with~$R$ and~$S$ substituted by the corresponding residuals,
assuming that approximate solutions to the Poisson equations~\eqref{eq:poisson_eqs_cv} are available.

To conclude this section,
we briefly discuss the choice of the weight function~$w$.
Various weight functions have been proposed in the literature,
see~\cref{table:weights}.
In general, there is a trade-off between minimizing the bias and minimizing the variance of the estimator~\eqref{eq:einstein_est}.
Since the statistical error usually dominates the systematic error in applications,
most weight functions proposed in the literature assign more weight to short correlations (i.e., small time lags) than to long correlations (i.e., large time lags).
Finally, let us mention that it is possible to consider a more general weight of the form~$\widetilde w(t - s, T)$ instead of $w\left(\frac{t-s}{T}\right)$ in~\eqref{eq:einstein} and~\eqref{eq:einstein_est},
see~\cite{lu2019}.
This makes it possible to devise estimators such that the variance is not only bounded, 
but decreases to zero in the limit as~$T \to +\infty$.

\begin{table}[htb!]
    \footnotesize
    \centering
    \caption{
        Expression for various weight functions $w(t)$ for $0\leq t\leq 1$,
        otherwise defined as zero; see \cite{wang2012} and references therein.
    }
    \label{table:weights_num_illustration}
    \renewcommand{\arraystretch}{2.0}
    \begin{tabular}{ll}
        \toprule
        \textbf{Weight function} & \textbf{Expression for $w(t)$}\\\midrule
        Constant & $w(t) = 1$ \\
        Bartlett & $w(t) = 1 - t$ \\
        Parzen & $w(t) = \begin{cases} 1 - 6t^2 + 6t^3, &t\leq 0.5 \\ 2(1-t)^3, &t > 0.5\end{cases}$ \\
        Tukey--Hanning & $w(t) = \dfrac{1+\cos(\pi t)}{2}$ \\
        Parzen--Riesz & $w(t) = 1 - t^2$ \\
        Parzen--Geometric & $w(t) = \dfrac{1}{1+t}$ \\
        Parzen--Cauchy & $w(t) = \dfrac{1}{1+t^2}$ \\
        \bottomrule
    \end{tabular}
    \label{table:weights}
\end{table}

\paragraph{Martingale product.} 
Recall the standard dynamics~\eqref{eq:formal SDE} and its perturbation~\eqref{eq:nemd_dynamics}.
In~\cref{subsec:nemd_dynamics}, it has been discussed how to compute transport coefficients using linear response through the expression~\eqref{eq:linear_response}. A different approach to estimating this expression is based on the so-called
\textit{likelihood ratio} formula, which is obtained by using the Girsanov theorem and formally passing first to the small
forcing limit, and then to the infinite-time limit. As shown in~\cite{glynn2019likelihood,plechac_nemd,plechac2022,wang2019steady}, for~$R\in L_0^2(\pi)$,
\begin{equation}\label{eq:CLR}
  \alpha=\underset{\eta\to 0}{\lim}\,\frac{\E_{\pi_\eta}\left[R\right]}{\eta} = \underset{t\to +\infty}{\lim}\,\E_{\pi}\left[\left(\frac1t\int_{0}^{t}R(x_s)\,\d s\right)Z_t\right], \qquad Z_t = \int_{0}^{t}u(x_s)^\top\,\d W_s,
\end{equation}
with~$\sigma u = F$ and~$F$ being the non-gradient forcing~\eqref{eq:non_conservative_force}. The product with the martingale term~$Z_t$ is what suggests the name of this method.
The two main advantages of this method are:
\begin{itemize}
\item As achieved by the Green--Kubo formula~\eqref{eq:green_kubo}, the coefficient~\eqref{eq:linear_response} can be computed from the dynamics at equilibrium~\eqref{eq:formal SDE} by reweighting the trajectory average of the observable with the martingale~$Z_t$.
    \item Since~$R\in L_0^2(\pi)$ in~\eqref{eq:CLR}, it is in fact the \textit{centered likelihood ratio method} that is used. Furthermore, if~$R$ is taken smooth, it is shown in~\cite{plechac_nemd} that there exists a constant $C>0$ such that
    \begin{equation}
    \forall t > 0, \qquad
    \mathbb{V}\left[\left(\frac{1}{t}\int_0^t R(x_s)\, \d s\right) Z_t\right] \leq C.
    \end{equation}
    This result shows that the variance of the estimator~$\left(\frac{1}{t}\int_0^t R(x_s)\, \d s\right) Z_t$ is bounded in terms of the integration time, which is a very interesting property for long time simulations.
\end{itemize}
Furthermore, the papers~\cite{plechac_nemd,plechac2022} propose a time discretization of this method for the overdamped Langevin dynamics~\eqref{eq:overdamped_langevin} and Langevin dynamics~\eqref{eq:langevin},
and they study the associated time-step bias.
More precisely, it is shown that in both cases, if the numerical scheme is of weak order~$a\in\{1,2\}$, then the time-step bias is also of order~$a$. 


\section{Transient methods}
\label{sec:transient}

In this section, we turn to transient methods, where transport coefficients are obtained by monitoring the relaxation of the system to its steady state, either starting under the equilibrium Boltzmann--Gibbs measure and converging to the nonequilibrium steady state (see Section~\ref{sec:TTCF}), or starting under a perturbed Boltzmann--Gibbs measure and relaxing to the equilibrium state (see Section~\ref{sec:transient_subtraction}).

\subsection{Relaxation to the nonequilibrium steady state}
\label{sec:TTCF}

We briefly review here the \emph{transient time correlation function} (TTCF) method for calculating transport coefficients and nonlinear response functions~\cite{morriss1987,evans1988}. This section is partly based on~\cite[Sections~2.4 and~3.3]{MSS24}. To formally derive the method, we consider the generic dynamics~\eqref{eq:formal SDE} with state space~$\mathcal{X}$, driven out of equilibrium by an additional forcing term~$F\colon \mathcal \cX \to \R^d$:
\begin{equation}\label{eq:noneq_dynamics_started_eq}
  \d y_t^\eta = \left( b(y_t^\eta) + \eta F(y_t^\eta) \right) \, \d t + \sigma(y_t^\eta) \, \d W_t, \qquad y_0^\eta \sim \mu.
\end{equation}
As in previous sections, we assume that~\eqref{eq:noneq_dynamics_started_eq} admits a unique invariant probability measure~$\pi_\eta$ (with~$\pi_0 = \pi$), and denote by~$\mathcal L_\eta = \mathcal L + \eta \widetilde {\mathcal L}$ the generator of the above dynamics, with~$\widetilde {\mathcal L} = F^\top \nabla$. By It\^o's formula and the definition of the operator semigroup~$t \mapsto \rme^{t \mathcal L_{\eta}}$, the average value of a sufficiently smooth response function~$R \colon \mathcal X \to \R$ satisfies
\begin{equation}
  \expect \left[ R(y^{\eta}_t) \right] - \expect \left[ R(y^{\eta}_0) \right]
  = \int_{0}^{t} \expect \left[ \mathcal L_{\eta} R (y^{\eta}_t) \right]\, \d t
  = \int_{0}^{t} \expect \left[ \rme^{t \mathcal L_{\eta}} (\mathcal L_{\eta} R) (y^{\eta}_0) \right]\, \d t.
\end{equation}
Rewriting the expectation as an inner product on~$L^2(\pi)$ and recalling the $\star$ notation for adjoints on this space (as used in~\eqref{eq:generator_with_adjoints} for instance), we deduce that
\begin{equation}
  \expect \left[ R(y^{\eta}_t) \right] - \expect \left[ R(y^{\eta}_0) \right]
  = \int_{0}^{t} \left\langle \1, \mathcal L_{\eta} \rme^{t \mathcal L_{\eta}}  R \right\rangle_{L^2(\pi)} \, \d t
  = \int_{0}^{t} \left\langle \mathcal L_{\eta}^\star \1, \rme^{t \mathcal L_{\eta}}  R \right\rangle_{L^2(\pi)} \, \d t.
\end{equation}
Since $\mathcal L^\star \1 = 0$ by invariance of~$\pi$ under the reference dynamics
and using that the conjugate response~\eqref{eq:conjugate_flux} is equal to $S = \widetilde {\mathcal L}^\star \1$,
we obtain
\begin{equation}
  \frac{\expect \left[ R(y^{\eta}_t) \right] - \expect \left[ R(y^{\eta}_0) \right]}{\eta}
  = \int_{0}^{t} \left\langle S, \rme^{t \mathcal L_{\eta}}  R \right\rangle_{L^2(\pi)} \, \d t 
  = \int_{0}^{t} \expect \left[ S(y^{\eta}_0) R(y^{\eta}_t) \right]\, \d t.
\end{equation}
By passing to the limit~$t \to +\infty$, and denoting by~$\E_\eta$ the expectation with respect to the steady-state probability measure of the nonequilibrium dynamics~$(y_t^\eta)_{t \geq 0}$, we conclude that
\begin{equation}\label{eq:nl_response}
  \frac{\E_\eta[R]-\E_0[R]}{\eta} = \int_0^{+\infty} \E\left[ S(y_0^\eta) R(y_t^\eta) \right] \, \d t,
\end{equation}
the expectation on the right hand side being taken with respect to equilibrium initial conditions~$y_0^\eta \sim \mu$ evolved through the nonequilibrium dynamics~\eqref{eq:noneq_dynamics_started_eq}.
The usual Green--Kubo formula is recovered in the limit of $\eta \to 0$ in~\eqref{eq:nl_response}. In view of~\eqref{eq:nl_response}, the natural TTCF estimator
\begin{equation}\label{eq:TTCF_estimator}
  \widehat A^{\rm TTCF}_{K,T} = \frac{1}{K} \sum_{k=1}^K \int_0^T R\left(y_t^{\eta,k}\right) S\left(y_0^{\eta,k}\right) \, \d t, \qquad y_0^{\eta,k} \stackrel{\rm i.i.d.}{\sim} \pi,
\end{equation}
is asymptotically unbiased as~$T \to +\infty$, for any value of~$\eta$. Therefore, this estimator makes it possible to capture the full nonlinear response, in contrast to the transient methods based on relaxation to the equilibrium steady state described in Section~\ref{sec:transient_subtraction}.

\subsection{Relaxation to equilibrium}
\label{sec:transient_subtraction}

We present here transient methods based on a relaxation to the equilibrium state, starting from a perturbation thereof. This account is based on~\cite{MSS24} (in particular Sections~2.3 and~3.1 to~3.3), building upon the \emph{subtraction technique} proposed in~\cite{ciccotti1975} and further explored in \cite{ciccotti1979}. Applications to the computation of the mobility and shear viscosity of Lennard-Jones fluids can be found in~\cite[Section~4.3]{MSS24}.

\paragraph{Presentation of the method.}
The computation of transport coefficients from a transient relaxation to equilibrium is performed by integrating the reference dynamics~\eqref{eq:formal SDE} for initial conditions distributed according to a perturbation~$\widetilde{\pi}_\eta$ of the equilibrium probability measure~$\pi$ (whose realizations are noted by~$(x_t^\eta)_{t \geq 0}$), and monitoring the return to stationarity via a time integral. The probability measure~$\widetilde{\pi}_\eta$ is of the form
\begin{equation}\label{eq:init_dist}
  \widetilde{\pi}_\eta = (1+\eta S)\pi + \mathrm{O}(\eta^2).
\end{equation}
Although not immediately clear, the time integral of the expectation of $R(x_t^\eta)$, when divided by $\eta$, converges to the transport coefficient $\alpha$ as $\eta$ goes to zero:
\begin{equation}\label{eq:transient_def}
  \alpha = \lim_{\eta\to 0} \frac{1}{\eta}\int_0^{+\infty} \E\bigl[R(x_t^\eta)\bigr] \, \rmd t.
\end{equation}
The equality~\eqref{eq:transient_def} can be formally derived as follows, since~$R$ (and hence~$\rme^{t \calL_0}R$) has average zero with respect to~$\pi$: 
\begin{equation}\label{eq:gk_equiv2}
\begin{aligned}
    \frac{1}{\eta}\int_0^{+\infty} \E\bigl[R(x_t^\eta)\bigr] \, \rmd t
    &= \frac{1}{\eta} \int_0^{+\infty} \int_\mathcal{X} \left(\rme^{t\calL_0} R\right) \, \rmd \widetilde{\pi}_\eta \, \rmd t \notag \\
    &= \frac{1}{\eta}\int_0^{+\infty} \int_\mathcal{X} \rme^{t\calL_0} R \, \rmd\pi \, \rmd t + \int_0^{+\infty} \int_\mathcal{X} \left(\rme^{t\calL_0} R\right) S \, \rmd\pi \, \rmd t + \mathrm{O}(\eta) \notag \\
    &= \int_0^{+\infty} \int_\mathcal{X} \left(\rme^{t\calL_0} R\right) S \, \rmd \pi \, \rmd t + \mathrm{O}(\eta) \notag \\
    &= \int_0^{+\infty} \E_\pi{\left[R(x_t^0)S(x_0^0)\right]} \, \rmd t + \mathrm{O}(\eta).
\end{aligned}
\end{equation}
The first term on the right-hand side of the last term is the transport coefficient~$\alpha$ expressed through the Green--Kubo formula~\eqref{eq:green_kubo}.

In practice, the perturbed initial distribution is constructed by taking the image measure of the equilibrium distribution~$\pi$ through some transformation~$\Phi_\eta$:
\begin{equation}\label{eq:deterministic_map}
  \widetilde{\pi}_\eta = \Phi_\eta \# \pi = (1+\eta S)\pi + \mathrm{O}(\eta^2).
\end{equation}
This means that, for any bounded measurable test function~$\varphi\colon \cX \to \R$,
\begin{equation}
  \int_\mathcal{X} \varphi \, \rmd \widetilde{\pi}_\eta = \int_\mathcal{X} \varphi \circ \Phi_\eta \, \rmd \pi.
\end{equation}
It is natural to look for a map $\Phi_\eta$ of the form
\begin{equation}
  \Phi_\eta(x) = x + \eta\varphi_1(x),
  \label{eq:Phi_map}
\end{equation}
where $\varphi_1$ is determined by the condition~\eqref{eq:deterministic_map}. It is in fact given by a solution to a partial differential equation, which can be analytically solved for instance for underdamped Langevin dynamics and forcings on the drift (see~\cite{MSS24}). In the latter simple case, $\varphi_1(q,p) = (0,F(q))$.

Numerically estimating~\eqref{eq:transient_def} requires first approximating the limit with (sufficiently small) finite $\eta$, truncating the time integral to finite $T$, and approximating the expectation with an average over $K$ realizations of the dynamics started from i.i.d.\ initial conditions $x_0^\eta = \Phi_\eta(x_0^0)$ with~$x_0^0 \sim \pi$. This leads to the following estimator for~\eqref{eq:transient_def}:
\begin{equation}
  \widehat{\alpha}_{T,K,\eta} = \frac{1}{\eta K}\sum_{k=1}^K \int_0^T R(x_t^{\eta,k}) \, \rmd t.	
  \label{eq:T_estimator}
\end{equation}
Although these approximations lead to several sources of bias in~\eqref{eq:T_estimator} (made precise in~\cite[Section~3.2.2]{MSS24}), the primary concern associated with~\eqref{eq:T_estimator} is its very large variance, of order~$T/\eta^2$. Indeed, since~\eqref{eq:T_estimator} involves a time integral of a function with average zero with respect to~$\pi$, computations similar to the ones considered in Section~\ref{eq:sampling_error_continuous} give
\begin{equation}
  \lim_{T\to+\infty} \frac{1}{T} \bbV\left(\widehat{\alpha}_{T,K,\eta}\right) = \frac{2}{K\eta^2}\int_\mathcal{X} R\left(-\calL_0^{-1} R\right) \, \rmd \pi.
  \label{eq:asym_var_transient}
\end{equation}
Unlike the usual NEMD or Green--Kubo estimators of transport coefficients discussed in Sections~\ref{sec:NEMD_error_analysis} and~\ref{subsec:error_green_kubo}, respectively, the variance of~\eqref{eq:T_estimator} is magnified by two distinct contributions: a division by a factor of~$\eta$, as in NEMD; and a time integral, as in Green--Kubo formulas. This leads to a variance of order~$T/\eta^{-2}$, much higher than its NEMD and Green--Kubo counterparts. 

\paragraph{Variance reduction with a coupling method.}
The very large variance~\eqref{eq:asym_var_transient} motivates modifying the estimator with the use of variance reduction techniques, in particular to get rid of the $\eta^{-2}$ contribution and obtain bounds uniform in~$\eta$ as for TTCF. To this end, \emph{couplings} are considered as a control variate. More precisely, we consider the coupling $(x_t^\eta,y_t^0)$, where the processes $x_t^\eta$ and $y_t^0$ are evolved according to the same underlying reference dynamics but have different initial conditions:
\begin{equation}
\begin{cases}
\begin{aligned}
  dy_t^0 &= b(y_t^0) \, \rmd t + \sigma(y_t^0) \, \rmd W_t, \qquad y_0^0 \sim \pi, \\
  dx_t^\eta &= b(x_t^\eta) \, \rmd t + \sigma(x_t^\eta) \, \rmd \widetilde{W}_t, \qquad x_0^\eta = \Phi_\eta(y_0^0),
\end{aligned}
\end{cases}
\label{eq:coupled_dynamics_2}
\end{equation}
where $W_t$ and $\widetilde{W}_t$ are standard $d$-dimensional Brownian motions, respectively. The transport coefficient $\alpha$ can then be computed as
\begin{equation}
  \alpha = \lim_{\eta\to 0} \frac{1}{\eta}\int_0^{+\infty} \E\left( R(x_t^\eta) - R(y_t^0) \right) \, \rmd t.
  \label{eq:tc_subtraction}
\end{equation}
Note that $\int_0^{+\infty} R(y_t^0) \, \rmd t$ acts as a control variate since $\E(R(y_t^0))=0$ for all $t\geq 0$. The expression~\eqref{eq:tc_subtraction} admits the following natural estimator, carried out with independent initial conditions for the couple $(x_t^{\eta,k},y_t^{0,k})_{t\geq 0}$ for $1\leq k\leq K$ and independent realizations of the dynamics~\eqref{eq:coupled_dynamics_2}:
\begin{equation}
  \widehat{\alpha}_{T,K,\eta}^{\rm sub} = \frac{1}{\eta K} \sum_{k=1}^K \int_0^T \left[ R(x_t^{\eta,k}) - R(y_t^{0,k}) \right] \, \rmd t.
  \label{eq:ts_estimator}
\end{equation}
A sufficient condition for~\eqref{eq:ts_estimator} to have smaller variance than the standard estimator~\eqref{eq:T_estimator} is
\begin{itemize}
\item for the trajectories to start with initial conditions that are at distance of order~$\eta$ one from another. This is ensured by construction from~\eqref{eq:Phi_map} and the fact that~$x_0^\eta = \Phi_\eta(y_0^0)$;
\item trajectories should stay close for times of order $1/\kappa$, where $\kappa$ is the relaxation rate of the system to the stationary state (see~\eqref{eq:exp_decay_semigroup}). In order to achieve this, one needs to couple the Brownian motions in~\eqref{eq:coupled_dynamics_2}, for instance through a synchronous coupling (see Section~\ref{subsec:NEMD_coupling}). 
\end{itemize}
When the trajectories of the reference dynamics do not decorrelate too fast (at most exponentially in time), one can show that the variance of the estimator~$\widehat{\alpha}_{T,K,\eta}^{\rm sub}$ in~\eqref{eq:ts_estimator} is uniformly bounded in~$\eta$, with however a scaling with respect to time which depends on the properties of the reference dynamics, and in particular its possible dissipativity; see~\cite[Section~3.2.3]{MSS24} for precise results.

\section{Extensions and perspectives}
\label{sec:extensions}

We conclude this review with some practical considerations, as well as perspectives for future work. In a nutshell, a lot remains to be explored from a mathematical and algorithmic perspective in order to achieve actual variance reduction in the estimation of transport coefficients.

\paragraph{Which class of methods to choose?}
Most practitioners rely either on NEMD techniques or Green--Kubo formulas. Transient techniques are less commonly employed. The reason why one method is used rather than another mostly depends on personal habits. We can only advocate cross comparing the estimators from different techniques, in order to improve the reliability of the results. 

Let us nonetheless discuss some possible (dis)advantages of NEMD techniques and methods based on Green--Kubo formulas. The main interest of Green--Kubo formulas is that multiple integrated correlation functions can be estimated for realizations of the equilibrium dynamics, which allows to estimate several transport coefficients at a time with a single simulation. In contrast, a different nonequilibrium perturbation has to be considered for every transport coefficient in the NEMD approach. Moreover, in order to ensure that the forcing magnitude is appropriate, one needs in principle to check the linearity of the response in nonequilibrium simulations, typically by computing the response for several values of the forcing magnitude, which further increases the computational cost (although these computations can be parallelized in a straightforward way).

However, if one wants to make a full use of a single long MD trajectory, the postprocessing of Green--Kubo simulations is less straightforward than for nonequilibrium simulations, for which plain time averages are considered. Moreover, it is possible to increase the magnitude of the forcings in nonequilibrium simulations in order to enhance the response to measure; whereas the correlation functions appearing in the Green--Kubo cannot be enhanced in such a way and may therefore require more computational time to emerge out of the statistical noise.

\paragraph{Extensions to other dynamics.}
Our paradigmatic examples of dynamics in this review are underdamped and overdamped Langevin dynamics. In practice, one can be led to considering dynamics with a noise which can be even more degenerate. One example of such systems are atom chains; see Section~\ref{ss:thermtransp}. Other possible dynamics are listed at the end of Section~\ref{sec:sampling_methods}. 

The variance reduction techniques we mentioned in this review can for some of them be adapted to these settings, but not all of them. For DPD for instance, which is important for a wide range of complex fluids and soft matter applications, one example of method that can be adapted is the stochastic Norton dynamics from Section~\ref{sec:NEMD_Norton}. It has been numerically demonstrated in~\cite{wu2025} that the responses profiles for both the NEMD and stochastic Norton dynamics approaches coincide in both linear and nonlinear regimes, which indicates that, similar to the situation for underdamped Langevin dynamics, the stochastic Norton dynamics can indeed act as an alternative approach for the computation of transport coefficients, including the mobility and the shear viscosity, as the NEMD dynamics in the DPD context. 

\paragraph{Acknowledgements.}
This work benefited from a CECAM/CCP5 sandpit funding for the academic year 2025, awarded to X.S. and G.S.. N.B., L.C., S.D., R.G., G.S. and U.V. were funded by the European Research Council (ERC) under the European Union's Horizon 2020 research and innovation programme (project EMC2, grant agreement No 810367), and, together with A.I., E.M., by Agence Nationale de la Recherche, under grant ANR-21-CE40-0006 (SINEQ). R.G., G.S. and U.V. also benefited from the Imperial-CNRS 2023 PhD joint program. X.S. and G.S. acknowledge the support of the Royal Society through the International Exchanges Scheme (reference number IES$\backslash$R3$\backslash$203007). X.S. also acknowledges the support of the Royal Society through the Research Grants (reference number RG$\backslash$R2$\backslash$232257), the Isaac Newton Institute for Mathematical Sciences through the Network Support for the Mathematical Sciences initiative (EPSRC reference number EP/V521929/1), and the London Mathematical Society through the Joint Research Groups in the UK scheme (reference number 32430).



\begin{thebibliography}{100}

\bibitem{abdulle2017}
A.~Abdulle, G.~A. Pavliotis, and U.~Vaes.
\newblock Spectral methods for multiscale stochastic differential equations.
\newblock {\em SIAM/ASA J. Uncertain. Quantif.}, 5(1):720--761, 2017.

\bibitem{ackley_1985}
D.~H. Ackley, G.~E. Hinton, and T.~J. Sejnowski.
\newblock A learning algorithm for {B}oltzmann machines.
\newblock {\em Cognitive Science}, 9(1):147--169, 1985.

\bibitem{albergo_2023}
M.~S. Albergo, N.~M. Boffi, and E.~Vanden-Eijnden.
\newblock Stochastic interpolants: {A} unifying framework for flows and
  diffusions.
\newblock {\em Journal of Machine Learning Research}, 26(209):1--80, 2025.

\bibitem{alder1970}
B.~J. Alder, D.~M. Gass, and T.~E. Wainwright.
\newblock Studies in molecular dynamics. {VIII}. {T}he transport coefficients
  for a hard-sphere fluid.
\newblock {\em The Journal of Chemical Physics}, 53(10):3813--3826, 1970.

\bibitem{Allen2017}
M.~P. Allen and D.~J. Tildesley.
\newblock {\em Computer Simulation of Liquids, 2nd ed.}
\newblock Oxford University Press, Oxford, 2017.

\bibitem{andersen_1983}
H.~C. Andersen.
\newblock {RATTLE}: {A} "velocity" version of the {SHAKE} algorithm for
  molecular dynamics calculations.
\newblock {\em Journal of Computational Physics}, 52(1):24--34, 1983.

\bibitem{andradottir_1993}
S.~Andrad\'{o}ttir, D.~P. Heyman, and T.~J. Ott.
\newblock Variance reduction through smoothing and control variates for
  {Markov} chain simulations.
\newblock {\em ACM Trans. Model. Comput. Simul.}, 3(3):167–189, 1993.

\bibitem{AC1999}
R.~Assaraf and M.~Caffarel.
\newblock Zero-variance principle for {M}onte {C}arlo algorithms.
\newblock {\em Phys. Rev. Lett.}, 83:4682--4685, 1999.

\bibitem{bakry-gentil-ledoux-14}
D.~Bakry, I.~Gentil, and M.~Ledoux.
\newblock {\em Analysis and Geometry of {M}arkov Diffusion Operators}.
\newblock Springer, 2014.

\bibitem{Balian}
R.~Balian.
\newblock {\em From Microphysics to Macrophysics. Methods and Applications of
  Statistical Physics}, volume I - II.
\newblock Springer, 2007.

\bibitem{barker_1971}
J.~A. Barker, R.~A. Fisher, and R.~O. Watts.
\newblock Liquid {A}rgon: {Monte Carlo} and molecular dynamics calculations.
\newblock {\em Molecular Physics}, 21(4):657--673, 1971.

\bibitem{Ciccotti_history_MD}
G.~Battimelli, G.~Ciccotti, and P.~Greco.
\newblock {\em Computer Meets Theoretical Physics. The New Frontier of
  Molecular Simulation}.
\newblock The Frontiers Collection. Springer, 2020.

\bibitem{Bennett_1975}
C.~H. Bennett.
\newblock Mass tensor molecular dynamics.
\newblock {\em Journal of Computational Physics}, 19(3):267--279, 1975.

\bibitem{bennett_1976}
C.~H. Bennett.
\newblock Efficient estimation of free energy differences from {Monte Carlo}
  data.
\newblock {\em Journal of Computational Physics}, 22(2):245--268, 1976.

\bibitem{BFLS20}
E.~Bernard, M.~Fathi, A.~Levitt, and G.~Stoltz.
\newblock Hypocoercivity with {S}chur complements.
\newblock {\em Annales Henri Lebesgue}, 5:523--557, 2022.

\bibitem{Bhattacharya_1982}
R.~N. Bhattacharya.
\newblock On the functional central limit theorem and the law of the iterated
  logarithm for {Markov} processes.
\newblock {\em Zeitschrift f{\"u}r Wahrscheinlichkeitstheorie und Verwandte
  Gebiete}, 60(2):185--201, 1982.

\bibitem{github_norton}
N.~Blassel.
\newblock {G}it{Hub} code repository implementing the {N}orton method for
  mobility and shear viscosity computations :
  \url{https://github.com/noeblassel/NortonMethod}, 2024.

\bibitem{github_cecam}
N.~Blassel, L.~Carillo, S.~Darshan, R.~Gastaldello, J.~Greener, A.~Iacobucci,
  and E.~Marini.
\newblock {G}it{H}ub code repository for the {CECAM} summer school:
  ``{S}ampling high-dimensional probability measures with applications in
  (non)equilibrium molecular dynamics and statistics'':
  \url{https://github.com/shiva-darshan/CECAMSummerSchool2025}, 2025.

\bibitem{blassel2024}
N.~Blassel and G.~Stoltz.
\newblock Fixing the flux: {A} dual approach to computing transport
  coefficients.
\newblock {\em Journal of Statistical Physics}, 191:17, 2024.

\bibitem{BLR2000}
F.~Bonetto, J.~L. Lebowitz, and L.~Rey-Bellet.
\newblock Fourier's law: {A} challenge for theorists.
\newblock In A.~Fokas, A.~Grigoryan, T.~Kibble, and B.~Zegarlinsky, editors,
  {\em Mathematical Physics 2000}, pages 128--151. Imperial College Press,
  2000.

\bibitem{BO10}
N.~Bou-Rabee and H.~Owhadi.
\newblock Long-run accuracy of variational integrators in the stochastic
  context.
\newblock {\em SIAM J. Numer. Anal.}, 48:278--297, 2010.

\bibitem{brooks_2011}
S.~Brooks, A.~Gelman, G.~Jones, and X-L. Meng.
\newblock {\em Handbook of Markov Chain Monte Carlo}.
\newblock CRC press, 2011.

\bibitem{Caflisch_1998}
R.~E. Caflisch.
\newblock {Monte Carlo and quasi-Monte Carlo methods}.
\newblock {\em Acta Numerica}, 7:1–49, 1998.

\bibitem{CDK+2003}
E.~Canc\`es, M.~Defranceschi, W.~Kutzelnigg, C.~{Le Bris}, and Y.~Maday.
\newblock Computational quantum chemistry: {A} primer.
\newblock In P.~G. Ciarlet and C.~{Le Bris}, editors, {\em Handbook of
  Numerical Analysis (Special Volume on Computational Chemistry)}, volume~X,
  pages 3--270. Elsevier, 2003.

\bibitem{cancrini2017}
N.~Cancrini and S.~Olla.
\newblock Ensemble dependence of fluctuations: {C}anonical/microcanonical
  equivalence of ensembles.
\newblock {\em Journal of Statistical Physics}, 168:707--730, 2017.

\bibitem{Carmona2007}
P.~Carmona.
\newblock Existence and uniqueness of an invariant measure for a chain of
  oscillators in contact with two heat baths.
\newblock {\em Stochastic Processes and their Applications}, 117:1076--1092,
  2007.

\bibitem{carter_1989}
E.~A. Carter, G.~Ciccotti, J.~T. Hynes, and R.~Kapral.
\newblock Constrained reaction coordinate dynamics for the simulation of rare
  events.
\newblock {\em Chemical Physics Letters}, 156(5):472--477, 1989.

\bibitem{Casas2022}
F.~Casas, J.~M. Sanz-Serna, and L.~Shaw.
\newblock Split {Hamiltonian Monte Carlo} revisited.
\newblock {\em Statistics and Computing}, 32(5):86, 2022.

\bibitem{C_rou_2011}
F.~C\'erou, A.~Guyader, T.~Lelièvre, and D.~Pommier.
\newblock A multiple replica approach to simulate reactive trajectories.
\newblock {\em The Journal of Chemical Physics}, 134(5):054108, 2011.

\bibitem{chen_1997}
M.-H. Chen and Q.-M. Shao.
\newblock {On Monte Carlo methods for estimating ratios of normalizing
  constants}.
\newblock {\em The Annals of Statistics}, 25(4):1563--1594, 1997.

\bibitem{chodera2011dynamical}
J.~D. Chodera, W.~C. Swope, F.~Noe, J.-H. Prinz, M.~R. Shirts, and V.~S. Pande.
\newblock Dynamical reweighting: {I}mproved estimates of dynamical properties
  from simulations at multiple temperatures.
\newblock {\em Journal of Chemical Physics}, 134(24):244107, 2011.

\bibitem{ciccotti2016}
G.~Ciccotti and M.~Ferrario.
\newblock Non-equilibrium by molecular dynamics: {A} dynamical approach.
\newblock {\em Molecular Simulation}, 42(16):1385--1400, 2016.

\bibitem{ciccotti1975}
G.~Ciccotti and G.~Jacucci.
\newblock Direct computation of dynamical response by molecular dynamics: The
  mobility of a charged {L}ennard--{J}ones particle.
\newblock {\em Physical Review Letters}, 35(12):789, 1975.

\bibitem{ciccotti1979}
G.~Ciccotti, G.~Jacucci, and I.~R. McDonald.
\newblock "{T}hought-experiments" by molecular dynamics.
\newblock {\em Journal of Statistical Physics}, 21:1--22, 1979.

\bibitem{ciccotti1976}
G.~Ciccotti, G.~Jacucci, and I.R. McDonald.
\newblock Transport properties of molten alkali halides.
\newblock {\em Physical Review A}, 13(1):426, 1976.

\bibitem{ciccotti2005}
G.~Ciccotti, R.~Kapral, and A.~Sergi.
\newblock Non-equilibrium molecular dynamics.
\newblock In S.~Yip, editor, {\em Handbook of Materials Modeling: Methods},
  pages 745--761. Springer, 2005.

\bibitem{cui_2025}
T.~Cui, X.~Tong, and O.~Zahm.
\newblock Optimal {Riemannian} metric for {Poincar\'e} inequalities and how to
  ideally precondition {Langevin} dynamics.
\newblock {\em arXiv preprint}, 2404.02554, 2025.

\bibitem{darshan2024}
S.~Darshan, A.~Eberle, and G.~Stoltz.
\newblock Sticky coupling as a control variate for sensitivity analysis.
\newblock {\em arXiv preprint}, 2409.15500, 2024.

\bibitem{darshan2025}
S.~Darshan and G.~Stoltz.
\newblock Equivalence of {N}orton and {T}h\'evenin ensembles for mean-field
  interacting particle systems.
\newblock {\em In preparation}, 2025.

\bibitem{davis_1984}
M.~H.~A. Davis.
\newblock Piecewise-deterministic {M}arkov processes: {A} general class of
  non-diffusion stochastic models.
\newblock {\em Journal of the Royal Statistical Society: Series B
  (Methodological)}, 46(3):353--376, 1984.

\bibitem{dSOG17}
L.~de~Sousa~Oliveira and P.~A. Greaney.
\newblock Method to manage integration error in the {Green-Kubo} method.
\newblock {\em Phys. Rev. E}, 95:023308, 2017.

\bibitem{Dhar2008}
A.~Dhar.
\newblock Heat transport in low-dimensional systems.
\newblock {\em Advances in Physics}, 57(5):457--537, 2008.

\bibitem{DG23}
M.~Dobson and A.~K.~A. Geraldo.
\newblock Convergence of nonequilibrium {L}angevin dynamics for planar flows.
\newblock {\em J. Stat. Phys.}, 190:91, 2023.

\bibitem{donati2017girsanov}
L.~Donati, C.~Hartmann, and B.~G. Keller.
\newblock Girsanov reweighting for path ensembles and {M}arkov state models.
\newblock {\em Journal of Chemical Physics}, 146(24):244112, 2017.

\bibitem{donati2018girsanov}
L.~Donati and B.~G. Keller.
\newblock Girsanov reweighting for metadynamics simulations.
\newblock {\em The Journal of Chemical Physics}, 149(7):072335, 2018.

\bibitem{DEMS25}
A.~Durmus, A.~Enfroy, E.~Moulines, and G.~Stoltz.
\newblock Uniform minorization condition and convergence bounds for
  discretizations of kinetic {L}angevin dynamics.
\newblock {\em Ann. Inst. H. Poincaré Probab. Statist.}, 61(1):629--664, 2025.

\bibitem{MR4807807}
M.~Dus and V.~Ehrlacher.
\newblock Numerical solution of {P}oisson partial differential equation in high
  dimension using two-layer neural networks.
\newblock {\em Math. Comp.}, 94(351):159--208, 2025.

\bibitem{espanol_1995}
P.~Español and P.~Warren.
\newblock Statistical mechanics of dissipative particle dynamics.
\newblock {\em Europhysics Letters}, 30(4):191, 1995.

\bibitem{EthKur86}
S.~N. Ethier and T.~G. Kurtz.
\newblock {\em Markov Processes}.
\newblock Probability and Mathematical Statistics. John Wiley \& Sons Inc.,
  1986.

\bibitem{evans1993}
D.~J. Evans.
\newblock The equivalence of {Norton} and {Thévenin} ensembles.
\newblock {\em Molecular Physics}, 80:221--224, 1993.

\bibitem{evans1986}
D.~J. Evans and J.~F. Ely.
\newblock Viscous flow in the stress ensemble.
\newblock {\em Molecular Physics}, 59(5):1043--1048, 1986.

\bibitem{hoover1983}
D.~J. Evans, W.~G. Hoover, B.~H. Failor, B.~Moran, and A.~J.~C. Ladd.
\newblock Nonequilibrium molecular dynamics via {Gauss}'s principle of least
  constraint.
\newblock {\em Physical Review A}, 28(2):1016--1021, 1983.

\bibitem{evans1983}
D.~J. Evans and G.~P. Morriss.
\newblock The isothermal/isobaric molecular dynamics ensemble.
\newblock {\em Physics Letters A}, 98(8):433--436, 1983.

\bibitem{evans1985}
D.~J. Evans and G.~P. Morriss.
\newblock Equilibrium-fluctuation expression for the resistance of a {Norton}
  circuit.
\newblock {\em Physical Review A}, 31(6):3817--3819, 1985.

\bibitem{evans1988}
D.~J. Evans and G.~P. Morriss.
\newblock Transient-time-correlation functions and the rheology of fluids.
\newblock {\em Phys. Rev. A}, 38(8):4142, 1988.

\bibitem{evans_morriss_2007}
D.~J. Evans and G.~P. Morriss.
\newblock {\em Statistical Mechanics of Nonequilibrium Liquids}.
\newblock Cambridge University Press, 2008.

\bibitem{fathi2017improving}
M.~Fathi and G.~Stoltz.
\newblock Improving dynamical properties of metropolized discretizations of
  overdamped {L}angevin dynamics.
\newblock {\em Numerische Mathematik}, 136:545--602, 2017.

\bibitem{FPU1955}
E.~Fermi, J.~Pasta, and S.~M. Ulam.
\newblock Studies of nonlinear problems, {{I}}.
\newblock {\em Los-Alamos Internal Report}, 1955.

\bibitem{fischer_2023}
J.~Fischer and M.~Wendland.
\newblock On the history of key empirical intermolecular potentials.
\newblock {\em Fluid Phase Equilibria}, 573:113876, 2023.

\bibitem{fishman_1996}
G.~S. Fishman.
\newblock {\em Monte Carlo: Concepts, Algorithms, and Applications}.
\newblock Springer New York, NY, 1996.

\bibitem{Frenkel2001}
D.~Frenkel and B.~Smit.
\newblock {\em Understanding Molecular Simulation: From Algorithms to
  Applications, 2nd edition}.
\newblock Academic Press, 2001.

\bibitem{Gallavotti2008}
G.~Gallavotti, editor.
\newblock {\em The {{Fermi-Pasta-Ulam}} Problem: {A} Status Report}.
\newblock Number 728 in Lecture Notes in Physics. Springer, 2008.

\bibitem{gastaldello2025dynamical}
R.~Gastaldello, G.~Stoltz, and U.~Vaes.
\newblock Dynamical reweighting for estimation of fluctuation formulas.
\newblock {\em arXiv preprint}, 2504.11968, 2025.

\bibitem{GS00}
O.~V. Gendelman and A.~V. Savin.
\newblock Normal heat conductivity of the one-dimensional lattice with periodic
  potential of nearest-neighbor interaction.
\newblock {\em Phys. Rev. Lett.}, 84:2381--2384, 2000.

\bibitem{GLPV00}
C.~Giardin\`a, R.~Livi, A.~Politi, and M.~Vassalli.
\newblock Finite thermal conductivity in {1D} lattices.
\newblock {\em Phys. Rev. Lett.}, 84:2144--2147, 2000.

\bibitem{girolami_2011}
M.~Girolami and B.~Calderhead.
\newblock Riemann manifold {Langevin and Hamiltonian Monte Carlo} methods.
\newblock {\em Journal of the Royal Statistical Society: Series B (Statistical
  Methodology)}, 73(2):123--214, 2011.

\bibitem{glynn2019likelihood}
P.~W. Glynn and M.~Olvera-Cravioto.
\newblock Likelihood ratio gradient estimation for steady-state parameters.
\newblock {\em Stochastic Systems}, 9(2):83--100, 2019.

\bibitem{goodfellow_2014}
I.~Goodfellow, J.~Pouget-Abadie, M.~Mirza, B.~Xu, D.~Warde-Farley, S.~Ozair,
  A.~Courville, and Y.~Bengio.
\newblock {Generative {A}dversarial {N}ets}.
\newblock In {\em {Advances in Neural Information Processing Systems}}, pages
  2672--2680, 2014.

\bibitem{goodman_2010}
J.~Goodman and J.~Weare.
\newblock Ensemble samplers with affine invariance.
\newblock {\em Communications in Applied Mathematics and Computational
  Science}, 5(1):65--80, 2010.

\bibitem{gms73}
E.~M. Gosling, I.~R. McDonald, and K.~Singer.
\newblock On the calculation by molecular dynamics of the shear viscosity of a
  simple fluid.
\newblock {\em Molecular Physics}, 26:1475--1484, 1973.

\bibitem{green1954}
M.~S. Green.
\newblock Markoff random processes and the statistical mechanics of
  time-dependent phenomena. {II}. {I}rreversible processes in fluids.
\newblock {\em The Journal of Chemical Physics}, 22(3):398--413, 1954.

\bibitem{Haario1999}
H.~Haario, E.~Saksman, and J.~Tamminen.
\newblock Adaptive proposal distribution for random walk {Metropolis}
  algorithm.
\newblock {\em Computational Statistics}, 14(3):375--395, 1999.

\bibitem{haario_2001}
H.~Haario, E.~Saksman, and J.~Tamminen.
\newblock {An adaptive Metropolis algorithm}.
\newblock {\em Bernoulli}, 7(2):223--242, 2001.

\bibitem{hairer_2006}
E.~Hairer, C.~Lubich, and G.~Wanner.
\newblock {\em {Geometric Numerical Integration}}, volume~31 of {\em Springer
  Series in Computational Mathematics}.
\newblock Springer-Verlag, Berlin, second edition, 2006.

\bibitem{hairer_2010}
M.~Hairer and A.~J. Majda.
\newblock A simple framework to justify linear response theory.
\newblock {\em Nonlinearity}, 23(4):909, 2010.

\bibitem{hairer_2011}
M.~Hairer and J.~C. Mattingly.
\newblock Yet another look at {H}arris' ergodic theorem for {M}arkov chains.
\newblock In {\em Seminar on {S}tochastic {A}nalysis, {R}andom {F}ields and
  {A}pplications {VI}}, volume~63 of {\em Progr. Probab.}, pages 109--117.
  Birkh\"auser/Springer, 2011.

\bibitem{hammersley_1964}
J.~M. Hammersley and D.~C. Handscomb.
\newblock {\em Monte Carlo Methods}.
\newblock Springer Dordrecht, 1964.

\bibitem{hastings_1970}
W.~K. Hastings.
\newblock Monte {Carlo} sampling methods using {Markov} chains and their
  applications.
\newblock {\em Biometrika}, 57(1):97--109, 1970.

\bibitem{Henderson1997}
S.~G. Henderson.
\newblock {\em Approximating martingales for variance reduction in {M}arkov
  process simulation}.
\newblock PhD thesis. Stanford University, 1997.

\bibitem{henderson_2002}
S.~G. Henderson and P.~W. Glynn.
\newblock Approximating martingales for variance reduction in {Markov} process
  simulation.
\newblock {\em Mathematics of Operations Research}, 27(2):253--271, 2002.

\bibitem{hoogerbrugge_1992}
P.~J. Hoogerbrugge and J.~M. V.~A. Koelman.
\newblock Simulating microscopic hydrodynamic phenomena with dissipative
  particle dynamics.
\newblock {\em Europhysics Letters}, 19(3):155, 1992.

\bibitem{hoover1986}
W.~G. Hoover.
\newblock Nonequilibrium molecular dynamics.
\newblock In {\em Molecular Dynamics}, volume 258 of {\em Lecture Notes in
  Physics}, pages 92--131. Springer, 1986.

\bibitem{hoover1993}
W.~G. Hoover.
\newblock Nonequilibrium molecular dynamics: {T}he first 25 years.
\newblock {\em Physica A: Statistical Mechanics and its Applications},
  194(1-4):450--461, 1993.

\bibitem{hormander_67}
L.~H\"ormander.
\newblock Hypoelliptic second order differential equations.
\newblock {\em Acta Mathematica}, 119:147--171, 1967.

\bibitem{IOS2021}
A.~Iacobucci, S.~Olla, and G.~Stoltz.
\newblock Thermo-mechanical transport in rotor chains.
\newblock {\em J. Stat. Phys.}, 183(2):26, 2021.

\bibitem{ikeda_1981}
N.~Ikeda and S.~Watanabe.
\newblock {Chapter IV - Stochastic Differential Equations}.
\newblock In {\em {Stochastic Differential Equations and Diffusion Processes}},
  volume~24 of {\em North-Holland Mathematical Library}, pages 145--232.
  Elsevier, 1981.

\bibitem{irving_kirkwood_1950}
J.~H. Irving and J.~G. Kirkwood.
\newblock The statistical mechanical theory of transport processes. {IV}. {T}he
  equations of hydrodynamics.
\newblock {\em The Journal of Chemical Physics}, 18(6):817--829, 1950.

\bibitem{JMA+2025}
R.~Jacobs, D.~Morgan, S.~Attarian, J.~Meng, C.~Shen, Z.~Wu, C.~Y. Xie, J.~H.
  Yang, N.~Artrith, B.~Blaiszik, G.~Ceder, K.~Choudhary, G.~Csanyi, E.~D.
  Cubuk, B.~Deng, R.~Drautz, X.~Fu, J.~Godwin, V.~Honavar, O.~Isayev,
  A.~Johansson, B.~Kozinsky, S.~Martiniani, S.~P. Ong, I.~Poltavsky, K.~J.
  Schmidt, S.~Takamoto, A.~P. Thompson, J.~Westermayr, and B.~M. Wood.
\newblock A practical guide to machine learning interatomic potentials –
  {{Status}} and future.
\newblock {\em Current Opinion in Solid State and Materials Science},
  35:101214, 2025.

\bibitem{jardat_1999}
M.~Jardat, O.~Bernard, P.~Turq, and G.~R. Kneller.
\newblock Transport coefficients of electrolyte solutions from {Smart}
  {Brownian} dynamics simulations.
\newblock {\em The Journal of Chemical Physics}, 110(16):7993--7999, 1999.

\bibitem{js12}
R.~Joubaud and G.~Stoltz.
\newblock Nonequilibrium {shear} {viscosity} {computations} with {Langevin}
  {dynamics}.
\newblock {\em Multiscale Modeling \& Simulation}, 10:191--216, 2012.

\bibitem{khasminskii_2012_i}
R.~Khasminskii.
\newblock Markov processes and stochastic differential equations.
\newblock In {\em {Stochastic Stability of Differential Equations}}, pages
  59--98. Springer Berlin Heidelberg, Berlin, Heidelberg, 2012.

\bibitem{kingma_2013}
D.~P. Kingma and M.~Welling.
\newblock {Auto-Encoding Variational Bayes}.
\newblock In {\em 2nd International Conference on Learning Representations,
  {ICLR} 2014, Banff, AB, Canada}, 2014.

\bibitem{kirkwood_1935}
J.~G. Kirkwood.
\newblock Statistical mechanics of fluid mixtures.
\newblock {\em The Journal of Chemical Physics}, 3(5):300--313, 1935.

\bibitem{kliemann_1987}
W.~Kliemann.
\newblock Recurrence and invariant measures for degenerate diffusions.
\newblock {\em The Annals of Probability}, 15(2):690--707, 1987.

\bibitem{kobyzev_2021}
I.~Kobyzev, S.~D. Prince, and M.~A. Brubaker.
\newblock Normalizing {F}lows: {A}n introduction and review of current methods.
\newblock {\em IEEE Transactions on Pattern Analysis and Machine Intelligence},
  43(11):3964--3979, 2021.

\bibitem{KomorowskiLandimOlla2012}
T.~Komorowski, C.~Landim, and S.~Olla.
\newblock {\em Fluctuations in {M}arkov processes: Time symmetry and martingale
  approximation}, volume 345 of {\em Grundlehren der Mathematischen
  Wissenschaften [Fundamental Principles of Mathematical Sciences]}.
\newblock Springer, 2012.

\bibitem{kopec1}
M.~Kopec.
\newblock Weak backward error analysis for {L}angevin process.
\newblock {\em BIT Numerical Mathematics}, 55(4):1057–1103, 2015.

\bibitem{kopec2}
M.~Kopec.
\newblock Weak backward error analysis for overdamped {L}angevin processes.
\newblock {\em IMA J. Numer. Anal.}, 35(2):583--614, 2015.

\bibitem{kubo1957}
R.~Kubo.
\newblock Statistical-mechanical theory of irreversible processes. {I}.
  {G}eneral theory and simple applications to magnetic and conduction problems.
\newblock {\em Journal of the Physical Society of Japan}, 12(6):570--586, 1957.

\bibitem{kubo1957b}
R.~Kubo, M.~Yokota, and S.~Nakajima.
\newblock Statistical-mechanical theory of irreversible processes. {II}.
  {R}esponse to thermal disturbance.
\newblock {\em Journal of the Physical Society of Japan}, 12(11):1203--1211,
  1957.

\bibitem{KDN09}
A.~Kundu, A.~Dhar, and O.~Narayan.
\newblock The {G}reen-{K}ubo formula for heat conduction in open systems.
\newblock {\em J. Stat. Mech. - Theory E.}, page L03001, 2009.

\bibitem{lapeyre_2003}
B.~Lapeyre, E.~Pardoux, and R.~Sentis.
\newblock {\em Introduction to Monte-Carlo Methods for Transport and Diffusion
  Equations}.
\newblock Oxford University Press, 2003.

\bibitem{lebowitz1967}
J.~L. Lebowitz, J.~K. Percus, and L.~Verlet.
\newblock Ensemble dependence of fluctuations with application to machine
  computations.
\newblock {\em Physical Review}, 153(1):250, 1967.

\bibitem{le72}
A.~W. Lees and S.~F. Edwards.
\newblock The computer study of transport processes under extreme conditions.
\newblock {\em Journal of Physics C: Solid State Physics}, 5(15):1921, 1972.

\bibitem{leimkuhler_2015}
B.~Leimkuhler and C.~Matthews.
\newblock {\em Molecular Dynamics}.
\newblock Springer Cham, 2015.

\bibitem{leimkuhler2016computation}
B.~Leimkuhler, C.~Matthews, and G.~Stoltz.
\newblock The computation of averages from equilibrium and nonequilibrium
  {L}angevin molecular dynamics.
\newblock {\em IMA Journal of Numerical Analysis}, 36(1):13--79, 2016.

\bibitem{Leimkuhler2018}
B.~Leimkuhler, C.~Matthews, and J.~Weare.
\newblock Ensemble preconditioning for {Markov chain Monte Carlo} simulation.
\newblock {\em Statistics and Computing}, 28(2):277--290, 2018.

\bibitem{lelievre_2010}
T.~Leli{\`e}vre, M.~Rousset, and G.~Stoltz.
\newblock {\em {Free Energy Computations: A Mathematical Perspective}}.
\newblock {Imperial College Press}, 2010.

\bibitem{lelievre2024}
T.~Leli{\`e}vre, R.~Santet, and G.~Stoltz.
\newblock Unbiasing {Hamiltonian Monte Carlo} algorithms for a general
  {Hamiltonian} function.
\newblock {\em Foundations of Computational Mathematics}, 26:1--74, 2026.

\bibitem{lelievre2016}
T.~Leli\`evre and G.~Stoltz.
\newblock Partial differential equations and stochastic methods in molecular
  dynamics.
\newblock {\em Acta Numerica}, 25:681--880, 2016.

\bibitem{lelievre_2025_ii}
T.~Lelièvre, G.~A. Pavliotis, G.~Robin, R.~Santet, and G.~Stoltz.
\newblock Optimizing the diffusion coefficient of overdamped {L}angevin
  dynamics.
\newblock {\em Mathematics of Computation}, 2025.
\newblock to appear.

\bibitem{lelievre_2012}
T.~Lelièvre, M.~Rousset, and G.~Stoltz.
\newblock Langevin dynamics with constraints and computation of free energy
  differences.
\newblock {\em Mathematics of Computation}, 81:2071--2125, 2012.

\bibitem{lelievre_2019}
T.~Lelièvre, M.~Rousset, and G.~Stoltz.
\newblock Hybrid {{Monte Carlo}} methods for sampling probability measures on
  submanifolds.
\newblock {\em Numerische Mathematik}, 143(2):379--421, 2019.

\bibitem{lelievre_2025}
T.~Lelièvre, R.~Santet, and G.~Stoltz.
\newblock Improving sampling by modifying the effective diffusion.
\newblock {\em Journal of Computational Physics}, 541:114313, 2025.

\bibitem{lennard-jones_1924_i}
J.~E. Lennard-Jones and S.~Chapman.
\newblock On the determination of molecular fields. {I}. {F}rom the variation
  of the viscosity of a gas with temperature.
\newblock {\em Proceedings of the Royal Society of London. Series A},
  106(738):441--462, 1924.

\bibitem{lennard-jones_1924_ii}
J.~E. Lennard-Jones and S.~Chapman.
\newblock On the determination of molecular fields. {II}. {F}rom the equation
  of state of a gas.
\newblock {\em Proceedings of the Royal Society of London. Series A},
  106(738):463--477, 1924.

\bibitem{LLP2003}
S.~Lepri, R.~Livi, and A.~Politi.
\newblock Thermal conduction in classical low-dimensional lattices.
\newblock {\em Physics Reports}, 377(1):1--80, 2003.

\bibitem{LLP2016}
S.~Lepri, R.~Livi, and A.~Politi.
\newblock Heat transport in low dimensions: {{Introduction}} and phenomenology.
\newblock In S.~Lepri, editor, {\em Thermal {{Transport}} in {{Low
  Dimensions}}: {{From Statistical Physics}} to {{Nanoscale Heat Transfer}}},
  pages 1--37. Springer International Publishing, 2016.

\bibitem{verlet1970}
D.~Levesque and L.~Verlet.
\newblock Computer ``experiments'' on classical fluids. {III}. {T}ime-dependent
  self-correlation functions.
\newblock {\em Physical Review A}, 2(6):2514, 1970.

\bibitem{verlet1973}
D.~Levesque, L.~Verlet, and J.~K{\"u}rkijarvi.
\newblock Computer ``experiments'' on classical fluids. {IV}. {T}ransport
  properties and time-correlation functions of the {L}ennard--{J}ones liquid
  near its triple point.
\newblock {\em Physical Review A}, 7(5):1690, 1973.

\bibitem{li_2022}
R.~Li, M.~Tao, S.~S. Vempala, and A.~Wibisono.
\newblock The mirror {Langevin} algorithm converges with vanishing bias.
\newblock In S.~Dasgupta and N.~Haghtalab, editors, {\em Proceedings of The
  33rd International Conference on Algorithmic Learning Theory}, volume 167 of
  {\em Proceedings of Machine Learning Research}, pages 718--742. PMLR, 2022.

\bibitem{liu_2008}
J.~S. Liu.
\newblock {\em Monte Carlo Strategies in Scientific Computing}.
\newblock Springer Publishing Company, Incorporated, 2008.

\bibitem{lu2019}
Y.~Lu and J.~Y. Park.
\newblock Estimation of longrun variance of continuous time stochastic process
  using discrete sample.
\newblock {\em Journal of Econometrics}, 210(2):236--267, 2019.

\bibitem{MESVDDT24}
L.~Maffioli, J.~P. Ewen, E.~R. Smith, S.~Varghese, P.~J. Daivis, D.~Dini, and
  B.~D. Todd.
\newblock {TTCF4LAMMPS}: {A} toolkit for simulation of the non-equilibrium
  behaviour of molecular fluids at experimentally accessible shear rates.
\newblock {\em Computer Physics Communications}, 300:109205, 2024.

\bibitem{MSEDDT22}
L.~Maffioli, E.~R. Smith, J.~P. Ewen, P.~J. Daivis, D.~Dini, and B.~D. Todd.
\newblock Slip and stress from low shear rate nonequilibrium molecular
  dynamics: {T}he transient-time correlation function technique.
\newblock {\em The Journal of Chemical Physics}, 156(18):184111, 2022.

\bibitem{MSH02}
J.~C. Mattingly, A.~M. Stuart, and D.~J. Higham.
\newblock Ergodicity for {SDEs} and approximations: locally {Lipschitz} vector
  fields and degenerate noise.
\newblock {\em Stoch. Proc. Appl.}, 101(2):185--232, 2002.

\bibitem{metropolis_1953}
N.~Metropolis, A.~W. Rosenbluth, M.~N. Rosenbluth, A.~H. Teller, and E.~Teller.
\newblock Equation of state calculations by fast computing machines.
\newblock {\em The Journal of Chemical Physics}, 21(6):1087--1092, 1953.

\bibitem{Metzner_2006}
P.~Metzner, C.~Schütte, and E.~Vanden-Eijnden.
\newblock Illustration of transition path theory on a collection of simple
  examples.
\newblock {\em J. Chem. Phys.}, 125(08):084110, 2006.

\bibitem{michel_2014}
M.~Michel, S.~C. Kapfer, and W.~Krauth.
\newblock {Generalized event-chain {Monte Carlo}: {C}onstructing rejection-free
  global-balance algorithms from infinitesimal steps}.
\newblock {\em The Journal of Chemical Physics}, 140(5):054116, 2014.

\bibitem{MSI2013}
A.~Mira, R.~Solgi, and D.~Imparato.
\newblock Zero variance {M}arkov {C}hain {M}onte {C}arlo for {B}ayesian
  estimators.
\newblock {\em Statistics and Computing}, 23:653–662, 2013.

\bibitem{MSS24}
P.~Monmarch\'e, R.~Spacek, and G.~Stoltz.
\newblock Transient subtraction: {A} control variate method for computing
  transport coefficients.
\newblock {\em J. Stat. Phys.}, 192:53, 2025.

\bibitem{morriss1987}
G.~P. Morriss and D.~J. Evans.
\newblock Application of transient correlation functions to shear flow far from
  equilibrium.
\newblock {\em Phys. Rev. A}, 35(2):792, 1987.

\bibitem{noe_2019}
F.~Noé, S.~Olsson, J.~Köhler, and H.~Wu.
\newblock Boltzmann generators: {S}ampling equilibrium states of many-body
  systems with deep learning.
\newblock {\em Science}, 365(6457):eaaw1147, 2019.

\bibitem{OP11}
M.~Ottobre and G.~A. Pavliotis.
\newblock Asymptotic analysis for the generalized {L}angevin equation.
\newblock {\em Nonlinearity}, 24(5):1629--1653, 2011.

\bibitem{papamakarios_2021}
G.~Papamakarios, E.~Nalisnick, D.~J. Rezende, S.~Mohamed, and
  B.~Lakshminarayanan.
\newblock Normalizing {F}lows for probabilistic modeling and inference.
\newblock {\em J. Mach. Learn. Res.}, 22(1):57, 2021.

\bibitem{Park_2003}
S.~Park, M.~K. Sener, D.~Lu, and K.~Schulten.
\newblock Reaction paths based on mean first-passage times.
\newblock {\em The Journal of Chemical Physics}, 119(3):1313–1319, 2003.

\bibitem{patterson_2013}
S.~Patterson and Y.~W. Teh.
\newblock Stochastic gradient {Riemannian Langevin} dynamics on the probability
  simplex.
\newblock In C.~J. Burges, L.~Bottou, M.~Welling, Z.~Ghahramani, and K.~Q.
  Weinberger, editors, {\em Advances in Neural Information Processing Systems},
  volume~26. Curran Associates, Inc., 2013.

\bibitem{PSSV24}
G.~A. Pavliotis, R.~Spacek, G.~Stoltz, and U.~Vaes.
\newblock Neural network approaches for variance reduction in fluctuation
  formulas.
\newblock {\em SIAM/ASA Journal of Uncertainty Quantification}, 14(1):221--255,
  2026.

\bibitem{pavliotis2008multiscale}
G.~A. Pavliotis and A.~Stuart.
\newblock {\em Multiscale {M}ethods: {A}veraging and {H}omogenization},
  volume~53 of {\em Texts in Applied Mathematics}.
\newblock Springer Science \& Business Media, 2008.

\bibitem{pavliotis2023}
Grigorios~A. Pavliotis, Gabriel Stoltz, and Urbain Vaes.
\newblock Mobility estimation for {L}angevin dynamics using control variates.
\newblock {\em Multiscale Model. Simul.}, 21(2):680--715, 2023.

\bibitem{plechac_nemd}
P.~Plechac, G.~Stoltz, and T.~Wang.
\newblock Convergence of the likelihood ratio method for linear response of
  non-equilibrium stationary states.
\newblock {\em ESAIM:M2AN}, 55:S593--S623, 2021.

\bibitem{plechac2022}
P.~Plechac, G.~Stoltz, and T.~Wang.
\newblock Martingale product estimators for sensitivity analysis in
  computational statistical physics.
\newblock {\em IMA Journal of Numerical Analysis}, 43(6):3430--3477, 2022.

\bibitem{rahman_1964}
A.~Rahman.
\newblock Correlations in the motion of atoms in liquid {A}rgon.
\newblock {\em Phys. Rev.}, 136:A405--A411, 1964.

\bibitem{MR3881695}
M.~Raissi, P.~Perdikaris, and G.~E. Karniadakis.
\newblock Physics-informed neural networks: a deep learning framework for
  solving forward and inverse problems involving nonlinear partial differential
  equations.
\newblock {\em J. Comput. Phys.}, 378:686--707, 2019.

\bibitem{bellet_2006}
L.~Rey-Bellet.
\newblock {Ergodic Properties of {M}arkov Processes}.
\newblock In S.~Attal, A.~Joye, and C.-A. Pillet, editors, {\em {Open Quantum
  Systems II}}, volume 1881 of {\em Lecture Notes in Mathematics}, pages 1--39.
  Springer, 2006.

\bibitem{robert_2004}
C.~P. Robert and G.~Casella.
\newblock {\em Monte {Carlo} Statistical Methods}.
\newblock Springer Verlag, 2004.

\bibitem{Roberts_2009}
G.~O. Roberts and J.~S. Rosenthal.
\newblock Examples of adaptive {MCMC}.
\newblock {\em Journal of Computational and Graphical Statistics},
  18(2):349--367, 2009.

\bibitem{rodenhausen1989}
H.~Rodenhausen.
\newblock Einstein's relation between diffusion constant and mobility for a
  diffusion model.
\newblock {\em Journal of Statistical Physics}, 55:1065--1088, 1989.

\bibitem{RS19}
J.~Roussel and G.~Stoltz.
\newblock A perturbative approach to control variates in molecular dynamics.
\newblock {\em Multiscale Model. Simul.}, 17(1):552--591, 2019.

\bibitem{rowley_1975}
L.~A. Rowley, D.~Nicholson, and N.~G. Parsonage.
\newblock {Monte Carlo} grand canonical ensemble calculation in a gas-liquid
  transition region for 12-6 {A}rgon.
\newblock {\em Journal of Computational Physics}, 17(4):401--414, 1975.

\bibitem{rubinstein_2016}
R.~Y. Rubinstein and D.~P. Kroese.
\newblock {\em Simulation and the Monte Carlo Method}.
\newblock John Wiley \& Sons, 2016.

\bibitem{ryckaert_1977}
J-P. Ryckaert, G.~Ciccotti, and H.~J.~C. Berendsen.
\newblock Numerical integration of the {C}artesian equations of motion of a
  system with constraints: {M}olecular dynamics of n-alkanes.
\newblock {\em Journal of Computational Physics}, 23(3):327--341, 1977.

\bibitem{sasaki2025}
R.~Sasaki, Y.~Tateyama, and D.J. Searles.
\newblock Constant-current nonequilibrium molecular dynamics approach for
  accelerated computation of ionic conductivity including ion-ion correlation.
\newblock {\em PRX Energy}, 4(1):013005, 2025.

\bibitem{shirts_2008}
M.~R. Shirts and J.~D. Chodera.
\newblock Statistically optimal analysis of samples from multiple equilibrium
  states.
\newblock {\em The Journal of Chemical Physics}, 129(12):124105, 2008.

\bibitem{MR3874585}
J.~Sirignano and K.~Spiliopoulos.
\newblock D{GM}: a deep learning algorithm for solving partial differential
  equations.
\newblock {\em J. Comput. Phys.}, 375:1339--1364, 2018.

\bibitem{sohl-dickstein_2015}
J.~Sohl-Dickstein, E.~Weiss, N.~Maheswaranathan, and S.~Ganguli.
\newblock Deep unsupervised learning using nonequilibrium thermodynamics.
\newblock In F.~Bach and D.~Blei, editors, {\em Proceedings of the 32nd
  International Conference on Machine Learning}, volume~37 of {\em Proceedings
  of Machine Learning Research}, pages 2256--2265, Lille, France, 2015. PMLR.

\bibitem{spacek2023}
R.~Spacek and G.~Stoltz.
\newblock Extending the regime of linear response with synthetic forcings.
\newblock {\em Multiscale Modeling \& Simulation}, 21(4):1602--1643, 2023.

\bibitem{stephan_2020}
S.~Stephan, J.~Staubach, and H.~Hasse.
\newblock Review and comparison of equations of state for the {Lennard-Jones}
  fluid.
\newblock {\em Fluid Phase Equilibria}, 523:112772, 2020.

\bibitem{StoltzM2}
G.~Stoltz.
\newblock {\em An Introduction to Computational Statistical Physics}.
\newblock M2 lecture notes.

\bibitem{stoltz2024}
G.~Stoltz.
\newblock Error estimates and variance reduction for nonequilibrium stochastic
  dynamics.
\newblock In A.~Hinrichs, P.~Kritzer, and F.~Pillichshammer, editors, {\em
  Monte Carlo and Quasi-Monte Carlo Methods (MCQMC 2022)}, volume 460 of {\em
  Springer Proceedings in Mathematics \& Statistics}, pages 163--187, 2024.

\bibitem{tabak_2010}
E.~G. Tabak and E.~Vanden-Eijnden.
\newblock {Density estimation by dual ascent of the log-likelihood}.
\newblock {\em Communications in Mathematical Sciences}, 8(1):217--233, 2010.

\bibitem{TT90}
D.~Talay and L.~Tubaro.
\newblock Expansion of the global error for numerical schemes solving
  stochastic differential equations.
\newblock {\em Stoch. Proc. Appl.}, 8(4):483--509, 1990.

\bibitem{tee2023}
S.~R. Tee and D.~J. Searles.
\newblock Constant potential and constrained charge ensembles for simulations
  of conductive electrodes.
\newblock {\em Journal of Chemical Theory and Computation}, 19(10):2758--2768,
  2023.

\bibitem{todd_daivis_2017}
B.~D. Todd and P.~J. Daivis.
\newblock {\em Nonequilibrium Molecular Dynamics: Theory, Algorithms and
  Applications}.
\newblock Cambridge University Press, 2017.

\bibitem{tokdar_2010}
S.~T. Tokdar and R.~E. Kass.
\newblock Importance sampling: {A} review.
\newblock {\em WIREs Computational Statistics}, 2(1):54--60, 2010.

\bibitem{Tran2024}
J.~H. Tran and T.~S. Kleppe.
\newblock Tuning diagonal scale matrices for {HMC}.
\newblock {\em Statistics and Computing}, 34(6):196, 2024.

\bibitem{Tuckerman2023}
M.~Tuckerman.
\newblock {\em Statistical Mechanics: Theory and Molecular Simulation}.
\newblock Oxford University Press, 2010.

\bibitem{Villani09}
C.~Villani.
\newblock Hypocoercivity.
\newblock {\em Mem. Amer. Math. Soc.}, 202(950), 2009.

\bibitem{wang2012}
Q.~Wang and N.~Wu.
\newblock Long-run covariance and its applications in cointegration regression.
\newblock {\em The Stata Journal}, 12(3):515--542, 2012.

\bibitem{wang2019steady}
T.~Wang and P.~Plechac.
\newblock Steady-state sensitivity analysis of continuous time {M}arkov chains.
\newblock {\em SIAM Journal on Numerical Analysis}, 57(1):192--217, 2019.

\bibitem{wirnsberger_2020}
P.~Wirnsberger, A.~J. Ballard, G.~Papamakarios, S.~Abercrombie,
  S.~Racani{\`e}re, A.~Pritzel, D.~Jimenez~Rezende, and C.~Blundell.
\newblock Targeted free energy estimation via learned mappings.
\newblock {\em J. Chem. Phys.}, 153(14):144112, 2020.

\bibitem{wu2025}
X.~Wu and X.~Shang.
\newblock Stochastic {N}orton dynamics: {A}n alternative approach for the
  computation of transport coefficients in dissipative particle dynamics.
\newblock {\em Journal of Computational Physics}, 541:114316, 2025.

\bibitem{zhang_ibister_evans}
F.~Zhang, D.~J. Isbister, and D.~J. Evans.
\newblock Nonequilibrium molecular dynamics simulations of heat flow in
  one-dimensional lattices.
\newblock {\em Phys. Rev. E}, 61:3541, 2000.

\end{thebibliography}
\end{document}